\documentclass[preprint,12pt]{elsarticle}
\usepackage[a4paper,margin=2.5cm]{geometry}
\usepackage[utf8]{inputenc}
\usepackage{amsmath,amssymb,amsthm}
\usepackage{booktabs}
\usepackage{graphicx}
\usepackage{xcolor}
\usepackage[section]{placeins}

\usepackage{tikz}
\usetikzlibrary{decorations.pathreplacing}
\usepackage[colorlinks=true,linkcolor=black,citecolor=black,urlcolor=black]{hyperref}
\AtBeginDocument{\hypersetup{linkcolor=black,citecolor=black,urlcolor=black}}

\journal{Journal of Computational Physics}
\makeatletter
\def\ps@pprintTitle{\let\@oddhead\@empty\let\@evenhead\@empty
  \def\@oddfoot{\hfil\thepage\hfil}\let\@evenfoot\@oddfoot}
\makeatother

\newtheorem{theorem}{Theorem}
\newtheorem{proposition}[theorem]{Proposition}
\newtheorem{lemma}[theorem]{Lemma}

\theoremstyle{remark}
\newtheorem{remark}[theorem]{Remark}

\newcommand{\dx}{\partial_x}
\newcommand{\dt}{\partial_t}
\newcommand{\ub}{\bar u}
\newcommand{\hs}{\hat h}
\newcommand{\qs}{\hat q}
\newcommand{\Gb}{\bar\Gamma}
\newcommand{\aph}{{\alpha+1/2}}
\newcommand{\amh}{{\alpha-1/2}}

\begin{document}

\begin{frontmatter}

\title{Hyperbolic relaxation of the multilayer non-hydrostatic system with layer-independent wave
speeds, and an explicit staggered scheme with energy-preserving structure}

\author[inf]{Cipriano Escalante\corref{cor}}
\ead{escalante@uma.es}
\author[cie]{Manuel J. Castro}
\author[cie]{Tomás Morales de Luna}
\cortext[cor]{Corresponding author.}
\affiliation[inf]{organization={Departamento de Matemática Aplicada, Escuela Técnica Superior
de Ingeniería Informática, Universidad de Málaga},
            postcode={29080},
            city={Málaga},
            country={Spain}}
\affiliation[cie]{organization={Departamento de Análisis Matemático, Estadística e
Investigación Operativa y Matemática Aplicada, Facultad de Ciencias, Universidad de Málaga},
            postcode={29080},
            city={Málaga},
            country={Spain}}

\begin{abstract}
We design an explicit scheme for multilayer non-hydrostatic shallow water models that avoids the
global elliptic problem of projection methods. To this end, we propose a hyperbolic relaxation of
the layer-averaged system LDNH$_0$ of Fern\'andez-Nieto, Parisot, Penel and Sainte-Marie
(Commun.\ Math.\ Sci., 2018). If the constraint of each layer is relaxed through the pressure
jump across the layer, the wave speeds, and therefore the cost of an explicit scheme, grow
linearly with the number of layers. We avoid this by writing the constraints in an equivalent
form centered at the interfaces. The resulting relaxed system has an exact energy balance, and
its wave speeds at rest are bounded independently of the number of layers. It is discretized on
a staggered mesh: the hydrostatic part by a path-conservative finite volume method, and the
non-hydrostatic part by a centered scheme that only requires one tridiagonal system per water
column. We show that the discretization in space of this part conserves the energy exactly, and that the time step does it up to a term of order $\tau^2$, which gives its stability condition.
The influence of the relaxation celerity is studied by means of exact solitary waves, and the
scheme is validated against exact solutions, wet-dry tests and laboratory data. In two
dimensions, the GPU implementation is up to $50$ times faster than the projection method.
\end{abstract}


\end{frontmatter}

\section{Introduction}
We consider free-surface water waves in which frequency dispersion plays a relevant role, that
is, waves whose celerity depends on the wavelength. As is well known, the nonlinear shallow
water equations, which are the standard tool for the simulation of tsunamis, storm surges and
coastal flows, do not take into account dispersive effects. Indeed, according to the linear Airy
theory, the phase velocity of a wave of wave number $k$ over a still-water depth $H$ is given by
\[
C_{\mathrm{Airy}}^2=gH\,\frac{\tanh(kH)}{kH},
\]
whereas the phase velocity of the shallow water equations is $C_{\mathrm{SWE}}^2=gH$ for every
wavelength. As a consequence, the shallow water equations overestimate the celerity of the
shorter components of the wave field, and they are not able to represent the dispersive wave
trains that develop, for instance, when waves propagate over submerged obstacles or shoal on a
beach.

Several families of models have been derived in the shallow water framework in order to improve
the dispersive properties. Boussinesq-type models
\cite{boussinesq1872theorie,Peregrine1967,MadsenSorensen1992}, which include the
Serre--Green--Naghdi equations \cite{Serre1953,GreenNaghdi1976} (see the review
\cite{Brocchini2013}), are derived by an asymptotic
expansion of the velocity field and contain high-order derivatives, so that an elliptic problem
has to be solved at each time step
\cite{LannesMarche2015}. Non-hydrostatic pressure models
\cite{Casulli1999,StellingZijlema2003,Bristeau2011} avoid the
hydrostatic assumption and only contain first-order derivatives. Multilayer non-hydrostatic
models go one step further: the water column is divided into several layers, each one with its
own velocities and pressures
\cite{BaiCheung2018,FNPPSM2018,Escalante2019JSC}.
With a few layers, they approximate the dispersion relation and the vertical profiles of the
Euler equations much better than single-layer models, at a much lower cost than a
three-dimensional solver \cite{Escalante2021JSC}.
Models with exact linear dispersion \cite{Zakharov1968,CraigSulem1993}
involve non-local operators, which are difficult to use with general bathymetries and wet-dry
fronts.

In non-hydrostatic models, the non-hydrostatic pressure plays the role of the Lagrange multiplier
associated with an incompressibility constraint. For this reason, these models are usually solved
by projection methods, which go back to the works by Chorin \cite{Chorin1968} and Temam
\cite{Temam1969} for the incompressible Navier--Stokes equations (see also
\cite{GuermondMinevShen2006}), and which have been
widely used for non-hydrostatic shallow flows
\cite{Casulli1999,StellingZijlema2003,Escalante2018AMC,Escalante2019JSC}. At each time step, an
elliptic problem that couples all the cells (and all the layers for multilayer models) has to be solved. This problem is
the most expensive part of the computational cost. Moreover, since it is usually solved by iterative methods
whose number of iterations grows when the mesh is refined, it is not easy to
solve efficiently on massively parallel architectures such as GPUs. This is a serious drawback,
since GPUs have proved to be extremely efficient, for instance, in explicit finite volume schemes for
shallow water equations \cite{Castro2011CRM}.

An alternative that avoids the elliptic problem is given by hyperbolic relaxation, in the spirit
of the classical works by Cattaneo \cite{Cattaneo1958}, the
constraint is replaced by an evolution equation for the pressure, which depends on a relaxation
celerity $c$, and the resulting hyperbolic system can be solved by explicit schemes. After the
first-order hyperbolic approach to dispersion of Mazaheri et al.\ \cite{Mazaheri2016} and the
approach of Favrie and Gavrilyuk \cite{FavrieGavrilyuk2017}, this idea has been applied to dispersive equations, such as the Serre--Green--Naghdi
equations and single-layer non-hydrostatic models, in several works
\cite{Escalante2019JCP,EscalanteMorales2020,Bassi2020,Busto2021,Guermond2022,Tkachenko2023,DuranRichard2026}.

For single-layer models, this approach is by now well established. For a relaxed model with improved
dispersive properties, an energy-stable splitting of the hyperbolic and dispersive parts has been
proposed in \cite{DuranRichard2026}, in which the energy is dissipated at each step. In our scheme,
the non-hydrostatic step does not dissipate energy. Energy-conserving discretizations of the
Serre--Green--Naghdi equations, also in their hyperbolic form, have been proposed in
\cite{RanochaRicchiuto2025}. They use summation-by-parts operators in space and relaxation
Runge--Kutta methods \cite{Ketcheson2019,RanochaKetcheson2020} in time, for one layer and periodic
boundary conditions. For several layers, however, the
choice of the relaxed variables is not obvious, and it has a strong influence on the efficiency
of the resulting explicit schemes. A natural choice is to relax the constraint of each layer
through the pressure jump across the layer. Unfortunately, the interface pressures are then
non-local functions of the relaxed variables, and the wave speeds of the relaxed system grow
linearly with the number of layers, and so does the cost of an explicit scheme (see
Remark~\ref{rem:jumps}). In other words, the advantage of the explicit approach is lost precisely
when the number of layers is increased in order to improve the dispersive properties of the
model.

To the best of our knowledge, the relaxation of multilayer non-hydrostatic models such as LDNH$_0$ has not been considered before.
In this work, we focus on the multilayer non-hydrostatic system LDNH$_0$ introduced in \cite{FNPPSM2018}.
In order to avoid the problem of the growing wave speeds, we propose a hyperbolic relaxation of this system based on an equivalent form of the incompressibility constraints, centered at
the interfaces between the layers, which is the adjoint of the non-hydrostatic pressure force.
The interface pressures are relaxed through this form and transported on staggered layers, so
that the relaxed system inherits an exact energy balance. Moreover, the wave speeds of the
linearized system  around the lake at rest are bounded by $\sqrt{gH+c^2}$ for any number of
layers, $H$ being the still-water depth and $c$ the relaxation celerity. This suggests, as the
numerical tests will confirm, that, for a fixed $c$, the time step of an explicit scheme does not need
to be reduced when the number of layers grows.

The relaxed system is discretized on a staggered mesh. The hydrostatic part is solved by a
path-conservative finite volume method. The non-hydrostatic part is solved in a single step,
explicit in the horizontal velocities and Crank--Nicolson in the vertical velocities and the
pressures of each water column, so that only one tridiagonal system of size $N$ per column has to
be solved. We prove in one dimension that this step conserves the energy up to a term of order $\tau^2$ (the energy computed with the half-step velocities) and that it is
stable for uniform states under the condition $c\tau\le\Delta x$, $\tau$ being the time step of this part and
$\Delta x$ the mesh size. Wet-dry fronts are treated
without switching off the non-hydrostatic pressure. Since the scheme is not
asymptotic-preserving when $c\to\infty$, the influence of the relaxation celerity is studied by
means of exact solitary waves of the relaxed system and of LDNH$_0$. Finally, the scheme is extended to two dimensions and implemented
on GPU, where it is up to $50$ times faster than a projection method for the same problem and
mesh.

\section{The non-hydrostatic multilayer model}\label{sec:model}
In this Section, we introduce the layer-averaged non-hydrostatic system LDNH$_0$ that was proposed in
\cite{FNPPSM2018}, and we rewrite its incompressibility constraints in an equivalent form, which
will be the basis of the relaxation proposed in Section~\ref{sec:relax}.

\begin{figure}[htbp]
\centering
\begin{tikzpicture}[x=1cm,y=1cm,
  hu/.style={circle,fill=black,inner sep=0pt,minimum size=5pt},
  ww/.style={rectangle,fill=blue!70!black,inner sep=0pt,minimum size=5pt},
  pp/.style={circle,draw=red!75!black,fill=white,thick,inner sep=0pt,minimum size=6pt},
  lab/.style={font=\scriptsize,inner sep=1pt}]
\def\NL{3}
\pgfmathdeclarefunction{zb}{1}{\pgfmathparse{0.35+0.45*exp(-((#1-4.6)^2)/1.2)}}
\pgfmathdeclarefunction{eta}{1}{\pgfmathparse{3.70+0.30*sin(deg(0.9*#1))}}
\pgfmathdeclarefunction{zi}{2}{\pgfmathparse{zb(#1)+(#2)/\NL*(eta(#1)-zb(#1))}}
\def\xo{2.4}
\fill[brown!25] plot[domain=0:7,samples=60] (\x,{zb(\x)}) -- (7,0) -- (0,0) -- cycle;
\draw[very thick,brown!60!black] plot[domain=0:7,samples=60] (\x,{zb(\x)});
\foreach \a in {1,2} \draw[gray] plot[domain=0:7,samples=60] (\x,{zi(\x,\a)});
\draw[thick,blue!60!black] plot[domain=0:7,samples=60] (\x,{eta(\x)});
\node[lab,anchor=west] at (7.05,{zi(7,0)}) {$z_{1/2}=z_b$};
\node[lab,anchor=west] at (7.05,{zi(7,1)}) {$z_{3/2}$};
\node[lab,anchor=west] at (7.05,{zi(7,2)}) {$z_{5/2}$};
\node[lab,anchor=west] at (7.05,{zi(7,3)}) {$z_{7/2}=\eta$};
\draw[dashed,gray] (\xo,0) -- (\xo,{eta(\xo)+0.4});
\node[lab,below] at (\xo,0) {$x_0$};
\foreach \a in {1,...,\NL}{
  \node[hu] at (\xo,{zi(\xo,\a-0.5)}) {};
  \draw[->,thick] (\xo+0.1,{zi(\xo,\a-0.5)}) -- (\xo+0.75,{zi(\xo,\a-0.5)});
  \node[ww] at (\xo-0.3,{zi(\xo,\a-0.5)}) {};
  \node[lab,anchor=south west] at (\xo+0.1,{zi(\xo,\a-0.5)+0.03}) {$u_{\a}$};
  \node[lab,anchor=east] at (\xo-0.4,{zi(\xo,\a-0.5)}) {$w_{\a}$};
}
\foreach \a/\t in {0/{1/2},1/{3/2},2/{5/2}}{
  \node[pp] at (\xo,{zi(\xo,\a)}) {};
  \node[lab,anchor=south west] at (\xo+0.08,{zi(\xo,\a)+0.02}) {$p_{\t}$};
}
\node[pp] at (\xo,{zi(\xo,3)}) {};
\node[lab,anchor=south west] at (\xo+0.08,{zi(\xo,3)+0.04}) {$p_{7/2}=0$};
\draw[<->,thick] (6.3,{zb(6.3)}) -- node[lab,right] {$h$} (6.3,{eta(6.3)});
\draw[<->] (4.0,{zi(4.0,1)}) -- node[lab,right] {$h_2=h/N$} (4.0,{zi(4.0,2)});
\node[font=\small] at (-0.3,4.5) {(a)};
\pgfmathsetmacro{\zA}{zi(\xo,0)} \pgfmathsetmacro{\zB}{zi(\xo,1)}
\pgfmathsetmacro{\zC}{zi(\xo,2)} \pgfmathsetmacro{\zD}{zi(\xo,3)}
\begin{scope}[shift={(9.6,0)}]
  \node[font=\small] at (-0.5,4.5) {(b)};
  \foreach \s/\t in {0/u,1.7/w,3.4/p}{
    \draw[->,gray] (\s,\zA-0.15) -- (\s,\zD+0.4);
    \foreach \z in {\zA,\zB,\zC,\zD} \draw[gray,densely dotted] (\s-0.1,\z) -- (\s+1.3,\z);
    \node[lab,below] at (\s+0.6,\zA-0.15) {$\t$};
    \node[lab,left] at (\s,\zD+0.35) {$z$};
  }
  \draw[very thick] (0.55,\zA) -- (0.55,\zB) (0.85,\zB) -- (0.85,\zC) (1.1,\zC) -- (1.1,\zD);
  \draw[very thick,blue!70!black] (1.7+0.1,\zA) -- (1.7+0.4,\zB) (1.7+0.6,\zB) -- (1.7+0.8,\zC)
    (1.7+1.0,\zC) -- (1.7+1.15,\zD);
  \foreach \ya/\yb/\wa/\wb in {\zA/\zB/0.1/0.4,\zB/\zC/0.6/0.8,\zC/\zD/1.0/1.15}
    \node[ww] at ({1.7+(\wa+\wb)/2},{(\ya+\yb)/2}) {};
  \draw[very thick,red!75!black] (3.4+1.1,\zA) -- (3.4+0.75,\zB) -- (3.4+0.35,\zC) -- (3.4,\zD);
  \foreach \z/\v in {\zA/1.1,\zB/0.75,\zC/0.35,\zD/0} \node[pp] at (3.4+\v,\z) {};
\end{scope}
\begin{scope}[shift={(0.4,-0.8)}]
  \node[hu] at (0,0) {}; \node[lab,anchor=west] at (0.15,0) {$u_\alpha$: horizontal velocity of layer $\alpha$};
  \node[ww] at (5.0,0) {}; \node[lab,anchor=west] at (5.15,0) {$w_\alpha$: vertical velocity of layer $\alpha$};
  \node[pp] at (10.0,0) {}; \node[lab,anchor=west] at (10.15,0) {$p_\aph$: pressure at the interfaces};
\end{scope}
\end{tikzpicture}

\caption{Multilayer decomposition for $N=3$ layers. (a) Vertical section: the velocities
$u_\alpha$, $w_\alpha$ of each layer and the non-hydrostatic pressures $p_\aph$ at the
interfaces, at a section $x=x_0$. (b) Vertical profiles at $x=x_0$ in LDNH$_0$: the horizontal velocity
is constant in each layer, the vertical velocity is linear in each layer and discontinuous at the interfaces, with mean $w_\alpha$ (squares), while the pressure is continuous, linear in each layer and zero at
the free surface.}
\label{fig:layers}
\end{figure}
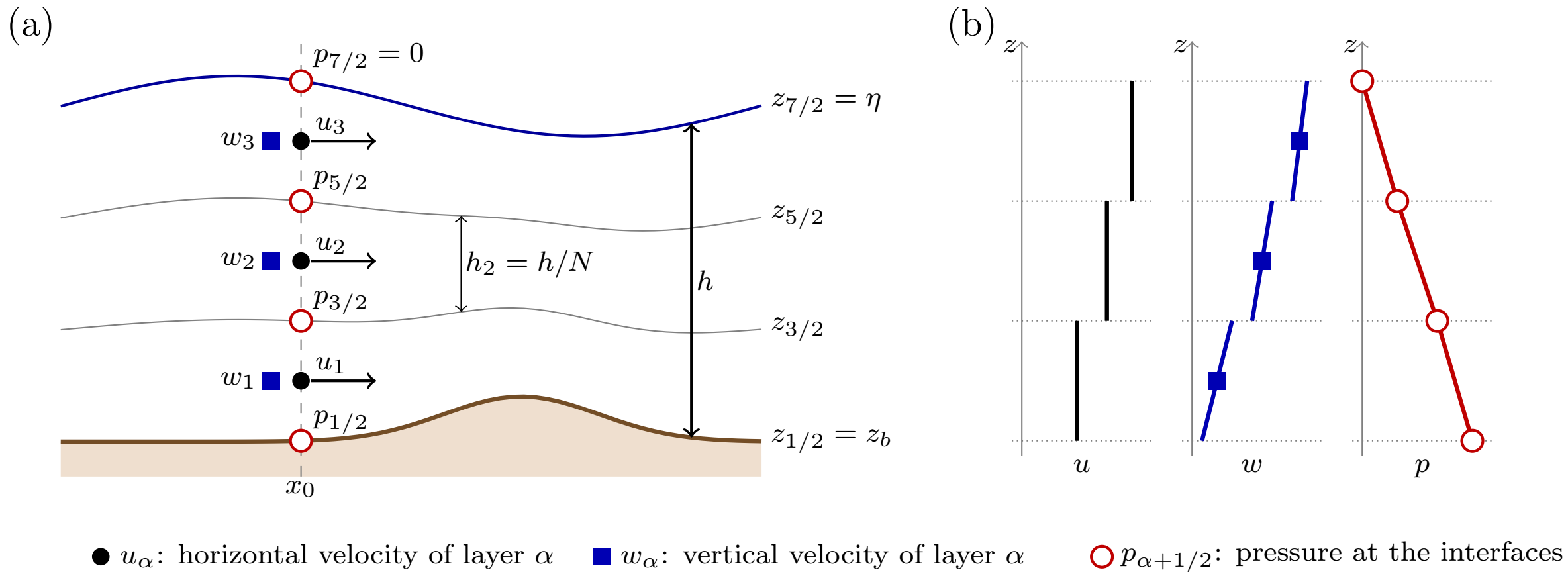

Let us first introduce the notation that will be used throughout the paper.  We consider a one-dimensional free-surface flow over a possibly moving bottom $z_b(x,t)$, and let us
denote by $h(x,t)$ the water depth and by $\eta=z_b+h$ the free surface. The water column is
decomposed along the vertical direction into $N$ layers of \emph{uniform} thickness
$h_\alpha=h/N$, $\alpha=1,\dots,N$, separated by the interfaces
\[
z_\aph=z_b+\frac{\alpha}{N}h,\qquad \alpha=0,\dots,N,
\]
so that $z_{1/2}=z_b$ and $z_{N+1/2}=\eta$ (see Figure~\ref{fig:layers}(a)). Within each layer, $u_\alpha$ and $w_\alpha$ denote
the layer averages of the horizontal and vertical velocities, respectively, and the column
depth-averaged velocity is given by
\[
\ub=\frac1N\sum_{\alpha=1}^N u_\alpha .
\]
We shall also denote by $q_\alpha=h_\alpha u_\alpha$ the momentum of layer $\alpha$. The
non-hydrostatic pressure is described by its values $p_\aph$ at the interfaces, which are
assumed to vanish at the free surface, $p_{N+1/2}=0$. We will denote by $p_\alpha$ the pressure at the center of each layer, defined by
\[
p_\alpha=\frac12\left(p_\aph+p_\amh\right).
\]
Finally, for any layer quantity $f_\alpha$ we write
\[
\tilde f_\aph=\frac12\left(f_\alpha+f_{\alpha+1}\right).
\]

\subsection{The \texorpdfstring{LDNH$_0$}{LDNH0} system}
The system LDNH$_0$ is obtained in \cite{FNPPSM2018} by averaging the Euler equations over
each layer, with a discontinuous Galerkin approximation in the vertical direction: the horizontal
velocity is constant within each layer, the vertical velocity is linear in each layer and discontinuous at the interfaces, and the non-hydrostatic pressure is
continuous and linear in each layer (see Figure~\ref{fig:layers}(b)).
Therefore, the pressure is defined by its values $p_\aph$ at the interfaces, and its value $p_\alpha$ at the center of each layer is the mean of the two values at the interfaces.
The vertical velocity of each layer is represented by its mean $w_\alpha$, which is also its value at the center of the layer, and its slope is given by the incompressibility condition (see below). 

The mass exchange between the layers $\alpha$ and $\alpha+1$ is defined by $\Gamma_\aph$ and the mass balance of each layer reads \begin{equation}\label{eq:layermass}
\dt h_\alpha+\dx(h_\alpha u_\alpha)+\Gamma_\aph-\Gamma_\amh=0 .
\end{equation}

Taking into account that $h_\alpha=h/N$, summing over all layers, one gets the mass balance of the whole water column,
\[
\dt h+\dx(h\ub)=0,
\]
and the mass exchange is therefore given by
\[
\Gamma_\aph=\sum_{\beta=\alpha+1}^{N}\dx\bigl(h_\beta(u_\beta-\ub)\bigr),\qquad
\Gamma_{1/2}=\Gamma_{N+1/2}=0.
\]

The incompressibility condition, integrated from the bottom, will give one
constraint per layer.

The model can be written as (see Section~3.4 and Appendix~C in \cite{FNPPSM2018})
\begin{equation}\label{eq:LDNH}
\left\{
\begin{aligned}
&\dt h+\dx(h\ub)=0,\\
&\dt(h_\alpha u_\alpha)+\dx(h_\alpha u_\alpha^2)+\tilde u_\aph\Gamma_\aph-\tilde u_\amh\Gamma_\amh
 +gh_\alpha\dx\eta=\mathcal D_{1,\alpha},\\
&\dt(h_\alpha w_\alpha)+\dx(h_\alpha u_\alpha w_\alpha)+\tilde w_\aph\Gamma_\aph-\tilde w_\amh\Gamma_\amh
 =\mathcal D_{2,\alpha},\\
&\mathcal C_\alpha=0,
\end{aligned}
\right.
\end{equation}
for $\alpha=1,\dots,N$, where the non-hydrostatic terms are given by
\begin{equation}\label{eq:D12}
\mathcal D_{1,\alpha}=-\dx(h_\alpha p_\alpha)+p_\aph\dx z_\aph-p_\amh\dx z_\amh,\qquad
\mathcal D_{2,\alpha}=-\left(p_\aph-p_\amh\right),
\end{equation}
and the incompressibility constraints are
\begin{equation}\label{eq:C}
\mathcal C_\alpha:=h_\alpha\dx u_\alpha+2\left(w_\alpha-w^+_\amh\right),\qquad
w^+_\amh=\dt z_b+u_\alpha\dx z_\amh-\sum_{\beta=1}^{\alpha-1}\dx(h_\beta u_\beta).
\end{equation}

Note that $w^+_\amh$ is not an unknown of the model.
Let us consider the vertical velocity $w(x,z,t)$ obtained by integrating the incompressibility condition $\dx u+\partial_zw=0$ from the bottom, where the kinematic condition $w=\dt z_b+u\dx z_b$ holds, for the horizontal velocity of the model, which is equal to $u_\alpha$ in the layer $\alpha$.
Unlike the unknown $w_\alpha$, this vertical velocity is not constant in the layers.
In the layer $\alpha$, it is linear in $z$,
\[
w(x,z,t)=w^+_\amh-\left(z-z_\amh\right)\dx u_\alpha,\qquad z_\amh<z<z_\aph,
\]
and it is discontinuous at the interfaces, where the horizontal velocity jumps.
Then, $w^+_\amh$ is its value just above the interface $z_\amh$.
On the other hand, the constraint $\mathcal C_\alpha=0$ reads
\[
w_\alpha=w^+_\amh-\frac{h_\alpha}2\,\dx u_\alpha,
\]
that is, it means that $w_\alpha$ is the value of this linear profile at the middle of the layer, which is also its average over the layer.
The $N$ interface pressures
$p_{1/2},\dots,p_{N-1/2}$ play the role of the Lagrange multipliers associated with the $N$
constraints.

The term $\dt z_b$ in \eqref{eq:C} allows bottom deformations that will then generate waves in the water column. In practice, for the proofs of the energy balance and for the description of the scheme, we will assume a fixed bottom, $\dt z_b=0$. Otherwise the evolution of the bottom should be prescribed and it will appear as a source term in the model.

\subsection{Equivalent form of the constraints}\label{sec:model:equiv}
Due to the term $w^+_\amh$, the constraint $\mathcal C_\alpha$ contains the sum of $\dx(h_\beta u_\beta)$
over all the layers below the layer $\alpha$. This will result in a coupling between the definition of the constraints, so that each constraint is coupled with all the layers underneath. To avoid this coupling, we propose an equivalent form of the constraints, which takes into account the half difference of two consecutive constraints. Since $\mathcal C_\alpha=0$, one has
\begin{equation}\label{eq:Kdef}
\mathcal K_\aph:=\frac12\left(\mathcal C_{\alpha+1}-\mathcal C_\alpha\right)=0,
\end{equation}
for $\alpha=0,\dots,N-1$, with $\mathcal C_0:=0$. Conversely, since $\mathcal C_\alpha=2\sum_{\beta=0}^{\alpha-1}\mathcal K_{\beta+1/2}$, having $\mathcal K_\aph=0$ for all $\alpha$ gives rise to $\mathcal C_\alpha=0$ for all $\alpha$.

Now, let us explicitly compute $\mathcal K_\aph$ in terms of the layer variables. From \eqref{eq:C}, and since
$z_\aph-z_\amh=h_\alpha$, we have
\[
\begin{aligned}
w^+_\aph-w^+_\amh&=u_{\alpha+1}\dx z_\aph-u_\alpha\dx z_\amh-\dx(h_\alpha u_\alpha)\\
&=-h_\alpha\dx u_\alpha+\left(u_{\alpha+1}-u_\alpha\right)\dx z_\aph ,
\end{aligned}
\]
so that the sum that previously appeared in the constraints now cancels, and we obtain
\begin{equation}\label{eq:K}
\mathcal K_\aph=\frac12\left(h_\alpha\dx u_\alpha+h_{\alpha+1}\dx u_{\alpha+1}\right)
+\left(w_{\alpha+1}-w_\alpha\right)-\left(u_{\alpha+1}-u_\alpha\right)\dx z_\aph .
\end{equation}
For $\alpha=0$ we have $w^+_{1/2}=u_1\dx z_b$, and \eqref{eq:K} also holds with the convention
$u_0=w_0=0$:
\[
\mathcal K_{1/2}=\frac12h_1\dx u_1+w_1-u_1\dx z_b .
\]

Unlike $\mathcal C_\alpha$, the condition $\mathcal K_\aph=0$ only involves the two layers
adjacent to the interface $z_\aph$, and it is a divergence-free condition centered at this
interface. Indeed, for $\alpha\ge1$ it reads
\[
\frac{w_{\alpha+1}-w_\alpha}{h/N}+\dx\Bigl(\frac{u_\alpha+u_{\alpha+1}}{2}\Bigr)
-\frac{u_{\alpha+1}-u_\alpha}{h/N}\,\dx z_\aph=0 ,
\]
which is nothing but the incompressibility condition on a vertically staggered grid, with the
velocities located at the centers of the layers and the pressures at the interfaces.

\subsection{Energy balance and pressure work}
As it was proved in \cite{FNPPSM2018} (see Proposition~1 and Appendix~C.2 there), smooth
solutions of \eqref{eq:LDNH} satisfy the energy balance
\begin{equation}\label{eq:energyLDNH}
\dt\Bigl[\sum_{\alpha}h_\alpha\Bigl(\frac{u_\alpha^2}{2}+\frac{w_\alpha^2}{2}\Bigr)+\frac12g\eta^2\Bigr]
+\dx\Bigl[\sum_\alpha h_\alpha u_\alpha\Bigl(\frac{u_\alpha^2}{2}+\frac{w_\alpha^2}{2}+g\eta+p_\alpha\Bigr)\Bigr]=0 .
\end{equation}
In the proof, the advective terms and the mass exchanges between the layers are handled by
means of the following identity, which goes back to the hydrostatic multilayer system of
\cite{ABPSM2011}: if $f_\alpha$, $\alpha=1,\dots,N$, satisfies
\[
\dt(h_\alpha f_\alpha)+\dx(h_\alpha u_\alpha f_\alpha)+\tilde f_\aph\Gamma_\aph-\tilde f_\amh\Gamma_\amh=S_\alpha ,
\]
then
\begin{equation}\label{eq:exchange}
\dt\Bigl(\sum_{\alpha}h_\alpha\frac{f_\alpha^2}{2}\Bigr)
+\dx\Bigl(\sum_{\alpha}h_\alpha u_\alpha\frac{f_\alpha^2}{2}\Bigr)=\sum_\alpha f_\alpha S_\alpha .
\end{equation}

In \cite{FNPPSM2018}, the work of the pressure is computed by using the constraints, and it
reduces to a telescoping sum that vanishes. The following identity holds instead for arbitrary
velocities and pressures.

\begin{proposition}[Pressure work]\label{prop:work}
For any smooth $u_\alpha,w_\alpha,p_\aph$ with $p_{N+1/2}=0$,
\[
\sum_{\alpha=1}^N\left(u_\alpha\mathcal D_{1,\alpha}+w_\alpha\mathcal D_{2,\alpha}\right)
=-\dx\Bigl(\sum_{\alpha=1}^N h_\alpha u_\alpha p_\alpha\Bigr)+\sum_{\alpha=0}^{N-1}p_\aph\,\mathcal K_\aph .
\]
\end{proposition}
\begin{proof}
Since
\[
u_\alpha\dx(h_\alpha p_\alpha)=\dx(h_\alpha u_\alpha p_\alpha)-h_\alpha p_\alpha\dx u_\alpha,
\qquad
h_\alpha p_\alpha\dx u_\alpha=\frac12h_\alpha\dx u_\alpha\left(p_\aph+p_\amh\right),
\]
the definition \eqref{eq:D12} of the non-hydrostatic forces gives
\[
\begin{aligned}
\sum_{\alpha=1}^N\left(u_\alpha\mathcal D_{1,\alpha}+w_\alpha\mathcal D_{2,\alpha}\right)
={}&-\dx\Bigl(\sum_{\alpha=1}^N h_\alpha u_\alpha p_\alpha\Bigr)
+\sum_{\alpha=1}^N\Bigl(\frac12h_\alpha\dx u_\alpha+u_\alpha\dx z_\aph-w_\alpha\Bigr)p_\aph\\
&+\sum_{\alpha=1}^N\Bigl(\frac12h_\alpha\dx u_\alpha-u_\alpha\dx z_\amh+w_\alpha\Bigr)p_\amh .
\end{aligned}
\]
Let us collect the terms that multiply each $p_\aph$, $\alpha=0,\dots,N-1$, recalling that
$p_{N+1/2}=0$. The pressure $p_\aph$ appears as the upper pressure of layer $\alpha$ and as the
lower pressure of layer $\alpha+1$, and the sum of both contributions is
\[
\frac12\left(h_\alpha\dx u_\alpha+h_{\alpha+1}\dx u_{\alpha+1}\right)
+\left(w_{\alpha+1}-w_\alpha\right)-\left(u_{\alpha+1}-u_\alpha\right)\dx z_\aph
=\mathcal K_\aph ,
\]
where, in agreement with the convention $u_0=w_0=0$, only layer $1$ contributes when
$\alpha=0$.
\end{proof}

In other words, the non-hydrostatic pressure force $(\mathcal D_1,\mathcal D_2)$ is the
adjoint of the operator $\mathcal K$, for arbitrary velocities and pressures. This is the
property that will be used for the relaxed system of Section~\ref{sec:relax}, which does not
satisfy the constraints. This property will be preserved at the discrete level in
Section~\ref{sec:scheme:stag}. When the constraints hold, the term $\sum_\alpha
p_\aph\mathcal K_\aph$ vanishes, and the energy balance \eqref{eq:energyLDNH} of LDNH$_0$ is
recovered from \eqref{eq:exchange} and Proposition~\ref{prop:work}.

\section{Hyperbolic relaxation}\label{sec:relax}
The hyperbolic relaxation of dispersive shallow water models has been studied in several previous works, mostly for a single layer.
Favrie and Gavrilyuk \cite{FavrieGavrilyuk2017} introduced an extended Lagrangian for the Serre--Green--Naghdi equations, with a new variable that approximates the water depth and a penalty parameter, whose Euler--Lagrange equations form a hyperbolic system.
The convergence of its solutions when the penalty parameter goes to infinity was proved in \cite{Duchene2019}.
In \cite{Escalante2019JCP}, the divergence constraint of a depth-averaged non-hydrostatic model is replaced by an evolution equation for the non-hydrostatic pressure, and the resulting system satisfies an energy balance.
This idea was extended in \cite{Bassi2020} to the Serre--Green--Naghdi equations with general topographies, by means of evolution equations for the depth-averaged and the bottom non-hydrostatic pressures.
A different approach, based on a weak compressibility of the fluid, was followed in \cite{Richard2021}, while the stiffness of the relaxed system for large relaxation celerities was dealt with in \cite{Macca2025} by means of IMEX Runge--Kutta methods.
For stratified flows, a hyperbolic approximation of a three-layer Green--Naghdi model, in which only the two outer layers are non-hydrostatic, was proposed in \cite{ChesnokovLiapidevskii2020}.
However, to the best of our knowledge, the relaxation of multilayer non-hydrostatic models such as LDNH$_0$ has not been considered so far.

In this Section, the incompressibility constraints of LDNH$_0$ are replaced by evolution
equations for the interface pressures, built on the equivalent constraints $\mathcal K_\aph$
introduced in Section~\ref{sec:model}. To obtain an exact energy balance, the pressures are
transported in a way that is compatible with the mass exchanges between the layers.

\subsection{The relaxed system}
In the relaxed system, the constraints are replaced by evolution equations for the interface pressures $p_\aph$, in which the pressures are transported by the flow.
In order to keep an exact energy balance, this transport is defined on \emph{staggered layers}: the staggered layer associated with $p_\aph$
goes from the center of layer $\alpha$ to the center of layer $\alpha+1$, while the lowest one is
the lower half of layer $1$. Their thickness and mass flux are defined as
\[
\hs_\aph=\frac12\left(h_\alpha+h_{\alpha+1}\right),\qquad
\qs_\aph=\frac12\left(q_\alpha+q_{\alpha+1}\right),
\]
while the exchange flux between two consecutive staggered layers, located at the center of
layer $\alpha$, is
\[
\Gb_\alpha=\frac12\left(\Gamma_\amh+\Gamma_\aph\right),
\]
with $h_0=q_0=0$ and $\Gb_0=0$. Thus,
\[
\hs_{1/2}=\frac{h}{2N},\qquad \hs_\aph=\frac{h}{N}\quad\text{for }\alpha\ge1,
\]
and we write
\[
\hs_\aph=\sigma_\alpha\frac{h}{N},\qquad \textrm{ where } \sigma_0=\frac12, \textrm{ and }  \sigma_\alpha=1\quad\text{for }\alpha\ge1 .
\]
Averaging \eqref{eq:layermass} over the layers $\alpha$ and $\alpha+1$, the mass balance of the
staggered layers is obtained:
\begin{equation}\label{eq:stagmass}
\dt\hs_\aph+\dx\qs_\aph+\Gb_{\alpha+1}-\Gb_\alpha=0,\qquad \alpha=0,\dots,N-1 .
\end{equation}

Let $c>0$ be a relaxation celerity, which is constant in space and time, and the same for all the layers. The relaxed
system is obtained by replacing the constraints in \eqref{eq:LDNH} by evolution equations for the
interface pressures, whose right-hand sides are given by the equivalent constraints
$\mathcal K_\aph$:
\begin{equation}\label{eq:relax}
\left\{
\begin{aligned}
&\dt h+\dx(h\ub)=0,\\
&\dt(h_\alpha u_\alpha)+\dx\bigl(h_\alpha u_\alpha^2+h_\alpha p_\alpha\bigr)+\tilde u_\aph\Gamma_\aph-\tilde u_\amh\Gamma_\amh
 +gh_\alpha\dx\eta\\
&\qquad\qquad=p_\aph\dx z_\aph-p_\amh\dx z_\amh,\\
&\dt(h_\alpha w_\alpha)+\dx(h_\alpha u_\alpha w_\alpha)+\tilde w_\aph\Gamma_\aph-\tilde w_\amh\Gamma_\amh
 =-\left(p_\aph-p_\amh\right),\\
&\dt(\hs_\aph p_\aph)+\dx(\qs_\aph p_\aph)+p_{\alpha+1}\Gb_{\alpha+1}-p_\alpha\Gb_\alpha=-c^2\,\mathcal K_\aph,
\end{aligned}
\right.
\end{equation}
where $\alpha=1,\dots,N$ in the momentum equations and $\alpha=0,\dots,N-1$ in the pressure
equations. The system has $3N+1$ unknowns, namely $h$, $(u_\alpha,w_\alpha)_{\alpha=1}^N$ and
$(p_\aph)_{\alpha=0}^{N-1}$.

\subsection{Energy balance}
The first property of the relaxed system is that it inherits an exact energy balance, in which
the energy of LDNH$_0$ in \eqref{eq:energyLDNH} is supplemented by a pressure term weighted by
$c^{-2}$.
The proof is based on the following analogue of identity~\eqref{eq:exchange} for the staggered layers.

\begin{lemma}\label{lem:stagexchange}
Let $p_\aph$, $\alpha=0,\dots,N-1$, satisfy
\[
\dt(\hs_\aph p_\aph)+\dx(\qs_\aph p_\aph)+p_{\alpha+1}\Gb_{\alpha+1}-p_\alpha\Gb_\alpha=S_\aph ,
\]
with $p_{N+1/2}=0$. Then
\[
\dt\Bigl(\sum_{\alpha=0}^{N-1}\hs_\aph\frac{p_\aph^2}{2}\Bigr)
+\dx\Bigl(\sum_{\alpha=0}^{N-1}\qs_\aph\frac{p_\aph^2}{2}\Bigr)=\sum_{\alpha=0}^{N-1}p_\aph S_\aph .
\]
\end{lemma}
\begin{proof}
Multiplying the equation by $p_\aph$ and subtracting $\frac12p_\aph^2$ times
\eqref{eq:stagmass}, we obtain
\[
\begin{aligned}
&\dt\Bigl(\hs_\aph\frac{p_\aph^2}{2}\Bigr)+\dx\Bigl(\qs_\aph\frac{p_\aph^2}{2}\Bigr)\\
&\qquad+\Gb_{\alpha+1}\Bigl(p_{\alpha+1}p_\aph-\frac12p_\aph^2\Bigr)
-\Gb_\alpha\Bigl(p_\alpha p_\aph-\frac12p_\aph^2\Bigr)=p_\aph S_\aph .
\end{aligned}
\]
When summing over $\alpha$, the exchange terms located at the center of layer $\alpha$,
$1\le\alpha\le N-1$, give
\[
\Gb_\alpha\Bigl(p_\amh p_\alpha-\frac12p_\amh^2-p_\aph p_\alpha+\frac12p_\aph^2\Bigr)=0,
\]
since $p_\alpha=\frac12(p_\aph+p_\amh)$. At the center of layer $N$ only the staggered layer
$N-1/2$ is present, and its contribution
\[
\Gb_N\Bigl(p_{N-1/2}\,p_N-\frac12p_{N-1/2}^2\Bigr)
\]
also vanishes, since $p_N=\frac12p_{N-1/2}$. Finally, the term at the bottom vanishes because
$\Gb_0=0$.
\end{proof}

\begin{theorem}[Energy]\label{thm:energy}
Smooth solutions of \eqref{eq:relax} satisfy $\dt\mathcal E+\dx\mathcal F=0$, with
\[
\mathcal E=\sum_{\alpha=1}^N h_\alpha\Bigl(\frac{u_\alpha^2}{2}+\frac{w_\alpha^2}{2}\Bigr)
+\sum_{\alpha=0}^{N-1}\hs_\aph\frac{p_\aph^2}{2c^2}+\frac12g\eta^2,
\]
\[
\mathcal F=\sum_{\alpha=1}^N h_\alpha u_\alpha\Bigl(\frac{u_\alpha^2}{2}+\frac{w_\alpha^2}{2}+g\eta+p_\alpha
{}+\frac{p_\amh^2+p_\aph^2}{4c^2}\Bigr).
\]
\end{theorem}
\begin{proof}
Applying identity~\eqref{eq:exchange} to $f_\alpha=u_\alpha$ and to $f_\alpha=w_\alpha$, and using
Proposition~\ref{prop:work}, we obtain
\[
\begin{aligned}
&\dt\Bigl(\sum_\alpha h_\alpha\frac{u_\alpha^2+w_\alpha^2}{2}\Bigr)
+\dx\Bigl(\sum_\alpha h_\alpha u_\alpha\frac{u_\alpha^2+w_\alpha^2}{2}\Bigr)\\
&\qquad=-\sum_\alpha gh_\alpha u_\alpha\dx\eta
-\dx\Bigl(\sum_\alpha h_\alpha u_\alpha p_\alpha\Bigr)
+\sum_{\alpha=0}^{N-1}p_\aph\mathcal K_\aph .
\end{aligned}
\]
Since the bottom is fixed, $\dt\eta=-\dx(h\ub)$, and the hydrostatic pressure term can be written
as
\[
\sum_\alpha gh_\alpha u_\alpha\dx\eta=gh\ub\,\dx\eta=\dt\Bigl(\frac12g\eta^2\Bigr)+\dx(gh\ub\eta).
\]
Contrary to the case of LDNH$_0$, the term $\sum_\alpha p_\aph\mathcal K_\aph$ does not vanish,
since the constraints are not imposed. On the other hand, Lemma~\ref{lem:stagexchange}, applied
to the pressure equations with
$S_\aph=-c^2\mathcal K_\aph$ and divided by $c^2$, gives
\[
\dt\Bigl(\sum_{\alpha=0}^{N-1}\hs_\aph\frac{p_\aph^2}{2c^2}\Bigr)
+\dx\Bigl(\sum_{\alpha=0}^{N-1}\qs_\aph\frac{p_\aph^2}{2c^2}\Bigr)
=-\sum_{\alpha=0}^{N-1}p_\aph\mathcal K_\aph .
\]
The result follows by combining these three identities, where the last sum of $\mathcal F$ has
been written layer by layer. More explicitly, since $\qs_\aph=\frac12(q_\alpha+q_{\alpha+1})$ with $q_0=0$, and
$p_{N+1/2}=0$, collecting the coefficient of each $q_\alpha$ gives
\[
\sum_{\alpha=0}^{N-1}\qs_\aph\,p_\aph^2=\sum_{\alpha=1}^N q_\alpha\,\frac{p_\amh^2+p_\aph^2}{2}.
\]
\end{proof}

\begin{remark}
Assuming that the solution of \eqref{eq:relax} remains smooth, the relaxed system formally reduces to LDNH$_0$ as $c\to\infty$. Indeed, provided that the variables remain smooth and bounded as $c\to\infty$, then
\[
\mathcal K_\aph=O(c^{-2})\qquad\text{for all }\alpha
\]
and, since $\mathcal C_\alpha=2\sum_{\beta=0}^{\alpha-1}\mathcal K_{\beta+1/2}$
(Section~\ref{sec:model:equiv}), it implies $\mathcal C_\alpha=O(c^{-2})$ and the result follows.
\end{remark}

For $N=1$, the only pressure unknown corresponds to the one at the bottom, $p_{1/2}=2p_1$, and the pressure
equation is precisely the single-layer relaxation of \cite{Escalante2019JCP}, with relaxation
celerity $c/\sqrt2$.

\subsection{Linearized system around the lake at rest}\label{sec:relax:rest}
In order to study the hyperbolicity of the relaxed system and the dependence of its wave speeds on
the number of layers, which motivates the proposed relaxation, we linearize
\eqref{eq:relax} around the lake at rest solution over a flat bottom,
\[
h=H,\qquad z_b=-H,\qquad u_\alpha=w_\alpha=p_\aph=0 .
\]
Since the bottom is flat and $u_\alpha=p_\aph=0$ at rest, the geometric terms
$p_\aph\dx z_\aph$ and $(u_{\alpha+1}-u_\alpha)\dx z_\aph$ are quadratic in the perturbation, and the
linearized equations read
\[
\begin{aligned}
&\dt\eta+\frac HN\sum_{\alpha=1}^N\dx u_\alpha=0,\\
&\dt u_\alpha+\frac12\left(\dx p_\aph+\dx p_\amh\right)+g\,\dx\eta=0,\\
&\dt w_\alpha=-\frac NH\left(p_\aph-p_\amh\right),\\
&\dt p_\aph+\frac{c^2}{2\sigma_\alpha}\left(\dx u_\alpha+\dx u_{\alpha+1}\right)
=-\frac{c^2N}{\sigma_\alpha H}\left(w_{\alpha+1}-w_\alpha\right),
\end{aligned}
\]
with $u_0=w_0=0$ and $p_{N+1/2}=0$.

In order to write these equations in compact form, let us introduce the vectors of the values at the layers,
\[
\mathbf u=(u_1,\dots,u_N)^T,\qquad \mathbf w=(w_1,\dots,w_N)^T,
\]
and the vector of the values at the interfaces,
\[
\mathbf p=(p_{1/2},\dots,p_{N-1/2})^T.
\]
Let $P$ and $D$ be the $N\times N$ matrices that give, from the values at the interfaces, their mean and their difference at the center of each layer:
\[
(P\mathbf p)_\alpha=\tfrac12\bigl(p_\amh+p_\aph\bigr),\qquad
(D\mathbf p)_\alpha=p_\amh-p_\aph,\qquad\alpha=1,\dots,N,
\]
with $p_{N+1/2}=0$.
That is, $P$ has $1/2$ on the diagonal and on the superdiagonal, $D$ has $1$ on the diagonal and $-1$ on the superdiagonal, and their remaining entries vanish.
Their transposes act on the values at the layers and give values at the interfaces, the component $\beta$ corresponding to the interface $\beta+\frac12$,
\[
(P^T\mathbf u)_\beta=\tfrac12\bigl(u_\beta+u_{\beta+1}\bigr),\qquad
(D^T\mathbf w)_\beta=w_{\beta+1}-w_\beta,\qquad\beta=0,\dots,N-1,
\]
with $u_0=w_0=0$.
Let us also define the diagonal matrix
$\Sigma=\mathrm{diag}(\sigma_0,\dots,\sigma_{N-1})$ and the vector $\mathbf e=(1,\dots,1)^T$.

The linearized system can be symmetrized by means of the energy.
At rest, the thicknesses are $h_\alpha=H/N$ and $\hs_\aph=\sigma_\alpha H/N$, and the quadratic part of the energy $\mathcal E$ of Theorem~\ref{thm:energy} is
\[
\frac12g\eta^2+\frac HN\sum_{\alpha=1}^N\frac{u_\alpha^2+w_\alpha^2}{2}
+\frac{H}{Nc^2}\sum_{\alpha=0}^{N-1}\sigma_\alpha\frac{p_\aph^2}{2}.
\]
Let us multiply the equations for $\eta$, $u_\alpha$, $w_\alpha$ and $p_\aph$ by the weights of these variables in this quadratic form, that is, by $g$, $H/N$, $H/N$ and $\sigma_\alpha H/(Nc^2)$, respectively.
Then, the linearized system for the vector
\[
\mathbf X=(\eta,\mathbf u,\mathbf w,\mathbf p)^T\in\mathbb R^{1+3N}
\]
reads
\begin{equation}\label{eq:lin}
M\,\dt\mathbf X+A\,\dx\mathbf X=G\,\mathbf X,
\end{equation}
where $M$ is the diagonal matrix
\[
M=\mathrm{diag}\Bigl(g,\frac HN I,\frac HN I,\frac{H}{Nc^2}\Sigma\Bigr),
\]
$I$ being the identity matrix of size $N$, and
\[
A=\frac HN\begin{pmatrix}
0&g\,\mathbf e^T&0&0\\
g\,\mathbf e&0&0&P\\
0&0&0&0\\
0&P^T&0&0
\end{pmatrix},\qquad
G=\begin{pmatrix}
0&0&0&0\\
0&0&0&0\\
0&0&0&D\\
0&0&-D^T&0
\end{pmatrix}.
\]
The matrix $M$ is diagonal and positive definite, $A$ is symmetric and $G$ is skew-symmetric, so
that the quadratic energy $\frac12\mathbf X^TM\mathbf X$ is conserved. The
eigenvalues that we need are those of the second difference matrix with a Neumann condition at
one end and a Dirichlet condition at the other one, whose eigenvectors are the discrete cosine
bases \cite{Strang1999}.

\begin{lemma}[Vertical modes]\label{lem:modes}
Let $\theta_m=(2m-1)\pi/(2N)$ and $\varphi_m=\theta_m/2$, for $m=1,\dots,N$, and let
\[
\mathbf v_m=\bigl(\cos(\theta_m\beta)\bigr)_{\beta=0}^{N-1},\qquad
\mathbf a_m=\Bigl(\cos\bigl(\theta_m(\alpha-\tfrac12)\bigr)\Bigr)_{\alpha=1}^N,\qquad
\mathbf d_m=\Bigl(\sin\bigl(\theta_m(\alpha-\tfrac12)\bigr)\Bigr)_{\alpha=1}^N .
\]
Then
\[
\begin{aligned}
P\mathbf v_m&=\cos\varphi_m\,\mathbf a_m,&\qquad P^T\mathbf a_m&=\cos\varphi_m\,\Sigma\mathbf v_m,\\
D\mathbf v_m&=2\sin\varphi_m\,\mathbf d_m,&\qquad D^T\mathbf d_m&=2\sin\varphi_m\,\Sigma\mathbf v_m .
\end{aligned}
\]
Therefore, the eigenvalues of $\Sigma^{-1}P^TP$ and of $P\Sigma^{-1}P^T$ are $\cos^2\varphi_m$, and
those of $\Sigma^{-1}D^TD$ are $4\sin^2\varphi_m$, $m=1,\dots,N$.
\end{lemma}
\begin{proof}
Since $\cos(\theta_mN)=0$, for $\alpha=1,\dots,N$ we have
\[
\begin{aligned}
(P\mathbf v_m)_\alpha&=\tfrac12\bigl[\cos\bigl(\theta_m(\alpha-1)\bigr)+\cos(\theta_m\alpha)\bigr]
=\cos\varphi_m\cos\bigl(\theta_m(\alpha-\tfrac12)\bigr),\\
(D\mathbf v_m)_\alpha&=\cos\bigl(\theta_m(\alpha-1)\bigr)-\cos(\theta_m\alpha)
=2\sin\varphi_m\sin\bigl(\theta_m(\alpha-\tfrac12)\bigr).
\end{aligned}
\]
In the same way, $(P^T\mathbf a_m)_\beta=\cos\varphi_m\cos(\theta_m\beta)$ and
$(D^T\mathbf d_m)_\beta=2\sin\varphi_m\cos(\theta_m\beta)$ for $\beta=1,\dots,N-1$, while for
$\beta=0$ only the first layer contributes, and we obtain $\frac12\cos\varphi_m$ and
$\sin\varphi_m$, that is, the same expressions multiplied by $\sigma_0=\frac12$. The $N$ values
$\cos^2\varphi_m$ are different, and so are the values $4\sin^2\varphi_m$, so that these are all
the eigenvalues.
\end{proof}

\subsection{Hyperbolicity and stiffness}\label{sec:relax:hyp}
The horizontal wave speeds of the linearized system are bounded independently of the number of
layers, while the vertical coupling is stiff.

\begin{proposition}[Hyperbolicity at rest]\label{prop:hyp}
The first-order part of the system \eqref{eq:relax} linearized around the lake at rest over a flat bottom is
symmetrizable, hence hyperbolic, and its wave speeds $\lambda$ satisfy
\[
\lambda^2\le gH+c^2\cos^2\Bigl(\frac{\pi}{4N}\Bigr)< gH+c^2\qquad\text{for every }N .
\]
\end{proposition}
\begin{proof}
In \eqref{eq:lin}, $M$ is symmetric and positive definite and $A$ is symmetric, so that the
first-order part is symmetrizable, and the wave speeds are the solutions of
$A\mathbf X=\lambda M\mathbf X$, that is,
\[
\lambda\,\eta=\tfrac HN\,\mathbf e^T\mathbf u,\qquad
\lambda\,\mathbf u=g\,\mathbf e\,\eta+P\mathbf p,\qquad
\lambda\,\mathbf w=0,\qquad
\lambda\,\Sigma\mathbf p=c^2P^T\mathbf u .
\]
The vertical velocities give $N$ zero speeds. For $\lambda\ne0$, eliminating $\eta$ and
$\mathbf p$,
\[
\lambda^2\mathbf u=\Bigl(\frac{gH}{N}\,\mathbf e\,\mathbf e^T+c^2P\,\Sigma^{-1}P^T\Bigr)\mathbf u,
\]
where both matrices are symmetric and positive semidefinite. By Weyl's inequality and
Lemma~\ref{lem:modes},
\[
\lambda^2\le\frac{gH}{N}\,|\mathbf e|^2+c^2\lambda_{\max}\bigl(P\,\Sigma^{-1}P^T\bigr)
=gH+c^2\cos^2\Bigl(\frac{\pi}{4N}\Bigr).
\]
\end{proof}

Note that the bound does not depend on $N$ and that it is almost reached for large $N$: for
$c=3\sqrt{gH}$ the maximal speed is $0.925$, $0.982$ and $0.998$ times $\sqrt{gH+c^2}$ for
$N=2,5,40$. This suggests that, for a fixed $c$, the time step of an explicit scheme for \eqref{eq:relax} does not
need to be reduced when the number of layers grows. Note, however, that
Proposition~\ref{prop:hyp} only concerns the linearization around the lake at rest over a flat bottom: for general
states we have no proof of such a bound, and the time step of the scheme of
Section~\ref{sec:scheme} is based instead on the relaxation celerity $c$ and on the stability
analysis of the non-hydrostatic step carried out in that section. The numerical tests of
Sections~\ref{sec:tests}
and~\ref{sec:gpu} show that this choice is stable for all the numbers of layers considered. On
the contrary, other possible natural choices of the relaxed variables may lead to speeds that grow with $N$
already at rest, as shown in the following remark.

\begin{remark}[Relaxing the pressure jumps]\label{rem:jumps}
A natural alternative consists in relaxing each constraint $\mathcal C_\alpha$ through the
pressure jump across the layer, by means of an equation of the form
\[
\frac{1}{2c_{\mathrm J}^2}\Bigl(\dt\bigl(h_\alpha[p]_\alpha\bigr)+\dots\Bigr)=\frac12\,\mathcal C_\alpha,
\qquad [p]_\alpha=p_\aph-p_\amh .
\]
The resulting system also has an energy and is hyperbolic at rest, but the interface pressures
\[
p_\aph=-\sum_{\beta=\alpha+1}^{N}[p]_\beta
\]
are non-local functions of the unknowns, and the maximal speed grows linearly with $N$: it is
about $1.7\,c_{\mathrm J}$, $9.0\,c_{\mathrm J}$ and $18.0\,c_{\mathrm J}$ for $N=2,10,20$. Relaxing the
interface pressures through $\mathcal K_\aph$ removes this growth.
\end{remark}

While, at rest, the horizontal wave speeds are bounded independently of $N$, the terms that couple the
vertical velocities and the pressures within each water column are stiff, as shown by the
following result.

\begin{proposition}[Stiff vertical coupling]\label{prop:stiff}
Let us consider the terms of \eqref{eq:relax} without $x$-derivatives that couple $w$ and $p$,
at rest:
\[
\frac HN\,\dt w_\alpha=-(p_\aph-p_\amh),\qquad
\hs_\aph\dt p_\aph=-c^2(w_{\alpha+1}-w_\alpha).
\]
This is a linear system that conserves the energy
\[
\sum_{\alpha=1}^N\frac HN\,\frac{w_\alpha^2}{2}+\sum_{\alpha=0}^{N-1}\hs_\aph\frac{p_\aph^2}{2c^2},
\]
and whose eigenvalues are
\[
\pm\,i\,\frac{2cN}{H}\sin\Bigl(\frac{(2m-1)\pi}{4N}\Bigr),\qquad m=1,\dots,N .
\]
\end{proposition}
\begin{proof}
This is the system $M\dt\mathbf X=G\mathbf X$ of \eqref{eq:lin} for $\mathbf w$ and
$\mathbf p$,
\[
\frac HN\,\dt\mathbf w=D\mathbf p,\qquad\frac{H}{Nc^2}\,\Sigma\,\dt\mathbf p=-D^T\mathbf w .
\]
Since $G$ is skew-symmetric, $\frac12\mathbf X^TM\mathbf X$ is conserved, which for $\mathbf w$ and
$\mathbf p$ is the energy of the statement. The eigenvalues $\mu$ satisfy
$\mu^2\,\Sigma\mathbf p=-(cN/H)^2D^TD\,\mathbf p$, so that, by Lemma~\ref{lem:modes},
$\mu^2=-(cN/H)^2\,4\sin^2\varphi_m$.
\end{proof}

Therefore, the frequencies of the vertical coupling grow like $cN/H$, while the horizontal speeds
remain bounded by $\sqrt{gH+c^2}$. An explicit treatment of this coupling would require
\[
\Delta t\lesssim\frac{H}{2cN},
\]
which is much smaller than the time step of the horizontal part when $\Delta x\gg H/N$. This is
the reason why the scheme proposed in Section~\ref{sec:scheme} is explicit in the horizontal
direction and implicit only in the vertical coupling, which is local to each water column.

\subsection{Linear dispersion}\label{sec:relax:lin}
Let us look for plane waves $\mathbf X=\hat{\mathbf X}e^{\mathrm i(kx-\omega t)}$ of the linearized
system \eqref{eq:lin}. Then
\[
\omega\,M\hat{\mathbf X}=\bigl(kA+\mathrm iG\bigr)\hat{\mathbf X},
\]
where $kA+\mathrm iG$ is a Hermitian matrix and $M$ is positive definite, so that all the frequencies
$\omega$ are real. The phase velocity of the physical branch, which is the smallest positive
frequency, is $C=\omega/k$.
The remaining $2N$ non-zero frequencies form the fast branch of the dispersion relation.
For $k=0$, they are the frequencies of the vertical coupling given by Proposition~\ref{prop:stiff}, of order $cN/H$, and they have no counterpart in LDNH$_0$.

These eigenvalues are computed numerically. When $c\to\infty$, the celerity of LDNH$_0$, whose
dispersion relation is given in \cite{FNPPSM2018}, is recovered with an error of order $c^{-2}$.

\begin{proposition}[Convergence of the celerity]\label{prop:disp}
Let $N$ and $k>0$ be fixed, and let $C_0$ be the phase velocity of LDNH$_0$. Then the phase velocity
of the physical branch of the relaxed system satisfies $C=C_0+O(c^{-2})$ as $c\to\infty$.
\end{proposition}
\begin{proof}
Let $\varepsilon=c^{-2}$. The plane wave equations read
\[
\begin{aligned}
\omega\,\eta&=k\tfrac HN\,\mathbf e^T\mathbf u,&\qquad
\omega\,\mathbf u&=kg\,\mathbf e\,\eta+kP\mathbf p,\\
\omega\tfrac HN\,\mathbf w&=\mathrm iD\mathbf p,&\qquad
\varepsilon\,\omega\tfrac HN\,\Sigma\mathbf p&=k\tfrac HN\,P^T\mathbf u-\mathrm iD^T\mathbf w,
\end{aligned}
\]
so that the frequencies are the roots of a polynomial $\chi(\omega,\varepsilon)$, the determinant of
this linear system. For $\varepsilon=0$ the last equation is the constraint, and the system is
the linearized LDNH$_0$. Since the matrix $D$ is invertible, eliminating $\mathbf w$ and $\mathbf p$
gives $\omega L\mathbf u=kg\,\mathbf e\,\eta$, with
\[
L=I+\Bigl(\frac{kH}{N}\Bigr)^2P\,(D^TD)^{-1}P^T
\]
symmetric and positive definite, and then
$\omega^2=\omega_0^2:=k^2\,\frac{gH}{N}\,\mathbf e^TL^{-1}\mathbf e>0$, with $\omega_0=kC_0$.
$\chi(\omega,0)=a\,\omega^{N-1}(\omega^2-\omega_0^2)$ with $a\neq0$, so that $\omega_0$ is a simple root
of $\chi(\cdot,0)$. By the implicit function theorem, there is a smooth root $\omega(\varepsilon)$ of
$\chi(\cdot,\varepsilon)$ with $\omega(0)=\omega_0$, and therefore
$\omega(\varepsilon)=\omega_0+O(\varepsilon)$. Finally, $\omega=0$ is a root of multiplicity $N-1$ for
every $\varepsilon$, because the matrix of the system at $\omega=0$ does not depend on
$\varepsilon$, and the remaining roots are unbounded as $\varepsilon\to0$, since $\chi(\cdot,0)$ has no
other roots. Hence, for $c$ large enough, $\omega(\varepsilon)$ is the smallest positive
frequency.
\end{proof}

\begin{figure}[htbp]
\centering
\includegraphics[width=\textwidth]{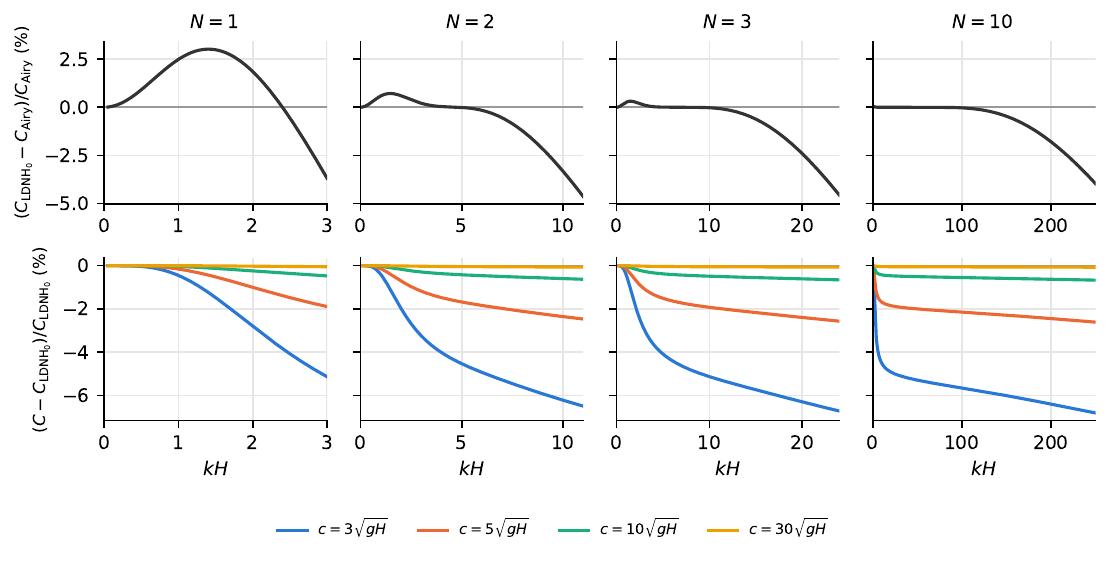}
\caption{Top: relative error of LDNH$_0$ with respect to Airy's formula versus $kH$, for $N=1,2,3$
and $10$ layers, in the range where it is below $5\%$. Bottom: relative error of the relaxed system
with respect to LDNH$_0$, in the same ranges, for several values of $c$.}
\label{fig:disp}
\end{figure}

Figure~\ref{fig:disp} (top) shows the relative error of LDNH$_0$, computed
from its dispersion relation \cite{FNPPSM2018}, with respect to Airy's formula, $C_{\mathrm{Airy}}^2=g\tanh(kH)/k$, for $N=1,2,3$ and $10$ layers. In each
case, the range of $kH$ is close to the one in which the error is below $5\%$, namely $kH\le3.2,11.3,24.7$ and $270$,
so that the range of validity grows quickly with $N$. The error is positive for small $kH$ and
at most about $3\%$ for all $N$. The bottom row shows
the relative error of the relaxed system with respect to LDNH$_0$, in the same ranges and for several
values of $c$. The relaxed system is slower than LDNH$_0$, and the difference grows with $kH$ and
decreases quickly with $c$: it is below $1\%$ for $c\ge10\sqrt{gH}$ and all the values of $N$.

\begin{remark}[Resonance with solitary waves]\label{rem:resonance}
The phase speeds $\omega/k$ of the fast branch are not bounded from below by $c$.
For large wave numbers, they approach the non-zero wave speeds of the first-order part of the linearized system, given in the proof of Proposition~\ref{prop:hyp}.
By Weyl's inequality and Lemma~\ref{lem:modes}, the smallest of these speeds is not smaller than
\[
c\cos\varphi_N=c\sin\Bigl(\frac{\pi}{4N}\Bigr)\approx0.8\,\frac cN,
\]
and it is close to this bound in the computations of the linearized system for $N=1,2,3,6$.
Therefore, a solitary wave of speed $C>\sqrt{gH}$ is in linear resonance with the fast branch whenever $0.8\,c/N<C$, and the relaxed system then only admits generalized solitary waves, with small oscillating tails of short wavelength.
Localized solitary waves require, roughly, $c\gtrsim1.25\,N\,C$, that is, $c\gtrsim1.4\,N\sqrt{gH}$ for $A/H\approx0.2$--$0.3$, although these tails are not visible in our tests.
This does not change the stability condition, but, if $c$ is increased with $N$, the number of sub-steps grows accordingly.
\end{remark}

\section{Numerical scheme}\label{sec:scheme}
The aim of this Section is to describe a numerical scheme for the relaxed system \eqref{eq:relax}
that is explicit in the horizontal direction and implicit only in the vertical coupling of each
water column, so that, except at $t=0$ (see Section~\ref{sec:scheme:init}), no linear system
larger than $N\times N$ has to be solved. The price to be paid is a time step restricted by the
relaxation celerity $c$. In order to reduce its impact, sub-steps are used, so that this
restriction only affects the non-hydrostatic part, which is linear and cheap. Note the difference
with the semi-implicit and IMEX relaxation schemes for multi-scale problems
\cite{PareschiRusso2005,ThomannIolloPuppo2023,Thomann2025},
in which the fast waves are treated implicitly and global linear systems arise: here the fast
horizontal waves are resolved explicitly by the sub-steps, and only the stiff vertical coupling,
which is local to each water column, is treated implicitly.

\subsection{Time splitting}\label{sec:split:time}
System \eqref{eq:relax} is split into two subsystems.
The first one is the hydrostatic part, that is, system \eqref{eq:relax} without the non-hydrostatic terms:
\begin{equation}\tag{H}\label{eq:H}
\left\{
\begin{aligned}
&\dt h+\dx(h\ub)=0,\\
&\dt(h_\alpha u_\alpha)+\dx(h_\alpha u_\alpha^2)+\tilde u_\aph\Gamma_\aph-\tilde u_\amh\Gamma_\amh
 +gh_\alpha\dx\eta=0,\\
&\dt(h_\alpha w_\alpha)+\dx(h_\alpha u_\alpha w_\alpha)+\tilde w_\aph\Gamma_\aph-\tilde w_\amh\Gamma_\amh=0,\\
&\dt(\hs_\aph p_\aph)+\dx(\qs_\aph p_\aph)+p_{\alpha+1}\Gb_{\alpha+1}-p_\alpha\Gb_\alpha=0.
\end{aligned}
\right.
\end{equation}
It corresponds to the hydrostatic multilayer
system of \cite{ABPSM2011}, in which $w_\alpha$ and $p_\aph$ are simply transported. By
identity~\eqref{eq:exchange} and Lemma~\ref{lem:stagexchange}, it conserves $\mathcal E$ for smooth solutions.
The second one is the non-hydrostatic part, in which the depth $h$ is frozen:
\begin{equation}\tag{NH}\label{eq:NH}
h_\alpha\dt u_\alpha=\mathcal D_{1,\alpha},\qquad h_\alpha\dt w_\alpha=\mathcal D_{2,\alpha},\qquad
\hs_\aph\dt p_\aph=-c^2\mathcal K_\aph .
\end{equation}
This is a linear system in $(u,w,p)$ which, by Proposition~\ref{prop:work}, conserves the energy
\[
\sum_{\alpha=1}^N h_\alpha\frac{u_\alpha^2+w_\alpha^2}{2}+\sum_{\alpha=0}^{N-1}\hs_\aph\frac{p_\aph^2}{2c^2},
\]
with flux $\sum_\alpha h_\alpha u_\alpha p_\alpha$.

The time step of each subsystem is restricted by its own characteristic speeds, which are given
by the following result.

\begin{proposition}[Characteristic speeds of the subsystems]\label{prop:subsystems}
(i) The characteristic speeds of \eqref{eq:H} are those of the hydrostatic multilayer system for the
unknowns $(h,h_\alpha u_\alpha)$, together with the transport speeds $u_\alpha$ of the vertical
velocities and $\hat u_\aph=\qs_\aph/\hs_\aph$ of the pressures.

(ii) Let $h>0$ be frozen. For any bottom and any state, \eqref{eq:NH} is a linear symmetrizable
hyperbolic system whose characteristic speeds are $0$, with multiplicity $N$, and
\[
\pm\,c\cos\Bigl(\frac{(2m-1)\pi}{4N}\Bigr),\qquad m=1,\dots,N .
\]
In particular, they are bounded by $c\cos(\pi/(4N))<c$, independently of the number of layers
and of the state.
\end{proposition}
\begin{proof}
Let us first prove (i). The equations for $h$ and $h_\alpha u_\alpha$ in \eqref{eq:H} do not involve $w$
and $p$. In the
equations for $w_\alpha$ and $p_\aph$, written in non-conservative form (see
Section~\ref{sec:scheme:H}), the only derivatives of $w$ and $p$ are those of the transport terms
$u_\alpha\dx w_\alpha$ and $\hat u_\aph\dx p_\aph$, since the exchange terms involve $w$ and $p$
but not their derivatives. Therefore, the matrix of the quasilinear form of \eqref{eq:H} is block lower
triangular, and its diagonal blocks are the Jacobian matrix of the hydrostatic multilayer system
and the diagonal matrices of the transport speeds.

Let us now prove (ii). The characteristic speeds only depend on the terms that contain derivatives of the unknowns $(u,w,p)$.
In \eqref{eq:NH}, these terms are $-h_\alpha\dx p_\alpha$ in $\mathcal D_{1,\alpha}$ and
$\frac12(h_\alpha\dx u_\alpha+h_{\alpha+1}\dx u_{\alpha+1})$ in $\mathcal K_\aph$, while the remaining
terms, which involve $\dx h$, $\dx z$ and $w$, do not contain derivatives of the unknowns.
Dividing all the equations by $h/N$ and keeping only the terms with derivatives of the unknowns, \eqref{eq:NH} reads
\[
\dt\mathbf u+P\,\dx\mathbf p=\cdots,\qquad
\Sigma\,\dt\mathbf p+c^2P^T\dx\mathbf u=\cdots,
\]
where the dots stand for the terms without derivatives of the unknowns, and $P$ and $\Sigma$ are the matrices of \eqref{eq:lin}, whose coefficients do not depend on $h$. As
in the proof of Proposition~\ref{prop:hyp}, it is symmetrizable, and its characteristic speeds
satisfy $\lambda^2\,\Sigma\mathbf p=c^2P^TP\,\mathbf p$, that is,
$\lambda=\pm c\cos\bigl((2m-1)\pi/(4N)\bigr)$ by Lemma~\ref{lem:modes}. The vertical
velocities have no $x$-derivatives and give $N$ zero speeds.
\end{proof}

Let us denote by $\mathcal S_{\mathrm H}(s)$ and $\mathcal S_{\mathrm{NH}}(s)$ the fully discrete operators
that advance the numerical solutions of \eqref{eq:H} and \eqref{eq:NH}, respectively, over a time interval of length $s>0$.
They are the hydrostatic step of Section~\ref{sec:scheme:H} and the non-hydrostatic step of Sections~\ref{sec:scheme:stag} and~\ref{sec:scheme:UPU}.
A time step from $t^n$ to $t^{n+1}=t^n+\Delta t$ is then given by
\begin{equation}\label{eq:strang}
Q^{n+1}=\mathcal S_{\mathrm H}\!\left(\tfrac{\Delta t}{2}\right)\circ
\bigl[\mathcal S_{\mathrm{NH}}(\tau)\bigr]^{n_s}
\circ\mathcal S_{\mathrm H}\!\left(\tfrac{\Delta t}{2}\right)Q^n,\qquad \tau=\frac{\Delta t}{n_s},
\end{equation}
where
\[
\Delta t=\mathrm{CFL}\,\frac{\Delta x}{\max_{i,\alpha}\bigl(|u_{\alpha,i}|+\sqrt{gh_i}\bigr)},\qquad
n_s=\Bigl\lceil\frac{c\,\Delta t}{\mathrm{CFL}_{\mathrm{NH}}\,\Delta x}\Bigr\rceil .
\]
Therefore, the time step is determined by the hydrostatic part, through an estimate of the
speeds of Proposition~\ref{prop:subsystems}(i), while the $n_s$ sub-steps of the
non-hydrostatic part satisfy
\[
c\,\tau\le\mathrm{CFL}_{\mathrm{NH}}\,\Delta x .
\]
By Proposition~\ref{prop:subsystems}(ii), this is a CFL condition for \eqref{eq:NH}, whose characteristic speeds are bounded by $c$ for any state
and any number of layers. The stability of the discrete step is studied in Theorem~\ref{thm:stab}. Note that both $\Delta t$ and $\tau$ are proportional to $\Delta x$, so
that the number of sub-steps,
\[
n_s\approx\frac{\mathrm{CFL}}{\mathrm{CFL}_{\mathrm{NH}}}\,\frac{c}{\max_{i,\alpha}\bigl(|u_{\alpha,i}|+\sqrt{gh_i}\bigr)},
\]
depends neither on the number of layers (for a fixed $c$) nor on the mesh size: refining the mesh does not
increase the number of non-hydrostatic sub-steps per time step. We take $\mathrm{CFL}=0.45$ and
$\mathrm{CFL}_{\mathrm{NH}}=0.9$, in one and two dimensions.

Here, $\max_\alpha|u_\alpha|+\sqrt{gh}$ is used as an estimate of the largest characteristic speed of the
hydrostatic multilayer system in Proposition~\ref{prop:subsystems}(i). It is exact for a single layer and
for states without vertical shear, $u_\alpha=\ub$, but the eigenvalues of the multilayer system are not known
in closed form when the velocity varies along the vertical. The estimate is motivated by the continuous
model. The multilayer system is a vertical discretization of the hydrostatic Euler equations, and for
these equations the speeds $\lambda$ of the long surface waves on a flow with vertical profile $u(z)$, with
$z$ the height above the bottom, satisfy the condition of Burns \cite{Burns1953}
\[
g\int_0^h\frac{dz}{\bigl(u(z)-\lambda\bigr)^2}=1 .
\]
For the fast wave, $\lambda>\max_z u$, the left-hand side lies between $gh/(\lambda-\min_z u)^2$ and
$gh/(\lambda-\max_z u)^2$, so that
\[
\min_z u+\sqrt{gh}\le\lambda\le\max_z u+\sqrt{gh},
\]
and similarly for the slow wave. The same speeds are used in the HLL fluxes of the hydrostatic step
(\ref{app:hydro}).

The composition \eqref{eq:strang} is second-order accurate in time, since both $\mathcal S_{\mathrm H}$ and
$\mathcal S_{\mathrm{NH}}(\tau)$ are second-order accurate: the symmetric composition removes the leading
term of the splitting error, and the sub-steps do not need to be symmetric. It is well known that
splitting methods may lose accuracy when one of the operators is stiff and is integrated with a
large time step. This is not the case here: since the sub-steps satisfy $c\tau\le\Delta x$, the
fast horizontal dynamics of \eqref{eq:NH} is resolved in time, while the stiff vertical coupling is
treated implicitly in each sub-step. Note, however, that the order is not uniform in $c$, since
the error constants grow with $c$. For this reason, a fixed and moderate value of $c$ is always
used (see Section~\ref{sec:scheme:props}).

\subsection{Staggered mesh}\label{sec:scheme:mesh}
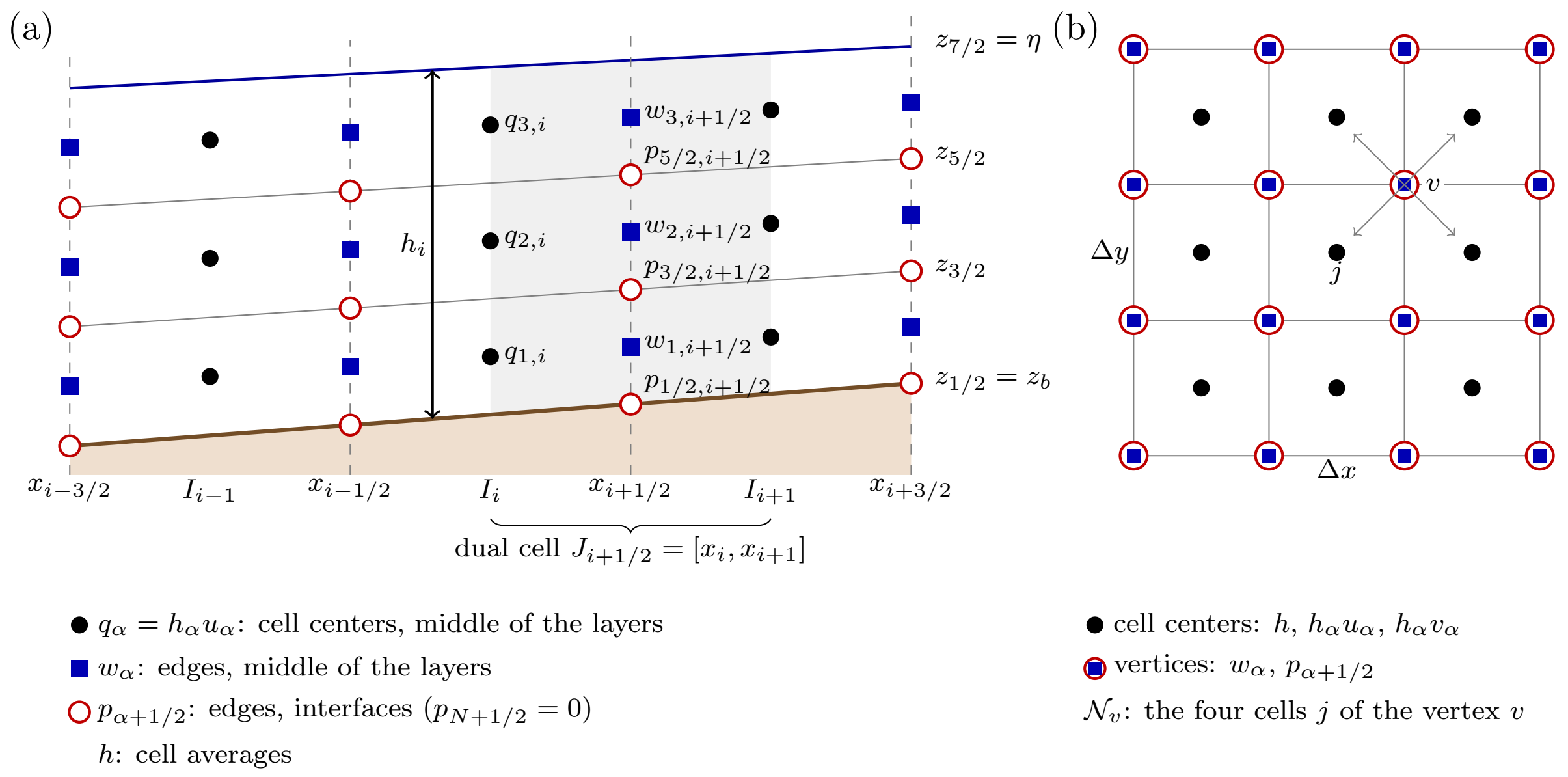
\begin{figure}[htbp]
\centering
\resizebox{\linewidth}{!}{\begin{tikzpicture}[x=1cm,y=1cm,
  hu/.style={circle,fill=black,inner sep=0pt,minimum size=5pt},
  ww/.style={rectangle,fill=blue!70!black,inner sep=0pt,minimum size=5pt},
  pp/.style={circle,draw=red!75!black,fill=white,thick,inner sep=0pt,minimum size=6pt},
  lab/.style={font=\scriptsize,inner sep=1pt}]
\def\NL{3}
\def\cw{2.9}   
\pgfmathsetmacro{\xa}{0}\pgfmathsetmacro{\xb}{\cw}\pgfmathsetmacro{\xc}{2*\cw}\pgfmathsetmacro{\xd}{3*\cw}
\pgfmathsetmacro{\xia}{0.5*\cw}\pgfmathsetmacro{\xib}{1.5*\cw}\pgfmathsetmacro{\xic}{2.5*\cw}
\pgfmathdeclarefunction{zb}{1}{\pgfmathparse{0.30+0.06*#1*3.6/\cw}}
\pgfmathdeclarefunction{eta}{1}{\pgfmathparse{4.00+0.04*#1*3.6/\cw}}
\pgfmathdeclarefunction{zi}{2}{\pgfmathparse{zb(#1)+(#2)/\NL*(eta(#1)-zb(#1))}}
\fill[gray!12] (\xib,{zb(\xib)}) -- (\xic,{zb(\xic)}) -- (\xic,{eta(\xic)}) -- (\xib,{eta(\xib)}) -- cycle;
\fill[brown!25] (0,{zb(0)}) -- (\xd,{zb(\xd)}) -- (\xd,0) -- (0,0) -- cycle;
\draw[very thick,brown!60!black] (0,{zb(0)}) -- (\xd,{zb(\xd)});
\foreach \a in {1,2} \draw[gray] (0,{zi(0,\a)}) -- (\xd,{zi(\xd,\a)});
\draw[thick,blue!60!black] (0,{eta(0)}) -- (\xd,{eta(\xd)});
\foreach \x in {\xa,\xb,\xc,\xd} \draw[dashed,gray] (\x,0) -- (\x,{eta(\x)+0.35});
\node[lab,anchor=west] at (\xd+0.2,{zi(\xd,0)}) {$z_{1/2}=z_b$};
\node[lab,anchor=west] at (\xd+0.2,{zi(\xd,1)}) {$z_{3/2}$};
\node[lab,anchor=west] at (\xd+0.2,{zi(\xd,2)}) {$z_{5/2}$};
\node[lab,anchor=west] at (\xd+0.2,{zi(\xd,3)}) {$z_{7/2}=\eta$};
\foreach \x in {\xia,\xib,\xic}
  \foreach \a in {1,...,\NL} \node[hu] at (\x,{zi(\x,\a-0.5)}) {};
\foreach \x in {\xa,\xb,\xc,\xd}{
  \foreach \a in {1,...,\NL} \node[ww] at (\x,{zi(\x,\a-0.5)}) {};
  \foreach \a in {0,1,2} \node[pp] at (\x,{zi(\x,\a)}) {};
}
\draw[<->,thick] (\xib-0.6,{zb(\xib-0.6)}) -- node[lab,left] {$h_i$} (\xib-0.6,{eta(\xib-0.6)});
\foreach \a in {1,...,\NL}
  \node[lab,anchor=west] at (\xib+0.1,{zi(\xib,\a-0.5)}) {$q_{\a,i}$};
\foreach \a in {1,...,\NL}
  \node[lab,anchor=west] at (\xc+0.1,{zi(\xc,\a-0.5)}) {$w_{\a,i+1/2}$};
\foreach \a/\t in {0/{1/2},1/{3/2},2/{5/2}}
  \node[lab,anchor=south west] at (\xc+0.1,{zi(\xc,\a)+0.02}) {$p_{\t,i+1/2}$};
\foreach \x/\t in {\xa/{x_{i-3/2}},\xb/{x_{i-1/2}},\xc/{x_{i+1/2}},\xd/{x_{i+3/2}}}
  \node[lab,below] at (\x,0) {$\t$};
\foreach \x/\t in {\xia/{I_{i-1}},\xib/{I_i},\xic/{I_{i+1}}}
  \node[lab,below] at (\x,0) {$\t$};
\draw[decorate,decoration={brace,mirror,amplitude=4pt}] (\xib,-0.45) -- node[lab,below=4pt]
  {dual cell $J_{i+1/2}=[x_i,x_{i+1}]$} (\xic,-0.45);
\node[font=\small] at (-0.4,4.6) {(a)};
\begin{scope}[shift={(0.1,-1.55)}]
  \node[hu] at (0,0) {}; \node[lab,anchor=west] at (0.15,0) {$q_\alpha=h_\alpha u_\alpha$: cell centers, middle of the layers};
  \node[ww] at (0,-0.45) {}; \node[lab,anchor=west] at (0.15,-0.45) {$w_\alpha$: edges, middle of the layers};
  \node[pp] at (0,-0.9) {}; \node[lab,anchor=west] at (0.15,-0.9) {$p_\aph$: edges, interfaces ($p_{N+1/2}=0$)};
  \node[lab,anchor=west] at (0.15,-1.35) {$h$: cell averages};
\end{scope}
\begin{scope}[shift={(11.0,0.2)}]
  \node[font=\small] at (-0.6,4.4) {(b)};
  \foreach \k in {0,1.4,2.8,4.2} {\draw[gray] (\k,0) -- (\k,4.2); \draw[gray] (0,\k) -- (4.2,\k);}
  \foreach \x in {0.7,2.1,3.5} \foreach \y in {0.7,2.1,3.5} \node[hu] at (\x,\y) {};
  \foreach \x in {0,1.4,2.8,4.2} \foreach \y in {0,1.4,2.8,4.2} {
    \node[pp,minimum size=8pt] at (\x,\y) {}; \node[ww,minimum size=4pt] at (\x,\y) {};}
  \node[lab,anchor=west,fill=white] at (2.98,2.8) {$v$};
  \foreach \x/\y in {2.1/2.1,3.5/2.1,2.1/3.5,3.5/3.5} \draw[->,gray] (2.8,2.8) -- ({0.75*\x+0.25*2.8},{0.75*\y+0.25*2.8});
  \node[lab,anchor=north] at (2.1,2.05) {$j$};
  \node[lab,below] at (2.1,0) {$\Delta x$};
  \node[lab,left] at (0,2.1) {$\Delta y$};
  \begin{scope}[shift={(-0.4,-1.75)}]
    \node[hu] at (0,0) {}; \node[lab,anchor=west] at (0.15,0) {cell centers: $h$, $h_\alpha u_\alpha$, $h_\alpha v_\alpha$};
    \node[pp] at (0,-0.45) {}; \node[ww,minimum size=4pt] at (0,-0.45) {}; \node[lab,anchor=west] at (0.15,-0.45) {vertices: $w_\alpha$, $p_\aph$};
    \node[lab,anchor=west] at (-0.15,-0.9) {$\mathcal N_v$: the four cells $j$ of the vertex $v$};
  \end{scope}
\end{scope}
\end{tikzpicture}
}
\caption{Location of the unknowns. (a) Vertical section in one dimension, for $N=3$ layers: the
momenta of the layers are located at the centers of the cells and at the middle of the layers,
the vertical velocities at the edges and at the middle of the layers, and the pressures at the
edges and at the interfaces between the layers. The vertical velocities and the pressures are
transported on the dual cells $J_{i+1/2}$ (shaded). (b) Plan view in two dimensions: the
vertical velocities and the pressures are located at the vertices of the mesh.}
\label{fig:stencil}
\end{figure}

Let us consider a uniform mesh of cells $I_i=[x_{i-1/2},x_{i+1/2}]$ of size $\Delta x$, and let us
introduce the notation
\[
\delta_{i+1/2}f=f_{i+1}-f_i,\qquad \{f\}_{i+1/2}=\frac12\left(f_i+f_{i+1}\right).
\]
When the edge is clear from the context, the subscript of $\delta_{i+1/2}$ will be dropped, so that $\delta f$ denotes $\delta_{i+1/2}f$ in the formulas written at the edge $x_{i+1/2}$.
The depth $h_i$ and the horizontal velocities $u_{\alpha,i}$ are cell averages, while the
vertical velocities $w_{\alpha,i+1/2}$ and the interface pressures $p_{\aph,i+1/2}$ are located at
the cell edges (staggered mesh; see Figure~\ref{fig:stencil}), where we define
\[
h_{i+1/2}=\{h\}_{i+1/2},\qquad h_{\alpha,i+1/2}=\frac{h_{i+1/2}}N,\qquad
\hs_{\aph,i+1/2}=\sigma_\alpha h_{\alpha,i+1/2}.
\]
Staggered meshes have also been used to build energy-stable schemes for the shallow water equations
\cite{DuranVilaBaraille2020}.

\subsection{Hydrostatic step}\label{sec:scheme:H}
The hydrostatic multilayer system is discretized by a well-balanced path-conservative finite
volume method of HLL type \cite{CastroMoralesPares2017}, combined with the SSP-RK2 method in
time. At second order, a MUSCL reconstruction with the monotonized central (MC) limiter \cite{vanLeer1977} is used,
together with the hydrostatic reconstruction of \cite{Audusse2004}. The non-conservative products
are the exchange terms and $gh_\alpha\dx\eta$. Two details matter when the layers move at different
velocities: the HLL speeds take into account the velocities of all the layers, and the exchange terms
use the same mass fluxes of the layers as the transport of momentum, also in the hydrostatic
reconstruction (see \ref{app:hydro}).
The values located at the edges are transported in non-conservative form,
\[
\begin{aligned}
&h_\alpha\bigl(\dt w_\alpha+u_\alpha\dx w_\alpha\bigr)
=-\frac12(w_{\alpha+1}-w_\alpha)\Gamma_\aph+\frac12(w_{\alpha-1}-w_\alpha)\Gamma_\amh,\\
&\hs_\aph\bigl(\dt p_\aph+\hat u_\aph\dx p_\aph\bigr)
=-(p_{\alpha+1}-p_\aph)\Gb_{\alpha+1}+(p_\alpha-p_\aph)\Gb_\alpha,
\end{aligned}
\]
with
\[
\hat u_\aph=\frac{\qs_\aph}{\hs_\aph}.
\]
These equations follow from \eqref{eq:layermass}, \eqref{eq:stagmass} and the conservative form
of \eqref{eq:relax}, and they are discretized on the dual mesh by a second-order upwind finite
difference scheme with MC-limited slopes, in which the upwind direction is given by the sign of
the mean of the transport velocity in the two neighboring cells, $\Gamma$ being computed at the
edges by compact differences. The details of the discretization of the hydrostatic step are
given in \ref{app:hydro}.

\subsection{Spatial discretization of the non-hydrostatic part}\label{sec:scheme:stag}
At the edge $x_{i+1/2}$, the equivalent constraint \eqref{eq:K} is discretized by compact
differences:
\begin{equation}\label{eq:Kstag}
\mathcal K_{\aph,i+1/2}=\frac{h_{\alpha,i+1/2}\,\delta u_\alpha+h_{\alpha+1,i+1/2}\,\delta u_{\alpha+1}}{2\Delta x}
+\bigl(w_{\alpha+1}-w_\alpha\bigr)_{i+1/2}
-\{u_{\alpha+1}-u_\alpha\}_{i+1/2}\,\frac{\delta z_\aph}{\Delta x},
\end{equation}
for $\alpha=0,\dots,N-1$, with $u_0=w_0=0$, where $\delta=\delta_{i+1/2}$ as announced in Section~\ref{sec:scheme:mesh}. The jumps of the
interfaces are computed by means of the hydrostatic reconstruction of the free surface,
\begin{equation}\label{eq:dzHR}
\delta_{i+1/2}z_\aph=\delta^*_{i+1/2}\eta-\Bigl(1-\frac{\alpha}{N}\Bigr)\delta_{i+1/2}h,\qquad
\delta^*_{i+1/2}\eta=\max(\eta_{i+1},b_{i+1/2})-\max(\eta_i,b_{i+1/2}),
\end{equation}
where
\[
b_{i+1/2}=\max(z_{b,i},z_{b,i+1}).
\]
Note that, where both cells are wet, \eqref{eq:dzHR} is the usual jump of
\[
z_\aph=\eta-\Bigl(1-\frac{\alpha}{N}\Bigr)h,
\]
while at a wet-dry front it is bounded by the depth (see \ref{app:wetdry}).

The discretization is based on the following idea. Let us write the discrete constraint as
\[
\mathcal K=B_u\mathbf u+B_w\mathbf w ;
\]
then the discrete pressure force is \emph{defined} as its transpose,
\begin{equation}\label{eq:forcestag}
\mathbb M_c\,\frac{d\mathbf u}{dt}=B_u^T\mathbf p,\qquad
\mathbb M_e\,\frac{d\mathbf w}{dt}=B_w^T\mathbf p,\qquad
\mathbb M_{\hat h}\,\frac{d\mathbf p}{dt}=-c^2\bigl(B_u\mathbf u+B_w\mathbf w\bigr),
\end{equation}
where $\mathbb M_c$, $\mathbb M_e$ and $\mathbb M_{\hat h}$ denote the diagonal matrices of the
thicknesses $h_{\alpha,i}$, $h_{\alpha,i+1/2}$ and $\hs_{\aph,i+1/2}$. It will also be convenient
to write
\[
U=(\mathbf u,\mathbf w),\qquad B=\bigl(B_u\;\;B_w\bigr),\qquad
\mathbb M=\mathrm{diag}\bigl(\mathbb M_c,\mathbb M_e\bigr),
\]
so that $\mathcal K=B\,U$, and the force on the velocities is $\mathbb M^{-1}B^T\mathbf p$. More
explicitly,
\[
(B_w^T\mathbf p)_{\alpha,i+1/2}=-\left(p_\aph-p_\amh\right)_{i+1/2}
\]
and
\[
(B_u^T\mathbf p)_{\alpha,i}=-\frac{(h_\alpha p_\alpha)_{i+1/2}-(h_\alpha p_\alpha)_{i-1/2}}{\Delta x}
+\Bigl\{p_\aph\,\frac{\delta z_\aph}{\Delta x}\Bigr\}_i-\Bigl\{p_\amh\,\frac{\delta z_\amh}{\Delta x}\Bigr\}_i,
\]
where $p_\alpha=\frac12(p_\aph+p_\amh)$ at the edges and
\[
\{f\}_i=\frac12\left(f_{i-1/2}+f_{i+1/2}\right).
\]
This is a compact discretization of the force $\mathcal D_{1,\alpha}$ in \eqref{eq:D12}.
The discretization in space of the non-hydrostatic part preserves the energy exactly: the force on the velocities is the transpose of the discrete constraint, as in the continuous model. By
construction, the semi-discrete system inherits the energy conservation of the continuous
non-hydrostatic part: the energy $E$ below is the discrete version of the energy of \eqref{eq:NH} in
Section~\ref{sec:split:time}, in which the thicknesses $h_\alpha$ and $\hs_\aph$ are replaced by
the diagonal matrices $\mathbb M_c$, $\mathbb M_e$ and $\mathbb M_{\hat h}$.

\begin{proposition}[Semi-discrete energy]\label{prop:stagenergy}
For periodic data and $h$ frozen, the solutions of \eqref{eq:forcestag} conserve exactly the energy
\[
E=\frac12\mathbf u^T\mathbb M_c\mathbf u+\frac12\mathbf w^T\mathbb M_e\mathbf w
+\frac1{2c^2}\mathbf p^T\mathbb M_{\hat h}\mathbf p .
\]
\end{proposition}
\begin{proof}
Since $h$ is frozen, the matrices $\mathbb M_c$, $\mathbb M_e$ and $\mathbb M_{\hat h}$ do not
depend on time, and \eqref{eq:forcestag} gives
\[
\frac{d}{dt}\Bigl(\frac12\mathbf u^T\mathbb M_c\mathbf u+\frac12\mathbf w^T\mathbb M_e\mathbf w\Bigr)
=\mathbf u^TB_u^T\mathbf p+\mathbf w^TB_w^T\mathbf p=\mathbf p^T\mathcal K,
\]
while
\[
\frac{d}{dt}\Bigl(\frac1{2c^2}\mathbf p^T\mathbb M_{\hat h}\mathbf p\Bigr)
=\frac1{c^2}\,\mathbf p^T\mathbb M_{\hat h}\frac{d\mathbf p}{dt}=-\mathbf p^T\mathcal K .
\]
\end{proof}

Moreover, the balance is local. For a set of consecutive cells and the edges between them, the energy of the
unknowns in the set changes only through the work of the pressure at its two ends.

Proposition~\ref{prop:stagenergy} is the discrete counterpart of Proposition~\ref{prop:work}, and
therefore no numerical viscosity is needed in the non-hydrostatic part. Moreover, the same pair
$(\mathcal K,B^T)$ is used to build a projection method for LDNH$_0$: the velocities are corrected
by $\mathbb M^{-1}B^T\boldsymbol\psi$, where the impulse $\boldsymbol\psi$ is the solution of the
block tridiagonal system
\[
B\,\mathbb M^{-1}B^T\boldsymbol\psi=-\mathcal K .
\]
This projection method provides the reference solutions of LDNH$_0$ used in
Section~\ref{sec:tests} (see \ref{app:proj}).

\subsection{Time integration of the non-hydrostatic part}\label{sec:scheme:UPU}
In \eqref{eq:forcestag}, the horizontal velocities couple neighboring edges, while $\mathbf w$ and
$\mathbf p$ only couple the unknowns located at the same edge. This suggests a step that is
explicit in $\mathbf u$ and implicit in $(\mathbf w,\mathbf p)$, which consists of a half update of
the horizontal velocities, the implicit midpoint (Crank--Nicolson) rule for $(\mathbf w,\mathbf p)$
with $\mathbf u$ fixed, and a second half update of the velocities:
\begin{equation}\label{eq:UPU}
\begin{aligned}
&\mathbf u^{*}=\mathbf u^n+\tfrac\tau2\,\mathbb M_c^{-1}B_u^T\mathbf p^n,\\
&\mathbb M_e\frac{\mathbf w^{n+1}-\mathbf w^n}{\tau}=B_w^T\mathbf p^{n+1/2},\qquad
\mathbb M_{\hat h}\frac{\mathbf p^{n+1}-\mathbf p^n}{\tau}=-c^2\bigl(B_u\mathbf u^{*}+B_w\mathbf w^{n+1/2}\bigr),\\
&\mathbf u^{n+1}=\mathbf u^{*}+\tfrac\tau2\,\mathbb M_c^{-1}B_u^T\mathbf p^{n+1},
\end{aligned}
\end{equation}
where, for $f=\mathbf w,\mathbf p$,
\[
f^{n+1/2}=\frac12\left(f^n+f^{n+1}\right).
\]
In what follows, this step will be referred to as the U--P--U step. Since $B_u\mathbf u^{*}$ is
known, the second line decouples into one small system per edge. Indeed, eliminating
$\mathbf w^{n+1/2}$, the midpoint pressures at the edge $x_{i+1/2}$ are the solution of the
tridiagonal system
\begin{equation}\label{eq:CNsys}
-\frac{\nu}{\sigma_\alpha}\,p^{n+1/2}_\amh+(1+2\nu)\,p^{n+1/2}_\aph
-\frac{\nu}{\sigma_\alpha}\,p^{n+1/2}_{\alpha+3/2}
=p^n_\aph-\frac{c^2\tau}{2\sigma_\alpha h_e}\Bigl[(B_u\mathbf u^{*})_\aph+w^n_{\alpha+1}-w^n_\alpha\Bigr],
\end{equation}
for $\alpha=0,\dots,N-1$, with $p^{n+1/2}_{-1/2}=p^{n+1/2}_{N+1/2}=0$, where
$h_e=h_{\alpha,i+1/2}$ denotes the thickness of the layers at the edge and
\[
\nu=\frac{c^2\tau^2}{4h_e^2}.
\]
Then, the pressures and the vertical velocities at the new time level are given by
\[
p^{n+1}_\aph=2p^{n+1/2}_\aph-p^n_\aph,\qquad
w^{n+1}_\alpha=w^n_\alpha-\frac{\tau}{h_e}\bigl(p^{n+1/2}_\aph-p^{n+1/2}_\amh\bigr).
\]
After scaling its rows by $\sigma_\alpha$, the matrix of \eqref{eq:CNsys} is symmetric and
strictly diagonally dominant, so that the system can be solved by the Thomas algorithm in $O(N)$
operations.

Note that the horizontal and the vertical couplings are not solved in separate steps: each
sub-step requires two explicit updates of the velocities and one tridiagonal solve of size $N$
per edge. The step is second-order accurate and symmetric. It is the leapfrog (St\"ormer--Verlet) method~\cite{HairerLubichWanner2003} for the semi-discrete system \eqref{eq:forcestag}, explicit in the horizontal direction, and its stability follows from a
discrete energy estimate, as stated in the following result.

\begin{theorem}[Stability of the U--P--U step]\label{thm:stab}
Let $h$ be frozen, let the data be periodic, let $\lambda_{\max}$ be the largest eigenvalue of
$\mathbb M_{\hat h}^{-1}B_u\mathbb M_c^{-1}B_u^T$, and let $E$ denote the energy of
Proposition~\ref{prop:stagenergy}.

(i) Let $\mathbf v^{n\pm1/2}=\mathbf u^n\pm\frac\tau2\,\mathbb M_c^{-1}B_u^T\mathbf p^n$ be the velocities after the half updates before and after $t^n$. The step \eqref{eq:UPU} conserves exactly the energy
\[
E_*^n=\frac12(\mathbf v^{n-1/2})^T\mathbb M_c\,\mathbf v^{n+1/2}+\frac12(\mathbf w^n)^T\mathbb M_e\mathbf w^n
+\frac1{2c^2}(\mathbf p^n)^T\mathbb M_{\hat h}\mathbf p^n ,
\]
in which the kinetic energy of the horizontal velocities is computed with the product of the velocities before and after $t^n$.

(ii) At $t^n$, $E_*^n=E_\tau(\mathbf u^n,\mathbf w^n,\mathbf p^n)$, with
\[
E_\tau(\mathbf u,\mathbf w,\mathbf p)=E(\mathbf u,\mathbf w,\mathbf p)-\frac{\tau^2}{8}\,(B_u^T\mathbf p)^T\mathbb M_c^{-1}(B_u^T\mathbf p),
\]
that is, $E_*$ is the energy $E$ minus a term of order $\tau^2$ and of fixed sign. This is the modified energy of second order of the leapfrog method~\cite{SkeelHardy2001}. Moreover, $E_\tau$ is positive definite if and only if
\[
c\,\tau\sqrt{\lambda_{\max}}<2,
\]
and in that case
\[
\Bigl(1-\frac{c^2\tau^2\lambda_{\max}}{4}\Bigr)E\le E_\tau\le E .
\]
Therefore, the iterates of the step satisfy
\[
E^n\le\Bigl(1-\frac{c^2\tau^2\lambda_{\max}}{4}\Bigr)^{-1}E^0\qquad\text{for all }n\ge0,
\]
that is, the step is stable and does not dissipate.

(iii) For a uniform state over a flat bottom,
\[
\lambda_{\max}=\frac{4}{\Delta x^2}\cos^2\Bigl(\frac{\pi}{4N}\Bigr),
\]
so that the condition reads
\[
c\,\tau<\frac{\Delta x}{\cos\bigl(\pi/(4N)\bigr)},
\]
and $c\tau\le\Delta x$ is enough for any $N$.
\end{theorem}
\begin{proof}
Let us denote by $\mathbf v^{n+1/2}=\mathbf u^{*}$ the velocity after the first half update.
Since the second half update of a step and the first half update of the next one use the same
pressure, we have
\[
\mathbf v^{n+1/2}-\mathbf v^{n-1/2}=\tau\,\mathbb M_c^{-1}B_u^T\mathbf p^n,
\qquad
\mathbf v^{n\pm1/2}=\mathbf u^n\pm\frac\tau2\,\mathbb M_c^{-1}B_u^T\mathbf p^n .
\]
We next show that
\[
E_*^n=\frac12(\mathbf v^{n-1/2})^T\mathbb M_c\,\mathbf v^{n+1/2}+\frac12(\mathbf w^n)^T\mathbb M_e\mathbf w^n
+\frac1{2c^2}(\mathbf p^n)^T\mathbb M_{\hat h}\mathbf p^n
\]
is conserved. Multiplying the Crank--Nicolson equations by $\mathbf w^{n+1/2}$ and by
$\mathbf p^{n+1/2}/c^2$, respectively, and adding them, we obtain
\[
\frac12\bigl|\mathbf w^{n+1}\bigr|^2_{\mathbb M_e}-\frac12\bigl|\mathbf w^n\bigr|^2_{\mathbb M_e}
+\frac{1}{2c^2}\bigl|\mathbf p^{n+1}\bigr|^2_{\mathbb M_{\hat h}}-\frac{1}{2c^2}\bigl|\mathbf p^n\bigr|^2_{\mathbb M_{\hat h}}
=-\tau\,(\mathbf p^{n+1/2})^TB_u\mathbf v^{n+1/2},
\]
since the terms involving $B_w$ cancel out. On the other hand,
\[
\frac12(\mathbf v^{n+1/2})^T\mathbb M_c\bigl(\mathbf v^{n+3/2}-\mathbf v^{n-1/2}\bigr)
=\frac\tau2(\mathbf v^{n+1/2})^TB_u^T\bigl(\mathbf p^{n}+\mathbf p^{n+1}\bigr)
=\tau\,(\mathbf p^{n+1/2})^TB_u\mathbf v^{n+1/2}.
\]
Since $\mathbb M_c$ is symmetric, combining both identities we obtain
\[
E_*^{n+1}-E_*^n
=\frac12(\mathbf v^{n+1/2})^T\mathbb M_c\bigl(\mathbf v^{n+3/2}-\mathbf v^{n-1/2}\bigr)
-\tau\,(\mathbf p^{n+1/2})^TB_u\mathbf v^{n+1/2}=0 .
\]
Moreover, replacing $\mathbf v^{n\pm1/2}$ by their expressions in terms of $\mathbf u^n$, and
using the symmetry of $\mathbb M_c$, we have
\[
\frac12(\mathbf v^{n-1/2})^T\mathbb M_c\,\mathbf v^{n+1/2}
=\frac12(\mathbf u^n)^T\mathbb M_c\mathbf u^n
-\frac{\tau^2}{8}\,(B_u^T\mathbf p^n)^T\mathbb M_c^{-1}(B_u^T\mathbf p^n),
\]
so that, with the energy $E$ of Proposition~\ref{prop:stagenergy},
\[
E_*^n=E(\mathbf u^n,\mathbf w^n,\mathbf p^n)-\frac{\tau^2}{8}\,\bigl|B_u^T\mathbf p^n\bigr|^2_{\mathbb M_c^{-1}}
=E_\tau(\mathbf u^n,\mathbf w^n,\mathbf p^n),
\]
which proves (i) and the first statement of (ii).

Let us now prove (ii). Since $E_\tau=E-\frac{\tau^2}{8}|B_u^T\mathbf p|^2_{\mathbb M_c^{-1}}$,
the modified energy is positive definite if and only if the quadratic form in $\mathbf p$
\[
\frac{1}{2c^2}\,\mathbf p^T\mathbb M_{\hat h}\mathbf p-\frac{\tau^2}{8}\,\mathbf p^TB_u\mathbb M_c^{-1}B_u^T\mathbf p
\]
is positive definite, that is, if and only if $\lambda_{\max}<4/(c^2\tau^2)$. Moreover, by
definition of $\lambda_{\max}$,
\[
\frac{\tau^2}{8}\,\bigl|B_u^T\mathbf p\bigr|^2_{\mathbb M_c^{-1}}
\le\frac{\tau^2\lambda_{\max}}{8}\,|\mathbf p|^2_{\mathbb M_{\hat h}}
=\frac{c^2\tau^2\lambda_{\max}}{4}\cdot\frac{1}{2c^2}\,|\mathbf p|^2_{\mathbb M_{\hat h}}
\le\frac{c^2\tau^2\lambda_{\max}}{4}\,E,
\]
which gives the lower bound for $E_\tau$; the upper bound is obvious. Finally, since $E_\tau$ is
conserved,
\[
\Bigl(1-\frac{c^2\tau^2\lambda_{\max}}{4}\Bigr)E^n\le E_\tau^n=E_\tau^0\le E^0 .
\]

Let us finally prove (iii). For a uniform state over a flat bottom, the geometric terms vanish
and all the thicknesses are equal to $h/N$, so that
\[
\mathbb M_c=\frac hN\,I,\qquad \mathbb M_{\hat h}=\frac hN\,\Sigma,
\]
with $\Sigma=\mathrm{diag}(\sigma_0,\dots,\sigma_{N-1})$. Let us consider the Fourier
modes
\[
u_{\alpha,i}=\hat u_\alpha\,e^{\mathrm ikx_i},\qquad p_{\aph,i+1/2}=\hat p_\aph\,e^{\mathrm ikx_{i+1/2}} .
\]
Since $e^{\mathrm ikx_{i+1}}-e^{\mathrm ikx_i}=2\mathrm i\sin(k\Delta x/2)\,e^{\mathrm ikx_{i+1/2}}$, and denoting
\[
\varsigma=\frac{2\mathrm i}{\Delta x}\sin\Bigl(\frac{k\Delta x}{2}\Bigr),
\]
the definitions of $B_u$ and of its transpose give
\[
(B_u\mathbf u)_\aph=\frac hN\,\varsigma\,\frac{\hat u_\alpha+\hat u_{\alpha+1}}2\,e^{\mathrm ikx_{i+1/2}},
\qquad
(B_u^T\mathbf p)_\alpha=-\frac hN\,\varsigma\,\frac{\hat p_\aph+\hat p_\amh}2\,e^{\mathrm ikx_i},
\]
with $\hat u_0=0$ and $\hat p_{N+1/2}=0$. In terms of the matrix $P$ of \eqref{eq:lin}, the
coefficients of these two expressions are $\frac hN\varsigma P^T\hat{\mathbf u}$ and
$-\frac hN\varsigma P\hat{\mathbf p}$. Therefore, since
$-\varsigma^2=\frac4{\Delta x^2}\sin^2(k\Delta x/2)$, the matrix
$\mathbb M_{\hat h}^{-1}B_u\mathbb M_c^{-1}B_u^T$ acts on the Fourier modes as
\[
\frac{4}{\Delta x^2}\sin^2\Bigl(\frac{k\Delta x}{2}\Bigr)\,\Sigma^{-1}P^TP .
\]
Since the largest eigenvalue of $\Sigma^{-1}P^TP$ is $\cos^2(\pi/(4N))$ (Lemma~\ref{lem:modes}), the
maximum is reached for
$k\Delta x=\pi$, which proves (iii).
\end{proof}

\begin{remark}\label{rem:taucross}
The term $\frac{\tau^2}{8}|B_u^T\mathbf p|^2_{\mathbb M_c^{-1}}$ in $E_\tau$ is not dissipation. The
pressure $\mathbf p^n$ is used in the second half update of step $n-1$ and in the first half update of
step $n$, so that the energy conserved by the step uses the product of the velocities at both sides of
$\mathbf p^n$, $\mathbf v^{n\pm1/2}=\mathbf u^n\pm\frac\tau2\mathbb M_c^{-1}B_u^T\mathbf p^n$, and not the square of the
saved velocity $\mathbf u^n$. Therefore $E$ stays above $E_\tau$, which is constant, with a difference of order
$\tau^2$, and it does not drift (see Section~\ref{sec:tests:energy}).
\end{remark}

Some comments are in order. Since $B_u$ is a discrete divergence, $\lambda_{\max}$ scales like
$\Delta x^{-2}$, but the stability condition involves $\sqrt{\lambda_{\max}}$ and is a standard CFL
condition $\tau=O(\Delta x)$ with speed $c$, in agreement with
Proposition~\ref{prop:subsystems}(ii); therefore the number of sub-steps does not depend on the mesh
size. The vertical coupling, which is stiff (Proposition~\ref{prop:stiff}), does not appear in the
condition, since it is treated implicitly. In general, $\lambda_{\max}$ depends on the state: the
geometric terms and the different depths at both sides of an edge add contributions of order
\[
\Bigl(\frac{N|\delta z|}{h\,\Delta x}\Bigr)^2
\qquad\text{and}\qquad
\frac{h_{i+1/2}}{h_i\,\Delta x^2},
\]
respectively. The first one acts on the explicit step as a stiff source of frequency
$cN|\dx z|/h$. It grows with $N$, so that the independence of the condition with respect to the
number of layers only holds exactly for uniform states; it is small for mild slopes and moderate
depths, and it becomes large near the shoreline, where $h\to0$.

Both contributions can be bounded by limiting the geometry in the coefficients of $B_u$. The
proof of Theorem~\ref{thm:stab} only uses that the same matrix $B_u$ appears in the constraint and
in the force \eqref{eq:forcestag}, so any change of its coefficients keeps the energy structure.
In \eqref{eq:Kstag}, and therefore in \eqref{eq:forcestag}, we replace the thickness
$h_{\alpha,i+1/2}=h_{i+1/2}/N$ and the jumps of the interfaces by
\begin{equation}\label{eq:limB}
\begin{aligned}
&h^{\mathrm{lim}}_{i+1/2}=\frac1N\min\bigl(h_{i+1/2},(1+\varepsilon_{\mathrm h})\min(h_i,h_{i+1})\bigr),\\
&\delta z^{\mathrm{lim}}_\aph=\max\bigl(-h^{\mathrm{lim}}_{i+1/2},\min(\delta z_\aph,h^{\mathrm{lim}}_{i+1/2})\bigr),
\end{aligned}
\end{equation}
with $\varepsilon_{\mathrm h}\ge0$ (we take $\varepsilon_{\mathrm h}=0.2$), while the matrices $\mathbb M_c$, $\mathbb M_e$ and $\mathbb M_{\hat h}$
are not changed. The coefficients are only modified where the depths of two neighboring cells
differ by more than a factor $1+2\varepsilon_{\mathrm h}$, or where an interface jumps by more than the
thickness of a layer. For smooth solutions, neither happens when the mesh is fine enough, so the
consistency of the scheme is not affected.

\subsection{Properties of the scheme}\label{sec:scheme:props}
Let us now make precise what is proved about the stability of the complete scheme. The
non-hydrostatic step $\mathcal S_{\mathrm{NH}}$, in which $h$ is frozen, is stable under the
condition of Theorem~\ref{thm:stab}, which for uniform states reads $c\tau\le\Delta x$; this is the
discrete counterpart of the CFL condition given by Proposition~\ref{prop:subsystems}(ii). The hydrostatic step $\mathcal S_{\mathrm H}$ is a standard
finite volume scheme, stable under the CFL condition of Section~\ref{sec:split:time}. However, the
stability of the composition \eqref{eq:strang}, in which $h$ changes between two consecutive
applications of $\mathcal S_{\mathrm{NH}}$, is not proved; it is observed in all our tests.

The lake at rest is preserved, since $\mathcal S_{\mathrm{NH}}$ is linear in $(u,w,p)$ and does not
modify it, and $\mathcal S_{\mathrm H}$ preserves it thanks to the hydrostatic reconstruction. The step
$\mathcal S_{\mathrm{NH}}$ has no numerical dissipation: it conserves exactly the energy $E_*$ of Theorem~\ref{thm:stab}, which differs from $E$ by a term of order $\tau^2$. The energy loss of
the scheme is the numerical diffusion of the upwind step $\mathcal S_{\mathrm H}$ (see
Section~\ref{sec:tests:energy}). The only implicit operation, apart from the initial data, is the
solution of the tridiagonal systems \eqref{eq:CNsys}, once per edge and per sub-step. The number of
sub-steps depends neither on $N$ nor on $\Delta x$, while the cost of each sub-step grows linearly
with $N$. Finally, the scheme is not asymptotic-preserving when $c\to\infty$, since the number of
sub-steps grows linearly with $c$ and the splitting error between \eqref{eq:H} and \eqref{eq:NH} grows with
$c\Delta t$ (see Section~\ref{sec:tests:soliton}); this is the reason why moderate values of $c$
are used.

Two further ingredients of the scheme are needed in practice: well-prepared initial data and a treatment of the wet-dry fronts.

\paragraph{Initial data}\label{sec:scheme:init}
The relaxed system has a fast branch of plane waves (Section~\ref{sec:relax:lin}), with frequencies of order $cN/H$
(Proposition~\ref{prop:stiff}), which does not exist in LDNH$_0$. Initial data satisfying
$\mathcal K=0$ and $p=0$ do not satisfy $\dt\mathcal K=0$, and they excite this branch with
oscillations whose amplitude is similar to that of the pressure itself, for any value of $c$.
Since the scheme has no numerical viscosity in the non-hydrostatic part, these oscillations are
not damped. They are avoided by using well-prepared initial data, that is, by doing at
$t=0$ what a step of the projection method does: the initial state is projected so that
$\mathcal K=0$, which fixes $\mathbf u$ and $\mathbf w$, and $p(0)$ is taken as the pressure of
LDNH$_0$, the Lagrange multiplier that keeps $\mathcal K=0$ in time ($\dt\mathcal K=0$). The second
choice is necessary because $p$ does not appear in $\mathcal K$: with $p=0$, or with any other
pressure, the relaxed system starts away from the slow solution. The pressure is computed with an auxiliary
projection step, which is discarded, and the first time step is already a step of the relaxed scheme.
This requires two solutions of the projection problem at $t=0$ (see \ref{app:init}). The effect of this choice is shown in Section~\ref{sec:tests:prepared}.

\paragraph{Wet-dry fronts}\label{sec:scheme:wd}
In the hydrostatic step, the treatment of the single-layer case \cite{Escalante2019JCP} is
followed: desingularized velocities \cite{KurganovPetrova2007}, reconstruction of the velocities
and zero momenta in dry cells (see also \cite{Castro2005}). In the non-hydrostatic part, the
interface jumps are computed with the hydrostatic reconstruction \eqref{eq:dzHR}, and the
coefficients of the constraint are limited as in \eqref{eq:limB}. The
non-hydrostatic pressure is never switched off. At the dry edges, $h_{i+1/2}<h_{\mathrm{dry}}$,
we set $w=p=0$. Finally, the layers are collapsed to a single layer where the water is very
shallow, $h<2\Delta x$: in these cells, the velocities of all the layers are replaced by their
mean (see \ref{app:wetdry}).

\section{Numerical tests in one dimension}\label{sec:tests}
Unless stated otherwise, the tests use the scheme of Section~\ref{sec:scheme} with second-order
reconstruction (MC limiter \cite{vanLeer1977}), $\mathrm{CFL}=0.45$, $\mathrm{CFL}_{\mathrm{NH}}=0.9$, and well-prepared initial
data. As usual in the literature, the still-water depth of each test is denoted by $H_0$, that
is, $H_0=H$ in the notation of the previous sections. In the tests without laboratory data, $H_0=1$\,m and $g=9.81$\,m\,s$^{-2}$, so that the lengths are in meters and the times in seconds; the units are not repeated. When no exact solution is available, the
reference solutions of LDNH$_0$ are computed with the projection method of
Section~\ref{sec:scheme:stag} on the same mesh, with the LSDIRK2 time integrator
\cite{BoscheriPareschi2021} (\ref{app:proj}).

We first compute exact solitary waves of LDNH$_0$ and of the relaxed system
(Section~\ref{sec:tests:tw}). They are used as reference solutions to measure the order of
convergence (Section~\ref{sec:tests:order}) and the error of the relaxation as a function of $c$
(Section~\ref{sec:tests:soliton}). Then, we study the evolution of the energy for the fully
discrete scheme (Section~\ref{sec:tests:energy}) and the influence of the initial data, in
particular of the use of well-prepared initial data (Section~\ref{sec:tests:prepared}). Finally,
the scheme is validated with two laboratory experiments: the run-up of a solitary wave on a plane
beach, which involves wet-dry fronts (Section~\ref{sec:tests:wd}), and periodic waves over a
submerged bar, in which dispersive effects are dominant (Section~\ref{sec:tests:bb}).

\subsection{Exact solitary waves and comparison with the Euler equations}\label{sec:tests:tw}
\begin{figure}[htbp]
\centering
\includegraphics[width=0.8\textwidth]{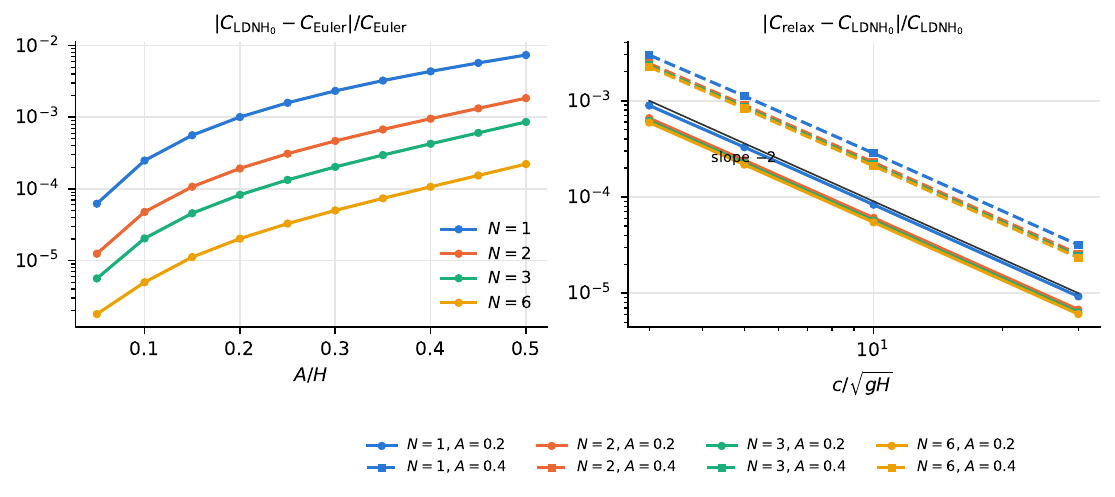}
\caption{Left: relative difference between the speeds of the solitary waves of LDNH$_0$ and of the
Euler equations versus the amplitude, for $N=1,2,3,6$. Right: relative difference between the speed of
the relaxed system and that of LDNH$_0$ versus $c$, for the same values of $N$ (same colors) and for
$A=0.2$ (solid lines) and $A=0.4$ (dashed lines); the thin line has slope $-2$.}
\label{fig:tw:speed}
\end{figure}

\begin{figure}[htbp]
\centering
\includegraphics[width=0.85\textwidth]{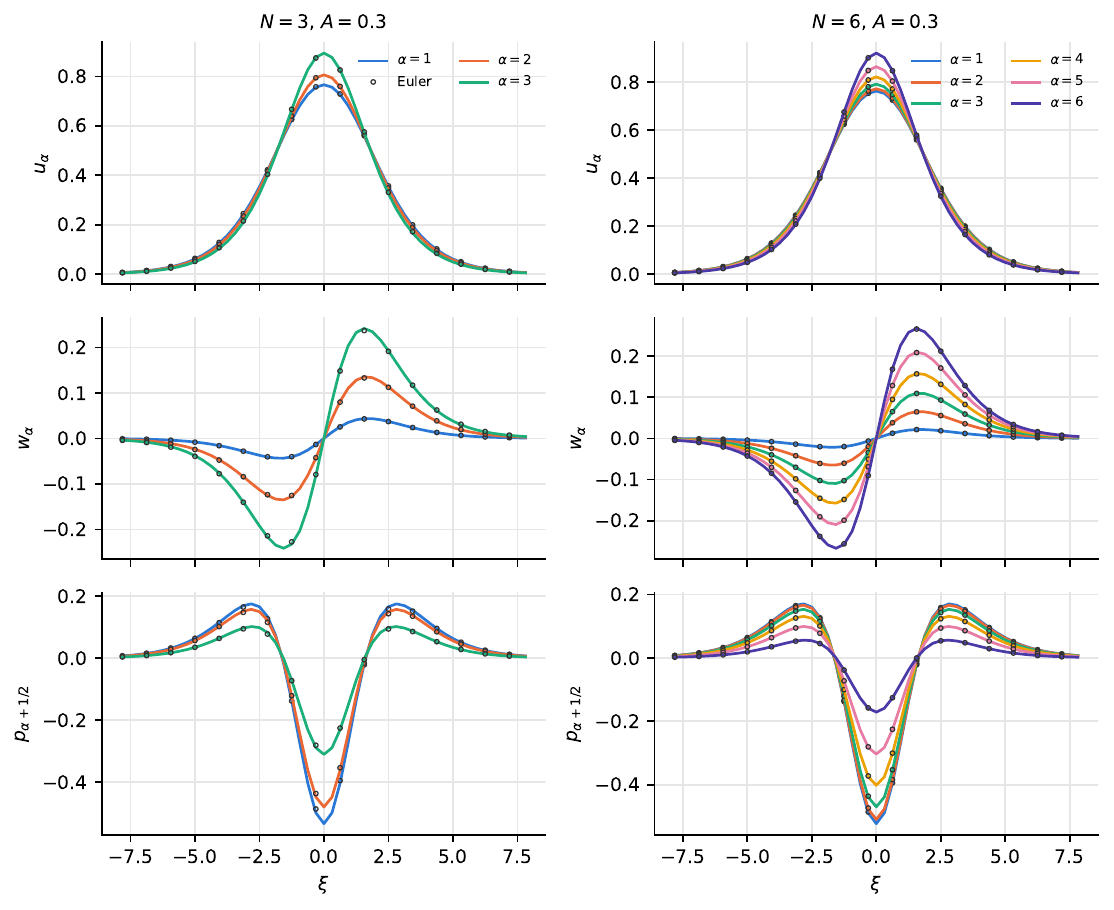}
\caption{Vertical structure of the solitary wave $A=0.3$ in the moving coordinate
$\xi=x-Ct$: horizontal velocity $u_\alpha$ and vertical velocity $w_\alpha$ of each layer, and pressure
$p_{\alpha+1/2}$ at the interfaces, for LDNH$_0$ with $N=3$ and $N=6$ (lines, one color for each
layer, with $\alpha=1$ at the bottom) and layer means of the Euler solution
(circles).}
\label{fig:tw:vertical}
\end{figure}

We compute exact solitary waves of amplitude $A$ of LDNH$_0$ and of the relaxed system with a spectral Newton
method, and solitary waves of the Euler equations with the method of \cite{ClamondDutykh2013}
(\ref{app:tw}). The speed of the solitary waves of LDNH$_0$ converges to that of the
Euler equations when $N$ grows (Figure~\ref{fig:tw:speed}). The
relative error is of order $10^{-3}$ for $N=1$ and $2\cdot10^{-5}$ for $N=6$ at $A=0.2$, and it
decreases roughly like $N^{-2}$. The relaxed system converges to LDNH$_0$ like $c^{-2}$ for every
$N$ (Table~\ref{tab:tw:relax}). The vertical structure of LDNH$_0$ is also very close to that of
the Euler equations. Figure~\ref{fig:tw:vertical} compares, for $A=0.3$, the horizontal and vertical
velocities of each layer and the non-hydrostatic pressure at the interfaces of LDNH$_0$ with the layer
means of the Euler solution. For $N=3$ the maximum difference is below $3\%$ of the maximum of $u$, $w$ and $p$, and
for $N=6$ it is below $1\%$. As for the
speed, it decreases roughly like $N^{-2}$, and the vertical velocity and the pressure, which are the
ones that depend on the vertical structure, are the most sensitive. Table~\ref{tab:tw:relax} also gives
the minimum phase speed of the fast branch. For $c=3\sqrt{gH_0}$ with $N\ge3$, and for
$c=5\sqrt{gH_0}$ with $N=6$, it is below the speed of the wave (Remark~\ref{rem:resonance}); there
the computed waves are those of the truncated Fourier problem, which does not resolve the short
resonant waves.

\begin{table}[htbp]
\centering
\small
\setlength{\tabcolsep}{5pt}
\begin{tabular}{cccccc}
\toprule
$N$ & & $c=3\sqrt{gH_0}$ & $c=5\sqrt{gH_0}$ & $c=10\sqrt{gH_0}$ & $c=30\sqrt{gH_0}$\\
\midrule
1 & $\min v_{\mathrm{fast}}/\sqrt{gH_0}$ & 2.35 & 3.67 & 7.14 & 21.24\\
 & $|C-C_{\mathrm{LDNH}_0}|/C_{\mathrm{LDNH}_0}$, $A=0.2$ & 8.9e-04 & 3.3e-04 & 8.3e-05 & 9.3e-06\\
\midrule
2 & $\min v_{\mathrm{fast}}/\sqrt{gH_0}$ & 1.20 & 1.95 & 3.85 & 11.49\\
 & $|C-C_{\mathrm{LDNH}_0}|/C_{\mathrm{LDNH}_0}$, $A=0.2$ & 6.6e-04 & 2.4e-04 & 6.1e-05 & 6.7e-06\\
\midrule
3 & $\min v_{\mathrm{fast}}/\sqrt{gH_0}$ & 0.81 & 1.32 & 2.60 & 7.77\\
 & $|C-C_{\mathrm{LDNH}_0}|/C_{\mathrm{LDNH}_0}$, $A=0.2$ & 6.2e-04 & 2.3e-04 & 5.7e-05 & 6.3e-06\\
\midrule
6 & $\min v_{\mathrm{fast}}/\sqrt{gH_0}$ & 0.41 & 0.67 & 1.32 & 3.93\\
 & $|C-C_{\mathrm{LDNH}_0}|/C_{\mathrm{LDNH}_0}$, $A=0.2$ & 5.9e-04 & 2.2e-04 & 5.4e-05 & 6.1e-06\\
\bottomrule
\end{tabular}

\caption{Relaxed system: minimum phase speed of the fast branch (linearized at rest) and relative
difference between the speed of the solitary wave and that of LDNH$_0$ ($A=0.2$).}
\label{tab:tw:relax}
\end{table}

\subsection{Order of convergence}\label{sec:tests:order}
The exact solitary waves are also exact solutions on a periodic domain: the solution at time $t$
is the initial profile translated by $Ct$. We take $A=0.2$, $H_0=1$, $N=3$ and the domain
$[-40,40]$, and we compute up to $T=4$ with the relaxed scheme ($c=10\sqrt{gH_0}$, starting from
the exact solitary wave of the relaxed system) and with the projection method (starting from
the exact solitary wave of LDNH$_0$). Both schemes converge with order $2$ in all the variables
(Table~\ref{tab:order:exact}). For large values of $c$ the asymptotic
order of the pressure is only reached on fine meshes, because the splitting error behaves like
$(c\,\Delta t)^2$ with $\Delta t\propto\Delta x$.

\begin{table}[htbp]
\centering
\small
\begin{tabular}{lrcccccc}
\toprule
scheme & cells & \multicolumn{2}{c}{$\eta$} & \multicolumn{2}{c}{$u_1$} & \multicolumn{2}{c}{$p_{1/2}$}\\
\midrule
relaxed, $c=10\sqrt{gH_0}$ & 100 & 1.34e-01 & -- & 3.23e-01 & -- & 5.37e+00 & --\\
 & 200 & 4.41e-02 & 1.61 & 1.08e-01 & 1.58 & 1.38e+00 & 1.96\\
 & 400 & 1.15e-02 & 1.94 & 2.84e-02 & 1.93 & 2.76e-01 & 2.32\\
 & 800 & 2.91e-03 & 1.98 & 7.18e-03 & 1.98 & 7.05e-02 & 1.97\\
 & 1600 & 7.37e-04 & 1.98 & 1.81e-03 & 1.99 & 1.70e-02 & 2.05\\
\midrule
projection (LDNH$_0$) & 100 & 1.41e-01 & -- & 3.17e-01 & -- & 5.34e-01 & --\\
 & 200 & 4.66e-02 & 1.60 & 1.03e-01 & 1.62 & 2.07e-01 & 1.37\\
 & 400 & 1.15e-02 & 2.02 & 2.63e-02 & 1.97 & 4.86e-02 & 2.09\\
 & 800 & 2.82e-03 & 2.03 & 6.53e-03 & 2.01 & 1.16e-02 & 2.06\\
 & 1600 & 6.98e-04 & 2.01 & 1.63e-03 & 2.01 & 2.83e-03 & 2.04\\
\bottomrule
\end{tabular}

\caption{Order of convergence against exact solitary waves ($A=0.2$, $N=3$, $T=4$): $L^1$ errors
and orders.}
\label{tab:order:exact}
\end{table}

\subsection{Solitary wave: influence of the relaxation celerity}\label{sec:tests:soliton}
The exact solitary wave of LDNH$_0$ gives a clean way to measure the error of the relaxation. We
take it with $A=0.2$ and $H_0=1$, for $N=1,3,6$ layers, as initial data of the relaxed scheme, with its
pressure. Since the exact solution of LDNH$_0$ is the same profile translated by $Ct$, the difference is the
error of the relaxation plus the discretization error, and the projection method, run with the same data,
gives the discretization error alone. We use the periodic domain $[-40,40]$ with $6400$ cells, $T=15$,
$c/\sqrt{gH_0}\in\{1,2,3,5,10,20,30\}$ and two CFL numbers, $0.45$ and $0.225$.

\begin{figure}[htbp]
\centering
\includegraphics[width=0.85\textwidth]{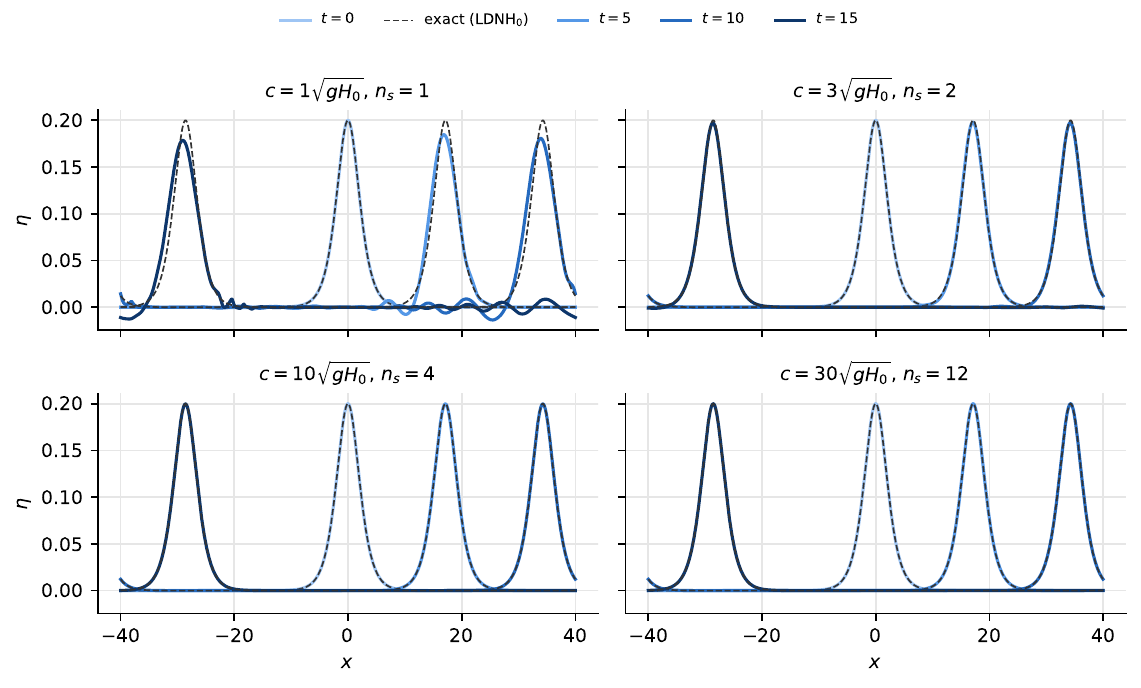}
\caption{Exact solitary wave of LDNH$_0$, $N=3$: free surface of the relaxed scheme at
$t=0,5,10,15$ for several relaxation celerities. Dashed lines: exact solution of LDNH$_0$.}
\label{fig:soliton:snapshots}
\end{figure}

Figure~\ref{fig:soliton:snapshots} shows the free surface for $N=3$. The scheme is stable for every
value of $c$, also for $c=\sqrt{gH_0}$, which is smaller than the speed of the wave: the relaxation
coupling conserves the energy and does not need a subcharacteristic condition. However, with
$c=\sqrt{gH_0}$ the fast branch is not separated from the gravity waves, the amplitude decreases from
$0.20$ to $0.18$ (Table~\ref{tab:soliton:c}) and the wave leaves a wake. From $c=3\sqrt{gH_0}$ on, the
free surface cannot be distinguished from the exact solution in the figure.

Table~\ref{tab:soliton:c} shows the errors for $N=3$. The error in the free surface decreases like
$c^{-2}$, which is the formal order of the relaxation, down to $3.1\cdot10^{-5}$ for $c=30\sqrt{gH_0}$,
close to the discretization error of the projection method ($2.1\cdot10^{-5}$). The error in the pressure
also decreases with $c$, but more slowly for $c\ge20\sqrt{gH_0}$. This remaining error comes from the
initial data: the solitary wave of LDNH$_0$ is not a travelling wave of the relaxed system, and the
difference excites fast oscillations of the pressure, of frequency of order $c/H_0$, which are not damped
since $\mathcal S_{\mathrm{NH}}$ conserves the energy. Starting from the exact solitary wave of the relaxed
system, the error in the pressure is the splitting error: during the hydrostatic step the flow evolves
without the pressure force, the constraint moves away from zero by $O(\Delta t)$, and the pressure reacts
with fast oscillations of amplitude proportional to $(c\Delta t)^2$. For $c=30\sqrt{gH_0}$ it decreases
from $1.4\cdot10^{-3}$ to $9.8\cdot10^{-5}$ when the CFL number is reduced from $0.45$ to $0.1$ (a value used only for this check); this is the
practical consequence of the fact that the scheme is not asymptotic-preserving. The results for $N=1$ and
$N=6$ are similar (not shown). These results support the choice $3\sqrt{gH_0}\le c\le10\sqrt{gH_0}$:
smaller values do not separate the gravity waves from the relaxation waves, and larger values increase the
cost and the splitting error in the pressure, unless the time step is reduced.
Note that, for $c=3\sqrt{gH_0}$ with $N\ge3$, and for $c=5\sqrt{gH_0}$ with $N=6$, the fast branch is in linear resonance with the solitary wave (Section~\ref{sec:tests:tw}), but the resulting tails are not visible in these tests.

\begin{table}[htbp]
\centering
\small
\begin{tabular}{rrccc}
\toprule
$c/\sqrt{gH_0}$ & $n_s$ & $\max_t\|\eta-\eta_{\mathrm{ex}}\|_\infty$ & $\max_t\|p_{1/2}-p_{1/2,\mathrm{ex}}\|_\infty$ & $\max_x\eta(T)$\\
\midrule
1 & 1 & 3.74e-02 & 1.61e-01 & 0.1785\\
2 & 1 & 1.12e-02 & 6.90e-02 & 0.1940\\
3 & 2 & 5.10e-03 & 3.99e-02 & 0.1972\\
5 & 2 & 1.86e-03 & 1.80e-02 & 0.1990\\
10 & 4 & 4.52e-04 & 7.16e-03 & 0.1998\\
20 & 8 & 9.70e-05 & 2.77e-03 & 0.1999\\
30 & 12 & 3.11e-05 & 2.40e-03 & 0.2000\\
\midrule
projection & -- & 2.10e-05 & 8.32e-05 & 0.2000\\
\bottomrule
\end{tabular}

\caption{Exact solitary wave of LDNH$_0$, $N=3$, $\mathrm{CFL}=0.45$: number of sub-steps, maximum
errors over the saved times against the exact solution of LDNH$_0$, and amplitude at $T=15$. Last
row: projection method (discretization error).}
\label{tab:soliton:c}
\end{table}

\subsection{Energy budget}\label{sec:tests:energy}
We measure the change of the total energy
\[
\begin{aligned}
E={}&\Delta x\sum_i\Bigl[\frac12\,g\,\eta_i^2+\sum_{\alpha=1}^N\frac12\,h_{\alpha,i}\,u_{\alpha,i}^2\Bigr]\\
&+\Delta x\sum_i\sum_{\alpha=1}^N\Bigl[\frac12\,h_{\alpha,i+1/2}\,w_{\alpha,i+1/2}^2
+\hs_{\alpha-1/2,i+1/2}\,\frac{p_{\alpha-1/2,i+1/2}^2}{2c^2}\Bigr],
\end{aligned}
\]
which is the discrete counterpart of the integral of $\mathcal E$ of Theorem~\ref{thm:energy}, in each
part of the step for the exact solitary wave of LDNH$_0$ of Section~\ref{sec:tests:soliton}, with $N=3$,
$c=5\sqrt{gH_0}$ and several meshes (Table~\ref{tab:energy}). Up to the factor $\Delta x$ and to the
potential energy, which does not change in $\mathcal S_{\mathrm{NH}}$, $E$ is the energy of
Proposition~\ref{prop:stagenergy}. The table gives three quantities, relative to the initial energy $E_0$.
The first one, $\Delta E_{\mathrm H}$, is the sum, over all the steps $\mathcal S_{\mathrm H}$ up to $T=15$, of the change of $E$ in each step.
The second one, $\Delta E_{\mathrm{NH}}$, is the same sum for the steps $\mathcal S_{\mathrm{NH}}$.
Finally, $\max|\Delta E_\tau|$ is the largest change of the modified energy $E_\tau$ of Theorem~\ref{thm:stab} during one step $\mathcal S_{\mathrm{NH}}$.

The step $\mathcal S_{\mathrm{NH}}$ conserves $E_\tau$ up to round-off errors. The change
$\Delta E_{\mathrm{NH}}$ is not a creation of energy. Since $\mathcal S_{\mathrm{NH}}$ conserves
$E_\tau=E-\frac{\tau^2}{8}|B_u^T\mathbf p|^2_{\mathbb M_c^{-1}}$, the change of $E$ in each step
$\mathcal S_{\mathrm{NH}}$ is the change of this $O(\tau^2)$ term (Remark~\ref{rem:taucross}), and we
have checked that $\Delta E_{\mathrm{NH}}$ is the sum of these changes up to round-off errors, on
the four meshes and for the two CFL numbers (the relative difference is below $2\cdot10^{-3}$ on the three coarser meshes and below $0.4$ on the finest one, where $\Delta E_{\mathrm{NH}}$ is itself close to the round-off level). It can have either sign (it is negative on the finest
mesh) and it decreases quickly as the mesh is refined. Its contribution to the change of $E$ is
between three and six orders of magnitude smaller than that of $\mathcal S_{\mathrm H}$. Therefore, the
energy loss of the scheme is the numerical diffusion of the hydrostatic step.
This is the expected behavior: over a flat bottom, the energy is an entropy of the hydrostatic system, and the HLL solver of the hydrostatic step dissipates it.
This loss decreases with order
three, as expected for a second-order upwind scheme, and it does not depend on the CFL number (the
results for $\mathrm{CFL}=0.225$ are almost the same).

\begin{table}[htbp]
\centering
\small
\begin{tabular}{r cccc}
\toprule
cells & $\Delta E_{\mathrm H}/E_0$ & order & $\Delta E_{\mathrm{NH}}/E_0$ & $\max|\Delta E_\tau|/E_0$ \\
\midrule
400 & -5.3e-03 & -- & 1.2e-06 & 5.9e-16 \\
800 & -5.7e-04 & 3.22 & 8.1e-08 & 7.0e-16 \\
1600 & -6.6e-05 & 3.12 & 3.9e-09 & 9.4e-16 \\
3200 & -7.8e-06 & 3.08 & -3.9e-11 & 1.4e-15 \\
\bottomrule
\end{tabular}

\caption{Exact solitary wave of LDNH$_0$, $N=3$, $c=5\sqrt{gH_0}$, $\mathrm{CFL}=0.45$: changes of the energy up to
$T=15$ in the hydrostatic steps ($\Delta E_{\mathrm H}$) and in the non-hydrostatic steps
($\Delta E_{\mathrm{NH}}$), relative to the initial energy $E_0$, order of $\Delta E_{\mathrm H}/E_0$,
and largest change of the modified energy $E_\tau$ during one non-hydrostatic step.}
\label{tab:energy}
\end{table}

\subsection{Well-prepared initial data}\label{sec:tests:prepared}
\begin{figure}[htbp]
\centering
\includegraphics[width=\textwidth]{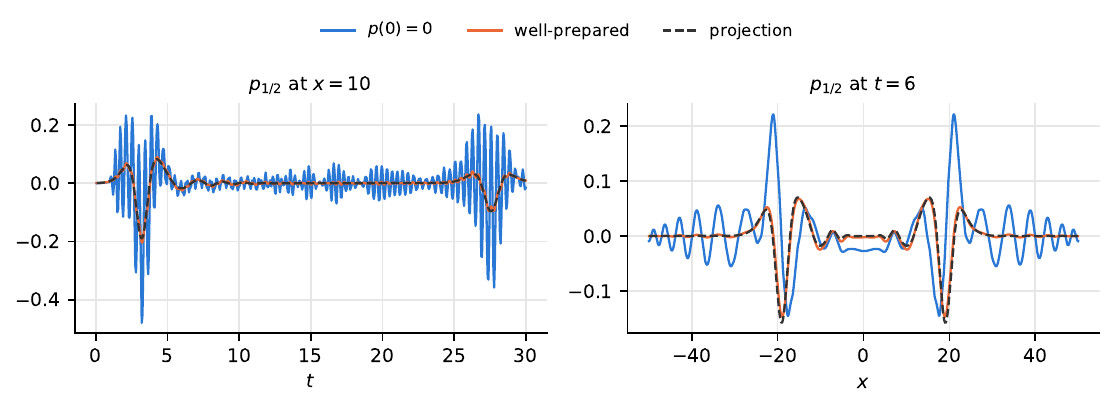}
\caption{Hump at rest, $N=3$, $c=3\sqrt{gH_0}$: bottom non-hydrostatic pressure with $p(0)=0$,
with well-prepared initial data and with the projection method.}
\label{fig:prepared}
\end{figure}

The aim of this test is to show the importance of using well-prepared initial data, as defined in Section~\ref{sec:scheme:init}.
As explained there, initial data satisfying $\mathcal K=0$ and $p=0$ do not satisfy $\dt\mathcal K=0$, and they excite the fast branch of the relaxed system with oscillations whose amplitude is similar to that of the pressure itself, for any value of $c$.
To show this, we consider a hump of the free surface at rest over a flat bottom,
$\eta=A\,\mathrm{sech}^2(x/\ell)$ with $A=0.3$, $H_0=1$ and $\ell=H_0\sqrt{(H_0+A)/A}$, in the
periodic domain $[-50,50]$ with $2000$ cells, $N=3$, $c=3\sqrt{gH_0}$ and $T=30$. With $p(0)=0$ the bottom
pressure shows fast oscillations with the frequency of the lowest relaxation mode
(Figure~\ref{fig:prepared}). Their amplitude is similar to the pressure itself, and they are not
damped. With well-prepared initial data they disappear: the maximum difference with the
projection method decreases from $0.60$ to $0.056$, while $\max|p|\approx0.61$. The free surface
is almost not affected (differences of $5\cdot10^{-3}$ in both cases).

\subsection{Run-up of a solitary wave on a plane beach}\label{sec:tests:wd}
\begin{figure}[htbp]
\centering
\includegraphics[width=\textwidth]{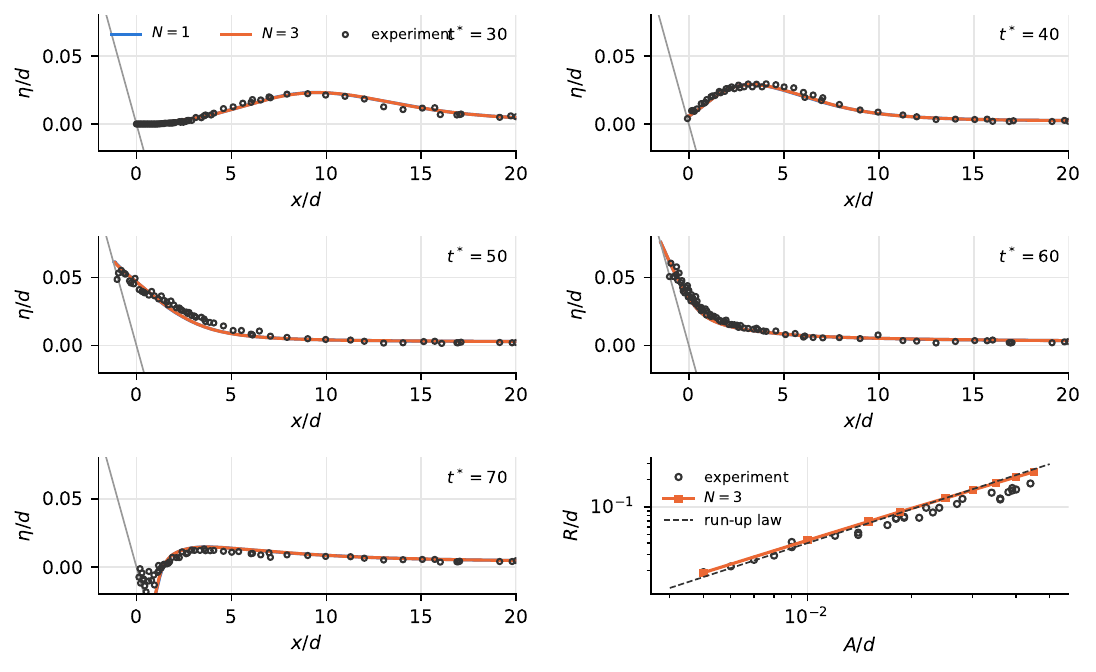}
\caption{Solitary wave on a plane beach, $A/d=0.0185$: free surface at
$t^*=t\sqrt{g/d}=30,\dots,70$ for $N=1$ and $N=3$, and laboratory data \cite{Synolakis1987}.
Bottom right: maximum run-up versus amplitude.}
\label{fig:wd:synolakis}
\end{figure}

In order to validate the wet-dry treatment of Section~\ref{sec:scheme:wd}, we consider the
solitary wave on the plane beach $1{:}19.85$ of \cite{Synolakis1987}, with the laboratory data of
the NTHMP benchmark BP4 \cite{SynolakisBenchmark2008}. As in \cite{Synolakis1987}, the
still-water depth is denoted by $d$, that is, $d=H_0$. We use $N=3$ (and $N=1$), $c=5\sqrt{gd}$
and $d=1$. The domain is $[-10,60]$ (in units of $d$), with $1400$ cells and open boundaries, and
the initial condition is the solitary wave
\[
\eta=A\,\mathrm{sech}^2\Bigl(\frac{x-x_1}\ell\Bigr),\qquad \ell=d\sqrt{\frac{d+A}{A}},
\]
with the velocity
\[
u=-C\,\frac{\eta}{d+\eta},\qquad C=\sqrt{g(d+A)},
\]
in all the layers, so that the wave travels towards the beach. The shoreline is at $x=0$ and the toe of the beach at $x=19.85$. The crest of the wave is initially at $x_1=19.85+\operatorname{arccosh}(\sqrt{20})/\gamma$, with $\gamma=\sqrt{3A/(4d^3)}$ as in \cite{Synolakis1987}, that is, $x_1\approx38.3$ for $A/d=0.0185$. For $A/d=0.0185$, the free surface agrees very well with the measurements at
all times, with one and with three layers (Figure~\ref{fig:wd:synolakis}). The maximum run-up is
$R/d=0.0886$, while the run-up law
\[
\frac Rd=2.831\sqrt{\cot\beta}\,\Bigl(\frac Ad\Bigr)^{5/4}
\]
gives $0.086$, and the measurements give $0.074$--$0.078$. The run-up of the non-breaking waves
follows the law and the measurements for all the computed amplitudes, $0.005\le A/d\le0.045$.
The lake at rest with emerged areas is tested in two dimensions in Section~\ref{sec:gpu:wd}.

\subsection{Periodic waves over a submerged bar}\label{sec:tests:bb}
We consider the non-breaking experiment of \cite{BejiBattjes1993,Dingemans1994}. Periodic waves
of period $2.02$\,s and height $0.02$\,m propagate in water of depth $0.4$\,m over a trapezoidal
bar: a $1{:}20$ upslope from $x=6$\,m to $x=12$\,m, a crest of depth $0.1$\,m up to $x=14$\,m,
and a $1{:}10$ downslope to $x=17$\,m. Over the bar the waves become steeper and generate higher
harmonics, which are released behind the bar into deeper water, where they are dispersive. This
is a classical test of the dispersive properties of the models. The waves are generated at $x=0$
with a ghost state in the hydrostatic step. We use $1800$ cells in $[0,45]$\,m,
$c=5\sqrt{gH_0}$ and $T=38$\,s, with free outflow on the right; the domain is long enough to
avoid reflections from the right boundary before $T$.

\begin{figure}[htbp]
\centering
\includegraphics[width=0.9\textwidth]{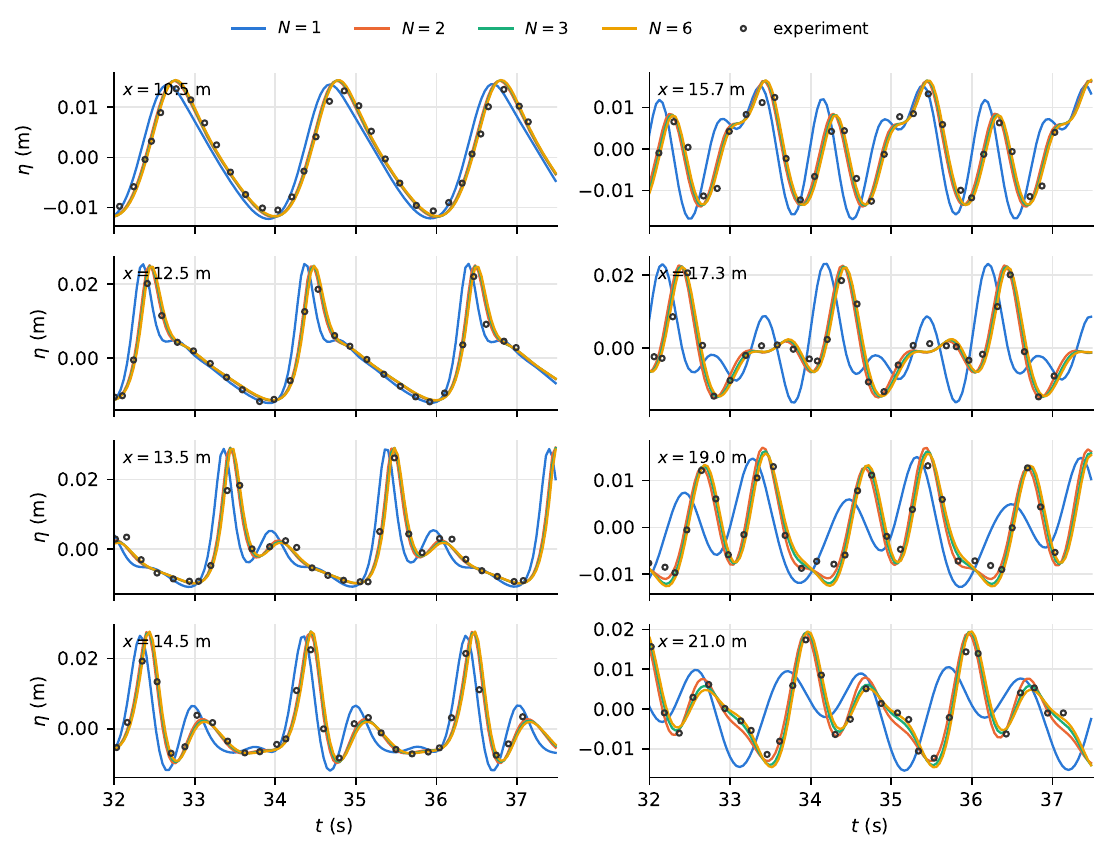}
\caption{Periodic waves over a submerged bar: free surface at the gauges for $N=1,2,3,6$ and
laboratory data.}
\label{fig:bb}
\end{figure}

Figure~\ref{fig:bb} and Table~\ref{tab:bb} compare the free surface with the laboratory data at
eight gauges. At each gauge, the difference with the data is measured by the root mean square
(RMS) difference
\[
\mathrm{RMS}=\Bigl(\frac1J\sum_{j=1}^J\bigl(\eta(t_j)-\eta_{\mathrm{lab}}(t_j)\bigr)^2\Bigr)^{1/2},
\]
where $t_1,\dots,t_J$ are the times of the laboratory record (shifted to the time origin of the
simulation) in the interval $[30.5,37.5]$\,s, about the last three and a half periods, the
computed free surface being interpolated at these times. Before the bar, all the models agree with each other and with the data. Behind the
bar, the single-layer model ($N=1$) is out of phase, with errors of about $1$\,cm. With $N\ge2$
layers the released harmonics are well reproduced: the errors are of about $3$\,mm, of the order of the accuracy of the digitized data, and the results for
$N=2$, $3$ and $6$ are almost the same.

\begin{table}[htbp]
\centering
\small
\begin{tabular}{ccccccccc}
\toprule
$N$ & 10.5 & 12.5 & 13.5 & 14.5 & 15.7 & 17.3 & 19 & 21\\
\midrule
1 & 3.09 & 4.89 & 6.84 & 4.88 & 9.38 & 12.70 & 8.54 & 10.08\\
2 & 1.34 & 1.41 & 2.30 & 2.34 & 3.13 & 3.03 & 2.23 & 2.08\\
3 & 1.16 & 1.55 & 1.89 & 2.99 & 2.29 & 2.29 & 1.73 & 1.56\\
6 & 1.10 & 1.82 & 1.80 & 3.45 & 1.98 & 2.20 & 1.98 & 1.98\\
\bottomrule
\end{tabular}

\caption{Periodic waves over a submerged bar: RMS difference (mm) with the laboratory data at the
gauges (position in m).}
\label{tab:bb}
\end{table}

\section{Extension to two dimensions on GPU}\label{sec:gpu}
We have extended both the relaxed scheme and the projection method to two dimensions on GPU. Since the projection method is expensive, we have implemented it only for comparison, and the new approach presented here has been implemented by taking as starting point the GPU code of the single-layer case \cite{Escalante2018AMC,Escalante2019JCP}. The
two codes are separated, but they share the hydrostatic part and the spatial discretization of the
non-hydrostatic terms, so that the comparison between them is fair. We begin this section with a general description of
the two codes.
Then, we check the order of convergence by means of a solitary wave that crosses the mesh obliquely, and we compare both codes for a Gaussian hump.
Next, the relaxed scheme is validated in the presence of wet-dry fronts and with two standard laboratory experiments: the solitary wave on a conical island and the periodic waves over a semicircular shoal.
Finally, the cost of both approaches is compared. As in Section~\ref{sec:tests}, the tests without laboratory data have $H=1$\,m and $g=9.81$\,m\,s$^{-2}$.

\subsection{Description of the 2D scheme}\label{sec:gpu:disc}
The depth $h$ and the momenta of the layers $(h_\alpha u_\alpha,\,h_\alpha v_\alpha)$, where
$v_\alpha$ is the velocity of the layer $\alpha$ in the $y$ direction, are cell averages on a
Cartesian mesh. We denote by $\mathbf u_\alpha=(u_\alpha,v_\alpha)$ the velocity of the layer $\alpha$. The vertical velocities $w_\alpha$ and the interface pressures $p_\aph$ are
placed at the vertices, the two-dimensional counterpart of the cell edges of
Section~\ref{sec:scheme:stag}. The hydrostatic part, which includes the exchange terms between the layers, is solved using the well-balanced path-conservative HLL
scheme of the hydrostatic multilayer system
\cite{ABPSM2011,Castro2011CRM}. At a
vertex $v$, let us denote by $\mathcal N_v$ the set of the four cells that share it.
The constraint \eqref{eq:K} is then discretized as
\[
\mathcal K_{\aph,v}=\sum_{j\in\mathcal N_v}\Bigl(\mathbf a_{\aph,vj}\cdot\mathbf u_{\alpha+1,j}
+\mathbf b_{\aph,vj}\cdot\mathbf u_{\alpha,j}\Bigr)+(w_{\alpha+1}-w_\alpha)_v,
\]
\[
\begin{aligned}
\mathbf a_{\aph,vj}&=\Bigl(\frac{s^x_{vj}h_{\alpha,v}}{4\Delta x}-\frac{\partial_xz_{\aph,v}}4,\
\frac{s^y_{vj}h_{\alpha,v}}{4\Delta y}-\frac{\partial_yz_{\aph,v}}4\Bigr),\\
\mathbf b_{\aph,vj}&=\Bigl(\frac{s^x_{vj}h_{\alpha,v}}{4\Delta x}+\frac{\partial_xz_{\aph,v}}4,\
\frac{s^y_{vj}h_{\alpha,v}}{4\Delta y}+\frac{\partial_yz_{\aph,v}}4\Bigr),
\end{aligned}
\]
where $s^x_{vj}=\pm1$ for the cells east and west of $v$ (and $s^y_{vj}$ likewise),
$h_{\alpha,v}$ is the mean of the four cell depths divided by $N$, and the gradients of the
interfaces are the means of the two jumps \eqref{eq:dzHR} across the vertex in each direction, divided by $\Delta x$ or $\Delta y$.
The pressure force on the cells is the transpose of this constraint, as in
\eqref{eq:forcestag}, so that the semi-discrete energy balance of
Proposition~\ref{prop:stagenergy} holds. In a setting uniform in $y$, this discretization
reduces to the one-dimensional scheme of Section~\ref{sec:scheme:stag}.

\subsubsection{Projection method}
We use the second-order LSDIRK2 scheme in time (\ref{app:proj}), so that both codes
are second-order accurate.
Each of its two implicit stages is a projection of Section~\ref{sec:scheme:stag}: it solves the
system $B\,\mathbb M^{-1}B^T\boldsymbol\psi=-\mathcal K$, where the
impulse $\boldsymbol\psi$ gives the pressure of the stage. The pressure force and the constraint use the
depth of the projected state, so the matrix is symmetric and positive definite, and the
projection is orthogonal in the energy norm. In two dimensions this system is a nine-point stencil
with $N\times N$ tridiagonal blocks. The coefficients are assembled at each stage, and the system
is solved by the conjugate gradient method preconditioned with the block Jacobi method (one
tridiagonal solve per vertex for the vertical coupling). The iteration starts from the pressure
of the previous step, and it stops when the residual scaled by the block diagonal is below the
tolerance in the maximum norm.

\subsubsection{Relaxed scheme}\label{sec:gpu:relax}
The time step is the splitting \eqref{eq:strang} with the U--P--U step
\eqref{eq:UPU}: the updates of the velocities are explicit in the cells, and the Crank--Nicolson
step for $(w,p)$ is a tridiagonal system at each vertex, with the horizontal part of the
constraint as a known source. Theorem~\ref{thm:stab} holds without changes. The divergence at
the vertices has the same maximal symbol as the one-dimensional divergence at the edges, so for
uniform states the step is stable for $c\tau\le\min(\Delta x,\Delta y)$. The well-prepared
initial data are computed with the solver of the projection method, which is only used at $t=0$.
At the shoreline, both codes use the limited coefficients of Section~\ref{sec:scheme:wd} at the
vertices, together with a single layer (the mean of the velocities) in the cells with
$h<2\Delta x$ (\ref{app:wetdry}); the pressure is never switched off. We have checked numerically, on periodic meshes of $8\times8$ square cells with $N=1$, $3$
and $6$, random depths between $10^{-3}$ and $4$ and random jumps of the free surface, that
$\lambda_{\max}\Delta x^2/4\le0.99$, while the stability condition with $\mathrm{CFL}_{\mathrm{NH}}=0.9$
only needs $\lambda_{\max}\Delta x^2/4<1/0.81$.

In the GPU code, the last half step of $\mathcal S_{\mathrm H}$ of each time step is merged with the first
half step of the next one, as usual in Strang splitting:
\[
\mathcal S_{\mathrm H}\!\left(\tfrac{\Delta t_{n+1}}{2}\right)\circ\mathcal S_{\mathrm H}\!\left(\tfrac{\Delta t_n}{2}\right)
\quad\text{is replaced by}\quad
\mathcal S_{\mathrm H}\!\left(\tfrac{\Delta t_n+\Delta t_{n+1}}{2}\right).
\]
If $\mathcal S_{\mathrm H}$ were exact, the two forms would be identical. With the second-order SSP
Runge--Kutta method they differ by $O(\Delta t^3)$ per step, so the scheme is still second order
(Table~\ref{tab:gpu:order}), and the hydrostatic step $\mathcal S_{\mathrm H}$ is applied once per time step instead of
twice. The half steps are separated only at the output times. The same is done in the non-hydrostatic part. Let $\mathcal U(s)$ denote the update $\mathbf u\leftarrow\mathbf u+s\,\mathbb M_c^{-1}B_u^T\mathbf p$ of the velocities with a given pressure. The second half update of a sub-step and the first half update of the next sub-step use the same pressure, so that
\[
\mathcal U\!\left(\tfrac\tau2\right)\circ\mathcal U\!\left(\tfrac\tau2\right)
\quad\text{is replaced by}\quad
\mathcal U(\tau),
\]
which gives the same result up to round-off. Thus, the $n_s$ sub-steps need $n_s+1$ updates of the velocities instead of $2n_s$.

With these two mergings, the scheme is equivalent to the Lie--Trotter splitting conjugated by $\mathcal S_{\mathrm H}(\Delta t/2)$. The half steps at the beginning, at the end and at the output times put the hydrostatic and the non-hydrostatic parts at the same time level, and they are needed to keep the second order at these times.

All the runs below use a single NVIDIA A100 GPU (40\,GB) in double precision, CFL $=0.45$, the
relaxed scheme with $\mathrm{CFL}_{\mathrm{NH}}=0.9$ and well-prepared initial data, and, unless stated
otherwise, $N=3$ layers. The tolerance of the iterative solver of the projection method is
$10^{-8}$ ($10^{-10}$ for the order tests). Both codes use the same time step, prescribed by the
hydrostatic part, $\Delta t=\mathrm{CFL}\,\Delta x/\bigl(2\max_{j,\alpha}(\max(|u_{\alpha,j}|,|v_{\alpha,j}|)+\sqrt{gh_j})\bigr)$ on square cells. The
wall times include the output of the results and do not include the initialization (initial
projection and well-prepared initial data). The codes use CUDA 12.4.

\subsection{Accuracy test}
To test the two-dimensional discretization, we rotate the exact solitary waves of
Section~\ref{sec:tests:tw} by $45^\circ$: the wave propagates along the diagonal of the square
$[-40\sqrt2,40\sqrt2]^2$ with periodic boundaries. The waves have period $80$ along their direction
of propagation, so the solution is exactly periodic in this square. We use $n_x\times n_x$ cells, $A=0.2$,
$N=3$, the exact wave of the relaxed system ($c=10\sqrt{gH}$) for the relaxed scheme and the exact
wave of LDNH$_0$ for the projection method, and we compute up to $T=4$.

\begin{figure}[htbp]
\centering
\includegraphics[width=\textwidth]{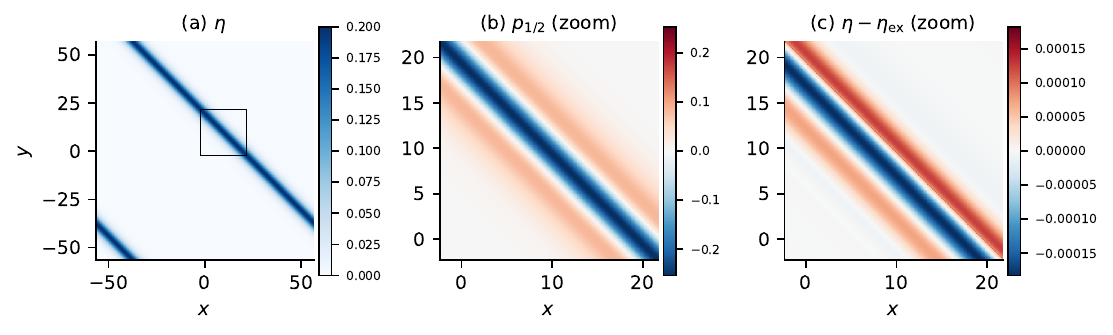}
\caption{Solitary wave propagating at $45^\circ$ in a periodic square, relaxed scheme with $N=3$,
$c=10\sqrt{gH}$ and $1600^2$ cells, at $T=4$: (a) free surface (the square marks the window of the
other panels); (b) non-hydrostatic pressure at the bottom; (c) error of the free surface.}
\label{fig:gpu:oblique}
\end{figure}

Both codes are second-order accurate in all the variables (Table~\ref{tab:gpu:order}). The errors
in the free surface and in the velocity are close to those of the one-dimensional codes with the
same number of cells per direction (Table~\ref{tab:order:exact}), and the errors in the pressure of
the relaxed scheme are a little smaller, because the two-dimensional time step is $1/\sqrt2$ times the one-dimensional one and the
splitting error behaves like $(c\Delta t)^2$. Figure~\ref{fig:gpu:oblique} shows the solution of
the relaxed scheme on the finest mesh. The error is small and mainly due to a small phase lag of the wave, and the solution does not show the orientation of the mesh.

\begin{table}[htbp]
\centering
\small
\begin{tabular}{lrcccccc}
\toprule
scheme & cells & \multicolumn{2}{c}{$\eta$} & \multicolumn{2}{c}{$u_1$} & \multicolumn{2}{c}{$p_{1/2}$}\\
\midrule
relaxed, $c=10\sqrt{gH}$ & $100^2$ & 1.61e-01 & -- & 3.48e-01 & -- & 2.78e+00 & --\\
 & $200^2$ & 4.96e-02 & 1.70 & 1.08e-01 & 1.69 & 9.48e-01 & 1.55\\
 & $400^2$ & 1.35e-02 & 1.88 & 3.24e-02 & 1.74 & 2.37e-01 & 2.00\\
 & $800^2$ & 3.34e-03 & 2.02 & 8.14e-03 & 1.99 & 6.06e-02 & 1.96\\
 & $1600^2$ & 8.27e-04 & 2.01 & 2.03e-03 & 2.00 & 1.41e-02 & 2.10\\
\midrule
projection (LDNH$_0$) & $100^2$ & 1.62e-01 & -- & 3.52e-01 & -- & 5.52e-01 & --\\
 & $200^2$ & 4.88e-02 & 1.73 & 1.01e-01 & 1.81 & 2.23e-01 & 1.31\\
 & $400^2$ & 1.31e-02 & 1.90 & 2.93e-02 & 1.78 & 5.74e-02 & 1.96\\
 & $800^2$ & 3.23e-03 & 2.02 & 7.35e-03 & 2.00 & 1.35e-02 & 2.08\\
 & $1600^2$ & 8.01e-04 & 2.01 & 1.84e-03 & 2.00 & 3.30e-03 & 2.03\\
\bottomrule
\end{tabular}

\caption{Two-dimensional GPU codes, $N=3$, solitary wave propagating at $45^\circ$: $L^1$ errors
per unit length of crest and orders against the exact solitary waves of the relaxed system
($c=10\sqrt{gH}$) and of LDNH$_0$ ($A=0.2$, $T=4$).}
\label{tab:gpu:order}
\end{table}

\subsection{Gaussian hump}
\begin{figure}[htbp]
\centering
\includegraphics[width=\textwidth]{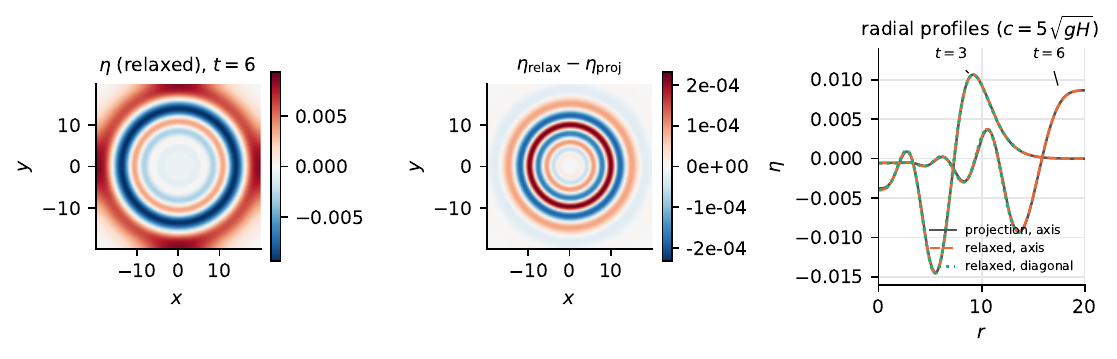}
\caption{Gaussian hump at $t=6$, $512^2$ cells, $N=3$: free surface of the relaxed scheme
($c=5\sqrt{gH}$), difference with the projection method, and profiles along the $x$ axis and along
the diagonal at $t=3$ and $t=6$.}
\label{fig:gpu:hump}
\end{figure}

We consider now a Gaussian hump defined by $\eta=0.1\,e^{-r^2/4}$ at rest over a flat bottom of unit depth, in the
periodic square $[-20,20]^2$. The final time is $T=6$. Both codes keep the symmetries of the mesh (axes and
diagonals) to round-off, and the solution is close to radially symmetric: the profiles along the
axes and along the diagonals almost coincide (Figure~\ref{fig:gpu:hump}). The differences between
the relaxed and the projection solutions are located on the dispersive wave train. They hardly
depend on the mesh and, as in Section~\ref{sec:tests:soliton}, they decrease like $c^{-2}$, from
less than $1\%$ of the amplitude for $c=3\sqrt{gH}$ to less than $0.1\%$ for $c=10\sqrt{gH}$
(Table~\ref{tab:gpu:hump}). The differences in the bottom pressure are larger, a few percent of its
maximum for $c\ge5\sqrt{gH}$; the largest ones, on the coarsest mesh, come
from the splitting error, which grows like $(c\Delta t)^2$.

\begin{table}[htbp]
\centering
\small
\begin{tabular}{rlrrcc}
\toprule
mesh & scheme & time (s) & speed-up & $\max|\Delta\eta|/A$ & $\max|\Delta p_{1/2}|/\max|p_{1/2}|$\\
\midrule
$256^2$ & projection (31 it.) & 3.94 & 1 & -- & --\\
 & relaxed, $c=3\sqrt{gH}$ & 0.232 & 17 & 6.4e-03 & 1.8e-01\\
 & relaxed, $c=5\sqrt{gH}$ & 0.245 & 16 & 2.2e-03 & 5.7e-02\\
 & relaxed, $c=10\sqrt{gH}$ & 0.259 & 15 & 4.5e-04 & 6.2e-02\\
\midrule
$512^2$ & projection (49 it.) & 28.3 & 1 & -- & --\\
 & relaxed, $c=3\sqrt{gH}$ & 1.37 & 21 & 6.6e-03 & 1.9e-01\\
 & relaxed, $c=5\sqrt{gH}$ & 1.54 & 18 & 2.3e-03 & 5.7e-02\\
 & relaxed, $c=10\sqrt{gH}$ & 1.62 & 18 & 5.5e-04 & 1.6e-02\\
\midrule
$1024^2$ & projection (72 it.) & 277 & 1 & -- & --\\
 & relaxed, $c=3\sqrt{gH}$ & 9.97 & 28 & 6.7e-03 & 1.8e-01\\
 & relaxed, $c=5\sqrt{gH}$ & 10.9 & 26 & 2.3e-03 & 5.7e-02\\
 & relaxed, $c=10\sqrt{gH}$ & 12.2 & 23 & 5.8e-04 & 1.5e-02\\
\bottomrule
\end{tabular}

\caption{Gaussian hump in two dimensions, $N=3$: wall times on an A100 GPU (the number of
iterations of the solver per stage in brackets), speed-up of the relaxed scheme, and maximum
differences $\Delta$ with the projection method at $T=6$, relative to the initial amplitude
$A=0.1$ and to the maximum bottom pressure.}
\label{tab:gpu:hump}
\end{table}

\subsection{Wet-dry fronts and laboratory experiments}\label{sec:gpu:wd}
Before the laboratory experiments, the treatment of the wet-dry fronts in two dimensions has been checked with two simple tests.
Firstly, both codes preserve a lake at rest with an emerged island, given by the bottom
\[
z_b=-1+1.5\,e^{-r^2/25}
\]
and the free surface $\eta=0$, on the periodic domain $[-10,10]^2$ with $200^2$ cells: the velocities and the free surface of the wet cells are zero up to round-off errors (below $10^{-14}$) after $T=10$.
Secondly, for the solitary wave on the plane beach of Section~\ref{sec:tests:wd}, the two-dimensional codes give the same run-up as in one dimension, $R/d=0.0886$, with the wave propagating along $x$ or along $y$ (channel $[-10,60]\times[-2\Delta x,2\Delta x]$, with $1400\times4$ cells and open boundaries).
Then, we consider two laboratory experiments: the solitary wave on a conical island, which involves the run-up and the diffraction of dispersive waves, and the periodic waves over a semicircular shoal, in which the dispersive effects are dominant.

\subsubsection{Conical island}\label{sec:gpu:island}
\begin{figure}[htbp]
\centering
\includegraphics[width=\textwidth]{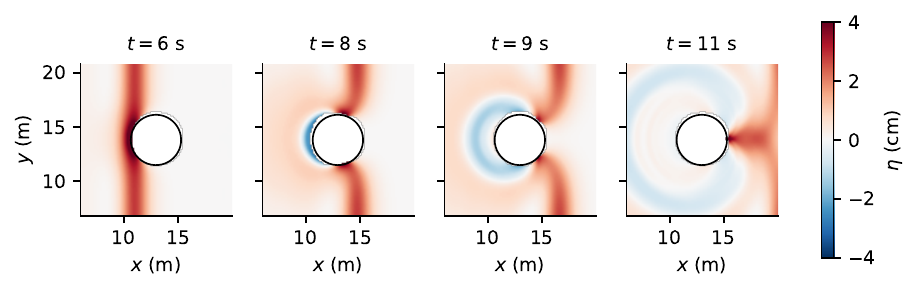}
\caption{Solitary wave on a conical island (BP6, case B), relaxed scheme with $N=3$: free surface
around the island (dry areas in white, still shoreline in black).}
\label{fig:gpu:islandplan}
\end{figure}

\begin{figure}[htbp]
\centering
\includegraphics[width=0.85\textwidth]{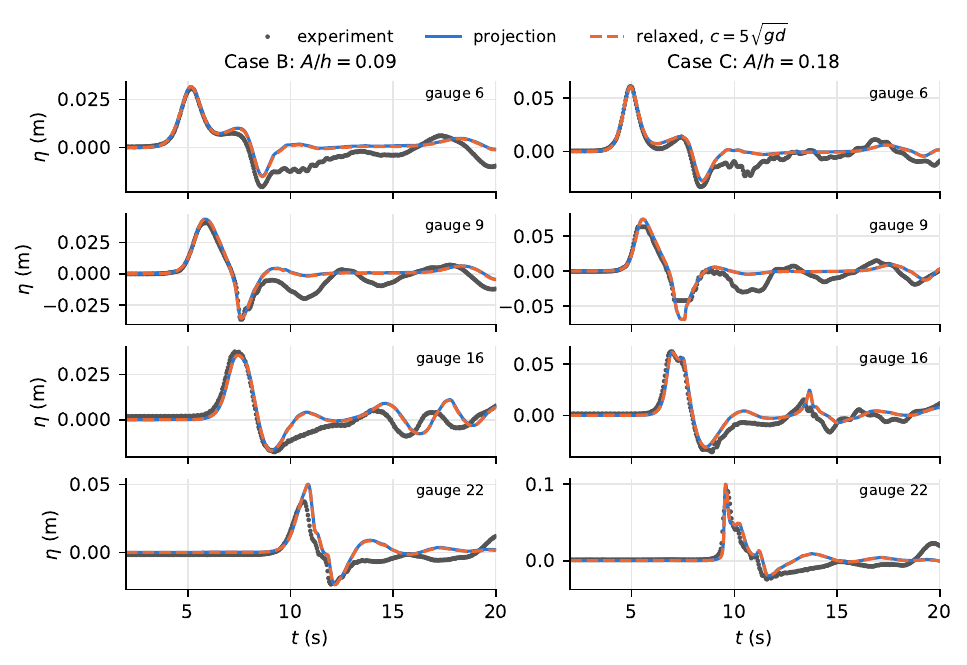}
\caption{Solitary wave on a conical island (BP6), set-up of~\cite{Kazolea2012}, $N=3$: free surface
at the gauges 6, 9, 16 and 22 for cases B and C. Original laboratory records, projection
method and relaxed scheme. The laboratory time is shifted so that the maximum at gauge 2
matches.}
\label{fig:gpu:island}
\end{figure}

\begin{figure}[htbp]
\centering
\includegraphics[width=0.65\textwidth]{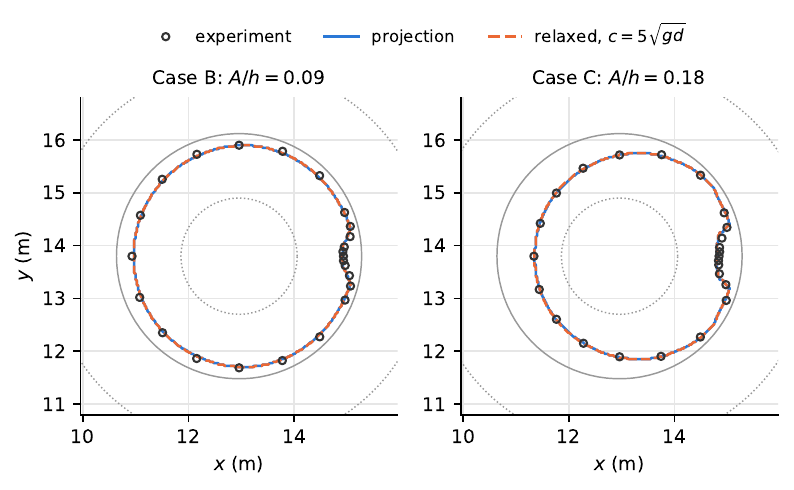}
\caption{Solitary wave on a conical island (BP6), set-up of~\cite{Kazolea2012}, $N=3$: limit of the
area wetted during the simulation (maximum depth above $2$\,mm), computed with the projection method and the
relaxed scheme, and measured run-up (original laboratory records) placed on the slope of the island. Gray lines: toe, still shoreline and
top of the island.}
\label{fig:gpu:runup}
\end{figure}

We reproduce the experiments of Briggs et
al.~\cite{Briggs1995} (NTHMP benchmark BP6), cases B and C, with the set-up of
Kazolea et al.~\cite{Kazolea2012}. As in the benchmark, the still-water depth is denoted by $d$, that is, $d=H=0.32$\,m, and the incident wave heights are
$A/d=0.09$ and $0.18$, slightly below the target values, because they fit better the record at
gauge 2. We use the original laboratory records of the NTHMP benchmark~\cite{SynolakisBenchmark2008}.
The domain is $[-5,28]\times[0,30]$\,m, with sponge layers of $3$\,m at both ends in $x$
and walls at $y=0$ and $y=30$\,m, and there is no bottom friction. We use $1320\times1200$ cells
($\Delta x= \Delta y = 2.5$\,cm) and $c=5\sqrt{gd}$. The initial condition is the solitary wave of
Section~\ref{sec:tests:wd}, centered at $x=0$ and travelling towards the island, that is, with $u=C\eta/(d+\eta)$. The incident wave travels around the
island and the two diffracted waves meet on the lee side (Figure~\ref{fig:gpu:islandplan}). Both
codes use the treatment of the shoreline of Section~\ref{sec:gpu:relax}: limited coefficients and a
single layer where $h<2\Delta x$, without switching off the pressure. The free
surface at the gauges (Figure~\ref{fig:gpu:island}) agrees with the measurements for the incident
wave and for the wave diffracted behind the island.
The largest difference is at gauge 22 in case B, on the lee side, where the computed crest is
higher than the measured one. The run-up around the island (Figure~\ref{fig:gpu:runup}; a point is
counted as wet when the maximum depth exceeds $2$\,mm) follows the measured distribution. With the
error measures of the NTHMP benchmark, the normalized root mean square deviation of the run-up is
$14.7\%$ (case B) and $8.5\%$ (case C), below the mean values of the eight models of the NTHMP benchmark,
$18\%$ and $11\%$~\cite{Macias2017}. The two codes give almost the same results: the differences at the
gauges are of a few millimeters, and the run-up errors are the same. For $T=20$\,s, the relaxed
code takes $51$ and $54$\,s, and the projection method $1248$ and $1601$\,s ($63$ and $75$
iterations of the solver per stage).

\subsubsection{Periodic waves over a semicircular shoal}\label{sec:gpu:whalin}
The experiments of Whalin~\cite{Whalin1971} are a standard test for the nonlinear refraction and
diffraction of dispersive waves. A semicircular shoal in a basin of width $6.096$\,m focuses
periodic waves, and the focusing transfers energy from the first harmonic to higher harmonics,
which are free waves of short length behind the shoal. The depth is $0.4572$\,m, and the shoal
rises with a slope $1{:}25$ to a depth of $0.1524$\,m. We use the three standard cases: periods
$T_{\mathrm p}=1,2,3$\,s with incident amplitudes $a=1.95$, $0.75$ and $0.68$\,cm. The laboratory data are
the amplitudes of the first three harmonics along the centerline; we digitized them from the
figures of the original report.

The basin is $[-L_{\mathrm p},30]\times[0,6.096]$\,m, where $L_{\mathrm p}$ is the wavelength, with
$\Delta x=2.5$\,cm ($3.6$--$4.0\cdot10^5$ cells with the sponge layer, depending on $T_{\mathrm p}$) and reflective side walls. The waves are generated in
a relaxation zone of one wavelength, $x<0$, where the state is relaxed at each time step towards
a second-order Stokes wave \cite{Mayer1998}; this zone also absorbs the waves reflected towards the
generator. A sponge layer of $5$\,m beyond $x=30$\,m, $30<x<35$\,m, absorbs the transmitted waves: there the
state is relaxed towards the state at rest at the rate $\kappa=5\,\vartheta^2$\,s$^{-1}$, where $\vartheta$ goes
from $0$ to $1$ across the layer. This rate does not depend on the time step, and the reflected
waves are small. We simulate
$50$\,s ($T_{\mathrm p}=1$\,s) and $45$\,s ($T_{\mathrm p}=2,3$\,s), and we compute the harmonics at the gauges over the
last four periods. The relaxed scheme uses $c=5\sqrt{gH}$ and $N=1,2,3$ (and $N=6,10,20$ for $T_{\mathrm p}=2$\,s); the projection
method uses $N=3$.

\begin{figure}[htbp]
\centering
\includegraphics[width=0.8\textwidth]{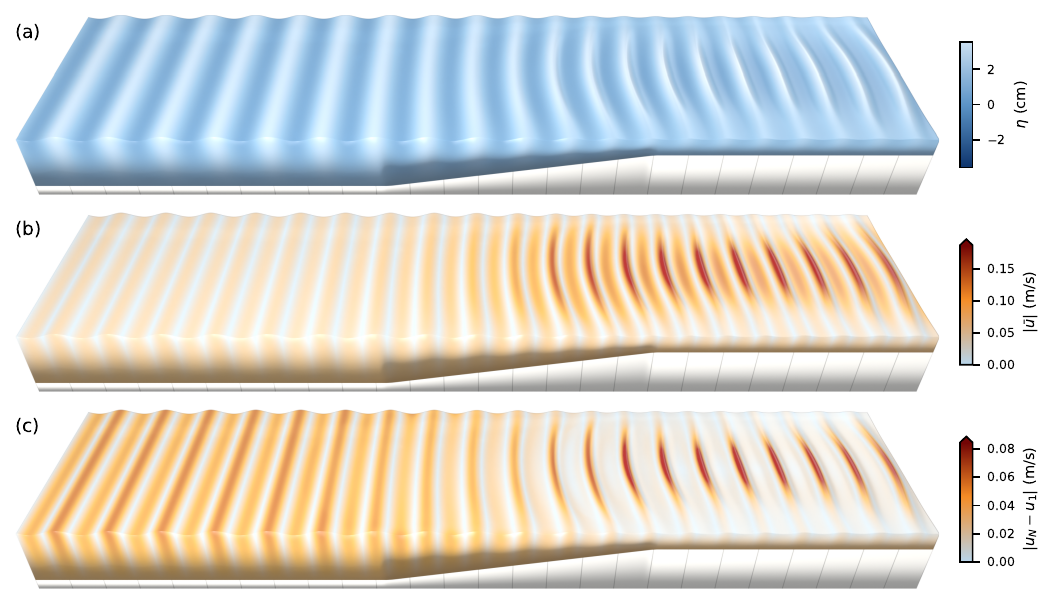}
\caption{Periodic waves over a semicircular shoal, $T_{\mathrm p}=1$\,s, relaxed scheme with $N=3$, at
$t=50$\,s, for $0\le x\le26$\,m: (a) free surface $\eta$; (b) speed of the depth-averaged
velocity; (c) vertical shear between the top and bottom layers, $|\mathbf u_N-\mathbf u_1|$. The
waves come from the left. Vertical scale exaggerated four times.}
\label{fig:gpu:whalin3d}
\end{figure}

Figure~\ref{fig:gpu:whalin3d} shows the solution for $T_{\mathrm p}=1$\,s: the shoal bends the wave fronts and
focuses them behind it, where the velocity is largest. The vertical shear between the top and the
bottom layers shows where the vertical structure of the flow matters: it is large for all the
waves in the deep part ($kH=1.9$), where the velocity decreases clearly with depth, and it almost
vanishes over the shoal, except at the focused crests. A single-layer model cannot represent this
shear.

\begin{figure}[htbp]
\centering
\includegraphics[width=\textwidth]{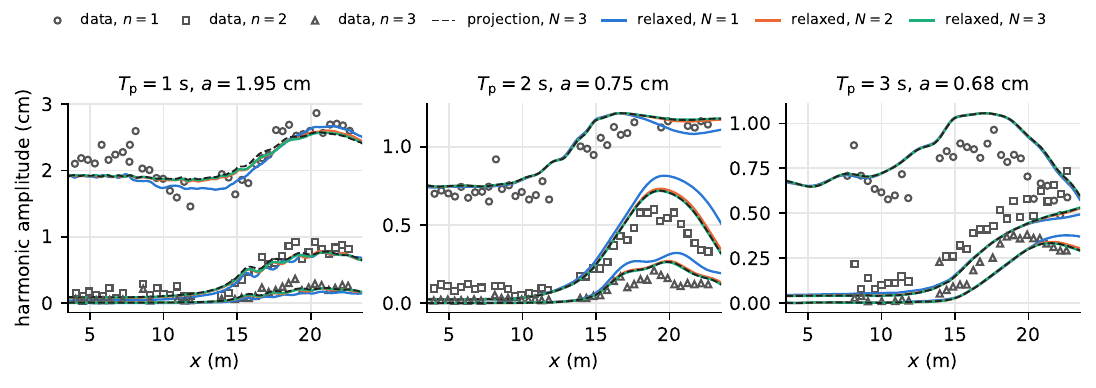}
\caption{Periodic waves over a semicircular shoal (Whalin): amplitudes of the first three
harmonics along the centerline. Symbols: laboratory data. Lines: relaxed scheme with $N=1,2,3$ and
projection method with $N=3$ (dashed).}
\label{fig:gpu:whalin}
\end{figure}

Figure~\ref{fig:gpu:whalin} compares the amplitudes of the harmonics along the centerline with the measurements.
Let us first discuss the effect of the number of layers.
For $T_{\mathrm p}=1$ and $2$\,s the waves are dispersive
($kH=1.9$ and $0.7$ in the deep part, and the higher harmonics are much shorter). With one layer,
the second and third harmonics behind the shoal are clearly overestimated for $T_{\mathrm p}=2$\,s (maxima of
$0.82$ and $0.32$\,cm, against about $0.6$ and $0.2$\,cm in the measurements). More layers reduce
them, and the results converge: the maximum of the second harmonic is $0.73$, $0.72$ and $0.71$\,cm
with $N=2$, $3$ and $6$, and it does not change with $N=10$ and $N=20$ ($0.71$\,cm). The
RMS difference with the measurements decreases from $0.16$\,cm ($N=1$) to $0.11$,
$0.099$ and $0.094$\,cm for the second harmonic, and from $0.088$ to $0.050$\,cm for the third
one. For $T_{\mathrm p}=3$\,s the waves are long and the number of layers has almost
no influence.

For $T_{\mathrm p}=2$\,s, the converged second harmonic is still about $18\%$
above the measurements. This is not caused by the dispersion of the model. Over the shoal the
second harmonic is generated by the first one, and the energy that it receives depends on the
phase mismatch $\Delta k=k(2\omega)-2k(\omega)$ between them: the smaller $\Delta k$, the larger the
second harmonic. With the linear dispersion relation of LDNH$_0$, at the depth of the top of the shoal, $\Delta k$ is $24\%$ smaller than
the exact value with $N=1$, $3\%$ smaller with $N=3$ and $0.7\%$ smaller with $N=6$. This explains
why the second harmonic decreases with $N$, and it shows that the converged result is close to the
result of a model with exact linear dispersion.

For $T_{\mathrm p}=3$\,s, the first harmonic is overestimated at the focus (about $1.05$\,cm
against $0.96$\,cm), and the second harmonic is too small before the shoal, where the measurements
already show a second harmonic of $0.1$--$0.2$\,cm. These differences do not change with a finer
mesh, a longer sponge layer, a longer simulation or the number of layers. We think that they come from the wave generation of the experiment,
which can produce free second harmonics that a second-order Stokes wave does not contain.

Finally, the curves of the relaxed scheme and of the
projection method with $N=3$ overlap for $T_{\mathrm p}=2$ and $3$\,s (the differences are negligible).
For $T_{\mathrm p}=1$\,s, the shortest waves, the differences are larger, but still much smaller than the differences with the measurements.
With $N=3$, the relaxed scheme takes $39$, $35$ and $37$\,s, and the
projection method $1928$, $1470$ and $1440$\,s ($142$, $119$ and $110$ iterations of the solver per
stage), so that the relaxed scheme is between $39$ and $50$ times faster.

\subsection{Performance study}\label{sec:gpu:perf}
\begin{figure}[htbp]
\centering
\includegraphics[width=\textwidth]{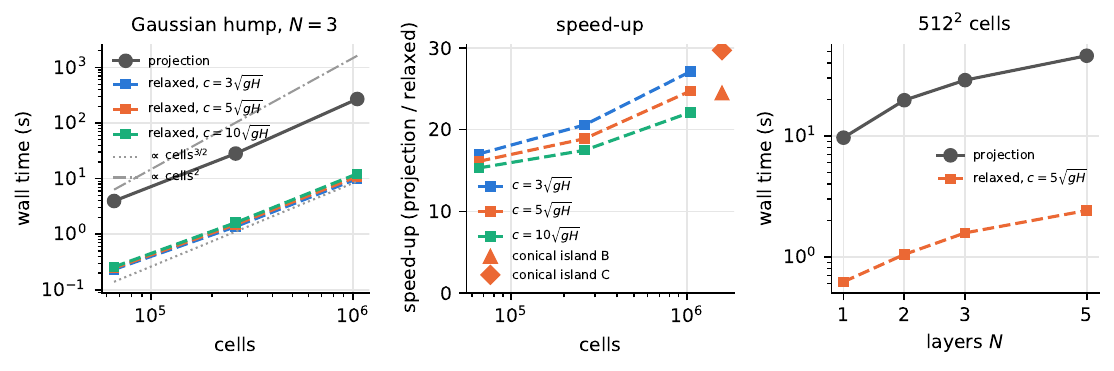}
\caption{Wall times on an A100 GPU. Left: Gaussian hump with $N=3$ versus the number of cells,
projection method and relaxed scheme for three values of $c$ (reference slopes
$(\text{cells})^{3/2}$ and $(\text{cells})^2$). Center: speed-up of the relaxed scheme (Gaussian
hump and conical island, cases B and C). Right: Gaussian hump on the $512^2$ mesh versus the
number of layers.}
\label{fig:gpu:times}
\end{figure}

Let us finally compare the wall times of both GPU codes for the same problems and the same meshes
(Figure~\ref{fig:gpu:times} and Table~\ref{tab:gpu:hump}).
On the $1024^2$ mesh with $N=3$, a time step of the relaxed
scheme takes $5.0$\,ms for $c=5\sqrt{gH}$, that is, about $2.1\cdot10^{8}$ cell updates per second. A time step of
the projection method takes $126$\,ms, because each of its two stages requires the
assembly of the block system and about $70$ iterations of the preconditioned conjugate gradient
solver.

For a fixed final time, the cost of the explicit relaxed
scheme grows like $(\text{cells})^{3/2}$, whereas the cost of the projection method grows faster,
because the number of iterations of the solver also grows with the resolution (from $31$ to $72$
per stage between $256^2$ and $1024^2$ cells): the condition number of the
non-hydrostatic operator grows like $(H/\Delta x)^2$. The speed-up of the relaxed scheme therefore
grows with the resolution, from $15$--$17$ on the $256^2$ mesh to $23$--$28$ on the $1024^2$
mesh. It is $24$--$30$ for the conical island (Section~\ref{sec:gpu:island}) and $39$--$50$ for the semicircular shoal (Section~\ref{sec:gpu:whalin}),
although the mesh of the shoal has fewer cells than the finest mesh of the hump: the speed-up follows the
number of iterations per stage and not the number of cells. It decreases moderately with $c$,
because the number of sub-steps $n_s$ is small (one to three for $c\le10\sqrt{gH}$ with $\mathrm{CFL}_{\mathrm{NH}}=0.9$).

The cost of the relaxed scheme grows a little more slowly than the number of layers
(Figure~\ref{fig:gpu:times}, right), because part of the work does not depend on $N$, whereas the cost of
the projection method grows almost linearly with $N$, since the number
of iterations hardly changes. On the $512^2$ mesh with $c=5\sqrt{gH}$ the
relaxed scheme is between $16$ and $19$ times faster.

These speed-ups are obtained with a standard preconditioned conjugate gradient solver.
Table~\ref{tab:gpu:solvers} compares it with the block Jacobi--Chebyshev iteration for the same
system and the same tolerance (both give the same solution). The
Chebyshev iteration needs no inner products and is the faster
solver on the $256^2$ mesh, but its number of iterations grows faster with the resolution, and on
the $1024^2$ mesh it is slower than the conjugate gradient. On a single GPU the
conjugate gradient is therefore the better choice; the Chebyshev iteration may still be of
interest on several GPUs, where the global reductions of the inner products are expensive.
A multigrid solver would reduce the number of iterations further, and therefore the speed-up, at
the price of a more complex implementation.

\begin{table}[htbp]
\centering
\small
\begin{tabular}{rrrrrrr}
\toprule
& \multicolumn{2}{c}{conjugate gradient} & \multicolumn{2}{c}{Chebyshev} & \multicolumn{2}{c}{speed-up of the relaxed scheme}\\
mesh & time (s) & iterations & time (s) & iterations & vs CG & vs Chebyshev\\
\midrule
$256^2$ & 3.94 & 31 & 2.96 & 68 & 16 & 12\\
$512^2$ & 28.3 & 49 & 56.6 & 169 & 18 & 37\\
$1024^2$ & 277 & 72 & 977 & 351 & 26 & 90\\
\bottomrule
\end{tabular}

\caption{Gaussian hump, $N=3$: wall times on an A100 GPU and iterations per stage of the projection
method with the block Jacobi preconditioned conjugate gradient and with the block Jacobi--Chebyshev
iteration (same system, tolerance $10^{-8}$), and speed-up of the relaxed scheme ($c=5\sqrt{gH}$)
against each of them. All the runs in the same job.}
\label{tab:gpu:solvers}
\end{table}

\section{Conclusions}\label{sec:conclusions}
We have proposed a hyperbolic relaxation of the non-hydrostatic multilayer system LDNH$_0$. It
is based on an equivalent form of the incompressibility constraints, centered at the interfaces
between the layers and adjoint to the non-hydrostatic pressure force. The relaxed system has an exact energy
balance, and the wave speeds of the system linearized around the lake at rest are bounded
independently of the number of layers.

The proposed scheme uses a staggered mesh and needs no numerical viscosity in the
non-hydrostatic part. The non-hydrostatic part is solved in one step, explicit in the horizontal
velocities and implicit only in the vertical coupling of each water column, with one tridiagonal
system of size $N$ per column. The space discretization conserves the energy exactly, and the time step does it up to a term of order $\tau^2$. For uniform states, the step
is stable when $c\tau\le\Delta x$ for any number of layers. The scheme is second-order accurate, well-balanced and robust
at wet-dry fronts. With well-prepared initial data and $3\sqrt{gH}\le c\le10\sqrt{gH}$, its
solutions are very close to those of LDNH$_0$ in our tests.
Nevertheless, avoiding the linear resonance of the fast branch of the relaxation with the solitary waves requires $c\gtrsim1.4\,N\sqrt{gH}$ (Remark~\ref{rem:resonance}), which is larger than $3\sqrt{gH}$ for $N\ge3$, although the resulting tails are not visible in our tests. The exact solitary waves show that LDNH$_0$ with a
few layers reproduces the solitary waves of the Euler equations very accurately, including their
vertical structure.

In two dimensions and on GPU, the relaxed scheme is between $15$ and $50$ times faster than the projection method for the same problem and mesh: from $15$--$17$ times faster on $256^2$ cells to $23$--$28$ times faster on $1024^2$ cells, $24$--$30$ times faster for the conical island experiment and $39$--$50$ times faster for the periodic waves over a semicircular shoal. In
this last test, two or three layers clearly improve the higher harmonics with respect to one layer. The speed-up grows with the resolution, because the number of iterations of the
elliptic solver of the projection method grows with $H/\Delta x$, while the relaxed scheme is
explicit.

Future work includes non-uniform layers, a breaking mechanism in the model, and applications of the
two-dimensional GPU code to large-scale problems.

\section*{Data availability}
The codes and the scripts that reproduce the numerical tests are available from the
corresponding author on reasonable request.

\section*{CRediT authorship contribution statement}
\textbf{Cipriano Escalante:} Conceptualization, Methodology, Software, Validation, Formal analysis, Investigation,
Visualization, Writing -- original draft. \textbf{Manuel J. Castro:} Conceptualization,
Formal analysis, Writing -- review \& editing, Funding acquisition. \textbf{Tom\'as Morales de Luna:}
Conceptualization, Formal analysis, Writing -- review \& editing.

\section*{Declaration of competing interest}
The authors declare that they have no known competing financial interests or personal
relationships that could have appeared to influence the work reported in this paper.

\section*{Declaration of generative AI and AI-assisted technologies\\in the writing process}
During the preparation of this work the authors used Claude (Anthropic) to assist in the
implementation of the numerical codes, the preparation of the figures and the editing of the
manuscript. After using this tool, the authors reviewed and edited the content as needed and take
full responsibility for the content of the publication.

\section*{Acknowledgements}
{\emergencystretch=2em
This work is supported by projects PCI2024-\allowbreak 155061-\allowbreak 2, PDC2025-\allowbreak 165273-\allowbreak I00 and PID2025-\allowbreak 170651NB-\allowbreak C21
funded by MCIN/\allowbreak AEI/\allowbreak 10.13039/\allowbreak
501100011033/ and ERDF A way of making Europe and by Grant PPRO-FQM216-G-2023
funded by Consejer\'ia de Universidad, Investigaci\'on e Innovaci\'on and by
ERDF Andalusia Program 2021-2027. Manuel J. Castro thanks the
IMUS-Mar\'ia de Maeztu grant CEX2024-001517-M --- Apoyo a Unidades de
Excelencia Mar\'ia de Maeztu for supporting this research, funded by
MICIU/\allowbreak AEI/\allowbreak 10.13039/\allowbreak 501100011033. 
\par}

\appendix

\section{Discretization of the hydrostatic step and wet-dry fronts}\label{app:hydro}\label{app:wetdry}
In this Appendix, the spatial discretization of the hydrostatic part \eqref{eq:H} of
Section~\ref{sec:split:time} is described in detail, together with the treatment of the wet-dry fronts in both parts of the time step. Let us first recall the location of the
unknowns, which is sketched in Figure~\ref{fig:stencil}. The depth $h_i$ and the momenta of the
layers
\[
q_{\alpha,i}=(h_\alpha u_\alpha)_i,\qquad \alpha=1,\dots,N,
\]
are cell averages: they are located at the center $x_i$ of the cell $I_i$ and, in the vertical
direction, at the middle of the layers. The vertical velocities $w_{\alpha,i+1/2}$,
$\alpha=1,\dots,N$, are located at the edges $x_{i+1/2}$ and at the middle of the layers, while
the pressures $p_{\aph,i+1/2}$, $\alpha=0,\dots,N-1$, are located at the edges and at the
interfaces between the layers, with $p_{N+1/2}=0$ at the free surface. The values at the edges
are point values, and it is convenient to associate them with the dual cells
\[
J_{i+1/2}=[x_i,x_{i+1}] .
\]
Two different discretizations are used: a path-conservative finite volume method for the depth
and the momenta, on the cells $I_i$, and an upwind finite difference scheme for the transport of
the vertical velocities and the pressures, on the dual cells $J_{i+1/2}$.

\paragraph{Depth and momenta} Let us denote by
\[
W=(h,q_1,\dots,q_N)^T
\]
the vector of the unknowns located at the cell centers. The equations for $h$ and $q_\alpha$ in
\eqref{eq:H} can be written in the form
\begin{equation}\label{eq:hydroW}
\dt W+\dx F(W)+\mathcal B(W)\,\dx W=-\frac{g\,h}{N}\,\hat{\mathbf e}\,\dx z_b,
\end{equation}
where
\[
F(W)=\Bigl(\sum_{\beta=1}^Nq_\beta,\ \frac{q_1^2}{h_1}+\frac{gh^2}{2N},\ \dots,\ \frac{q_N^2}{h_N}+\frac{gh^2}{2N}\Bigr)^{T},
\qquad
\hat{\mathbf e}=(0,1,\dots,1)^{T}.
\]
Here, the hydrostatic pressure term $gh_\alpha\dx\eta$ has been split into the conservative part
$\dx\bigl(gh^2/(2N)\bigr)$ and the bottom source term $-(gh/N)\dx z_b$. The non-conservative
products $\mathcal B(W)\dx W$ are the exchange terms: the component of $\mathcal B(W)\dx W$
corresponding to $h$ vanishes, and the one corresponding to $q_\alpha$ is
\[
\tilde u_\aph\Gamma_\aph-\tilde u_\amh\Gamma_\amh .
\]
Since $h_\beta\ub=\frac1N\sum_\gamma q_\gamma$, the exchange fluxes
\[
\Gamma_\aph=\sum_{\beta=\alpha+1}^N\dx q_\beta-\frac{N-\alpha}{N}\sum_{\beta=1}^N\dx q_\beta
\]
are linear combinations of the derivatives of the momenta, with coefficients that do not depend
on $W$, while $\tilde u_\aph=\frac12(u_\alpha+u_{\alpha+1})$ does. Note that the exchange terms
cancel when they are summed over the layers, since $\Gamma_{1/2}=\Gamma_{N+1/2}=0$, so that they
do not modify the total momentum.

\paragraph{Reconstruction} Second order in space is obtained by a MUSCL reconstruction with the
monotonized central (MC) limiter. For a quantity $f$ defined at the cell centers, the limited
difference in the cell $I_i$ is
\[
\delta^{\mathrm{mc}}_if=\mathrm{mc}\bigl(\delta_{i-1/2}f,\delta_{i+1/2}f\bigr),
\]
where
\[
\mathrm{mc}(a,b)=\mathrm{minmod}\Bigl(2a,\frac{a+b}2,2b\Bigr),
\]
and $\mathrm{minmod}$ returns the argument of smallest modulus if all of them have the same sign,
and zero otherwise. The reconstructed values at the edges of $I_i$ are
\[
f^-_{i+1/2}=f_i+\frac12\,\delta^{\mathrm{mc}}_if,\qquad
f^+_{i-1/2}=f_i-\frac12\,\delta^{\mathrm{mc}}_if,
\]
the superscripts $-$ and $+$ denoting the values on the left and on the right of an edge. The
reconstruction is applied to the depth, to the free surface and to the velocities $u_\alpha$,
which are desingularized as described at the end of this Appendix.
Note that $\delta^{\mathrm{mc}}_if$ vanishes unless $\delta_{i-1/2}f$ and $\delta_{i+1/2}f$ have the same sign, and that
\[
\bigl|\delta^{\mathrm{mc}}_if\bigr|\le2\min\bigl(|\delta_{i-1/2}f|,|\delta_{i+1/2}f|\bigr),
\]
so that the reconstructed values satisfy
\[
\min\bigl(f_{i-1},f_i,f_{i+1}\bigr)\le f^-_{i+1/2},\,f^+_{i-1/2}\le\max\bigl(f_{i-1},f_i,f_{i+1}\bigr).
\]
In particular, the reconstructed depths are non-negative, and no further limitation of the slopes is needed to preserve the positivity of the depth.
Some care is needed with the free surface: its
limited differences are computed from the hydrostatic jumps $\delta^*_{i\pm1/2}\eta$ of
\eqref{eq:dzHR}, so that the free surface of a dry cell is never used. At the edges, the bottom
and the momenta are then given by
\[
z^\pm_{b,i+1/2}=\eta^\pm_{i+1/2}-h^\pm_{i+1/2},\qquad
q^\pm_{\alpha,i+1/2}=\frac{h^\pm_{i+1/2}}N\,u^\pm_{\alpha,i+1/2},
\]
and the momenta are set to zero at the edges where $h^\pm_{i+1/2}<h_{\mathrm{dry}}$
(see below).

\paragraph{Path-conservative HLL scheme} Following \cite{CastroMoralesPares2017}, the
semi-discrete scheme for \eqref{eq:hydroW} reads
\[
\frac{dW_i}{dt}=-\frac1{\Delta x}\Bigl(D^-_{i+1/2}+D^+_{i-1/2}+\mathcal I_i\Bigr),
\]
where $D^-_{i+1/2}$ and $D^+_{i-1/2}$ are the fluctuations at the edges of $I_i$, and
\[
\begin{aligned}
\mathcal I_i={}&F\bigl(W^-_{i+1/2}\bigr)-F\bigl(W^+_{i-1/2}\bigr)
+\mathcal B(W_i)\bigl(W^-_{i+1/2}-W^+_{i-1/2}\bigr)\\
&+\frac{g\,h_i}{N}\bigl(z^-_{b,i+1/2}-z^+_{b,i-1/2}\bigr)\,\hat{\mathbf e}
\end{aligned}
\]
approximates, by the midpoint rule, the integral over $I_i$ of the terms of \eqref{eq:hydroW}
for the reconstructed state, which is linear in the cell.

The fluctuations are computed by means of the hydrostatic reconstruction \cite{Audusse2004}. At
the edge $x_{i+1/2}$, let
\[
b^*_{i+1/2}=\max\bigl(z^-_{b,i+1/2},z^+_{b,i+1/2}\bigr),\qquad
h^{*\pm}_{i+1/2}=\max\bigl(0,\eta^\pm_{i+1/2}-b^*_{i+1/2}\bigr),
\]
and let $W^{*\pm}_{i+1/2}$ be the states with depths $h^{*\pm}_{i+1/2}$ and the same
velocities as $W^\pm_{i+1/2}$. Then, we set
\[
\begin{aligned}
&D^-_{i+1/2}=\mathcal H^-\bigl(W^{*-}_{i+1/2},W^{*+}_{i+1/2}\bigr)
+F_{\mathrm{adv}}\bigl(W^{*-}_{i+1/2}\bigr)-F_{\mathrm{adv}}\bigl(W^-_{i+1/2}\bigr)\\
&\hphantom{D^-_{i+1/2}={}}{}+\mathcal B\bigl(W^-_{i+1/2}\bigr)\bigl(W^{*-}_{i+1/2}-W^-_{i+1/2}\bigr),\\
&D^+_{i+1/2}=\mathcal H^+\bigl(W^{*-}_{i+1/2},W^{*+}_{i+1/2}\bigr)
+F_{\mathrm{adv}}\bigl(W^+_{i+1/2}\bigr)-F_{\mathrm{adv}}\bigl(W^{*+}_{i+1/2}\bigr)\\
&\hphantom{D^+_{i+1/2}={}}{}+\mathcal B\bigl(W^+_{i+1/2}\bigr)\bigl(W^+_{i+1/2}-W^{*+}_{i+1/2}\bigr),
\end{aligned}
\]
where $F_{\mathrm{adv}}$ denotes the advective part of the flux, that is, $F$ without the
hydrostatic pressure $gh^2/(2N)$, and $\mathcal H^\pm$ are the HLL fluctuations of the system
without bottom,
\[
\mathcal H^\pm(W_L,W_R)=\frac12\Bigl((1\pm\rho_1)\,\mathcal A(W_L,W_R)\pm\rho_0\,(W_R-W_L)\Bigr),
\]
with
\[
\mathcal A(W_L,W_R)=F(W_R)-F(W_L)+\mathcal B\Bigl(\frac{W_L+W_R}2\Bigr)(W_R-W_L)
\]
and
\[
\rho_0=\frac{S_R|S_L|-S_L|S_R|}{S_R-S_L},\qquad \rho_1=\frac{|S_R|-|S_L|}{S_R-S_L}.
\]
The speeds $S_L$ and $S_R$ are estimates of the smallest and the largest characteristic speeds,
\[
\begin{aligned}
&S_L=\min_\alpha\min\Bigl(u_{\alpha,L}-\sqrt{gh_L},\,u_{\alpha,R}-\sqrt{gh_R}\Bigr),\\
&S_R=\max_\alpha\max\Bigl(u_{\alpha,L}+\sqrt{gh_L},\,u_{\alpha,R}+\sqrt{gh_R}\Bigr).
\end{aligned}
\]
They are the speeds used for the time step (Section~\ref{sec:split:time}). Without vertical shear
they reduce to $\ub\pm\sqrt{gh}$, but when the velocity varies along the vertical, the speeds of the mean
velocity can be smaller than the velocity of the fastest layer, which would then get too little
numerical viscosity.
Note that $\mathcal A$ uses the segment path between $W_L$ and $W_R$, with the midpoint rule for
the exchange terms, and that the scheme reduces to the upwind scheme when $S_L\ge0$ or
$S_R\le0$: for instance, if $S_L\ge0$, then $\rho_0=0$ and $\rho_1=1$, so that
$\mathcal H^-=0$ and $\mathcal H^+=\mathcal A$.

The last terms of $D^\pm_{i+1/2}$ are the exchange terms of the segments $W^-_{i+1/2}\to W^{*-}_{i+1/2}$
and $W^{*+}_{i+1/2}\to W^+_{i+1/2}$ of the path. Along them only the depth changes, so that
$\mathcal B$ is constant and the integral is exact. On these segments the mass flux of each layer
changes by $(h^*-h)\,u_\alpha/N$, as in the advective fluxes, and the exchange terms move the
corresponding momentum between the layers. In this way, the exchange terms always use the same
mass fluxes of the layers as the transport of momentum. Without them, at a step of the bottom the
layers would exchange mass without its momentum; in the extreme case of a step that acts as a wall,
$h^*=0$, no mass leaves the cell, while the exchange terms inside the cell would still distribute
the momentum as if it did.

Two observations are in order. Firstly, the scheme preserves the lake at rest, also with
emerged areas. Indeed, if $q_\alpha=0$ and $\eta$ is constant in the wet cells, then the
hydrostatic jumps of $\eta$ vanish, so that the reconstructed free surface is constant and
$W^{*-}_{i+1/2}=W^{*+}_{i+1/2}$ at every edge. Therefore, $\mathcal H^\pm=0$ and, since
$F_{\mathrm{adv}}$ and the exchange terms vanish when the momenta vanish, $D^\pm_{i+1/2}=0$. On the other hand, since
the reconstruction is linear, $h_i=\frac12\bigl(h^-_{i+1/2}+h^+_{i-1/2}\bigr)$, and the momentum
components of $\mathcal I_i$ are
\[
\frac g{2N}\Bigl(\bigl(h^-_{i+1/2}\bigr)^2-\bigl(h^+_{i-1/2}\bigr)^2\Bigr)
+\frac{gh_i}N\bigl(z^-_{b,i+1/2}-z^+_{b,i-1/2}\bigr)
=\frac{gh_i}N\bigl(\eta^-_{i+1/2}-\eta^+_{i-1/2}\bigr)=0 .
\]
Secondly, away from the wet-dry fronts and when the bottom is continuous, $W^{*\pm}=W^\pm$ up to
second-order terms, and the scheme is the second-order path-conservative HLL scheme of
\cite{CastroMoralesPares2017}.

\paragraph{Vertical velocities and pressures} The values at the edges are advanced with the
non-conservative equations of Section~\ref{sec:scheme:H}, which we write as
\[
\begin{aligned}
&\dt w_\alpha+u_\alpha\dx w_\alpha=-\frac1{2h_\alpha}\Bigl((w_{\alpha+1}-w_\alpha)\Gamma_\aph
-(w_{\alpha-1}-w_\alpha)\Gamma_\amh\Bigr),\\
&\dt p_\aph+\hat u_\aph\dx p_\aph=-\frac1{\hs_\aph}\Bigl((p_{\alpha+1}-p_\aph)\Gb_{\alpha+1}
-(p_\alpha-p_\aph)\Gb_\alpha\Bigr),
\end{aligned}
\]
where, since the layers have the same thickness, the transport velocities of the pressures are
\[
\hat u_{1/2}=u_1,\qquad \hat u_\aph=\frac12\left(u_\alpha+u_{\alpha+1}\right)\quad\text{for }\alpha\ge1 .
\]
Both equations are of the form
\[
\dt f+a\,\dx f=s,
\]
where $f$ is located at the edges, and they are discretized by the second-order upwind scheme
\[
\frac{df_{i+1/2}}{dt}=-\frac{a^+_{i+1/2}}{\Delta x}\bigl(f^-_{i+1}-f^-_i\bigr)
-\frac{a^-_{i+1/2}}{\Delta x}\bigl(f^+_{i+1}-f^+_i\bigr)+s_{i+1/2},
\]
with $a^+=\max(a,0)$ and $a^-=\min(a,0)$. The transport velocity at the edge is the mean of
the values in the two neighboring cells,
\[
a_{i+1/2}=\{u_\alpha\}_{i+1/2}\quad\text{for }w_\alpha,\qquad
a_{i+1/2}=\{\hat u_\aph\}_{i+1/2}\quad\text{for }p_\aph,
\]
so that the upwind direction is given by the sign of this mean. The values $f^\mp_i$ at the
point $x_i$, which separates the dual cells $J_{i-1/2}$ and $J_{i+1/2}$, are given by the
MC-limited linear reconstruction on the dual cells,
\[
f^-_i=f_{i-1/2}+\frac12\,\delta^{\mathrm{mc}}_{i-1/2}f,\qquad
f^+_i=f_{i+1/2}-\frac12\,\delta^{\mathrm{mc}}_{i+1/2}f,
\]
where, for the values at the edges,
\[
\delta^{\mathrm{mc}}_{i+1/2}f=\mathrm{mc}\bigl(\delta_if,\delta_{i+1}f\bigr),\qquad
\delta_if=f_{i+1/2}-f_{i-1/2}.
\]
For instance, if $a_{i+1/2}>0$, the scheme reads
\[
\frac{df_{i+1/2}}{dt}=-\frac{a_{i+1/2}}{\Delta x}\Bigl(f_{i+1/2}-f_{i-1/2}
+\frac12\bigl(\delta^{\mathrm{mc}}_{i+1/2}f-\delta^{\mathrm{mc}}_{i-1/2}f\bigr)\Bigr)+s_{i+1/2},
\]
that is, the first-order upwind scheme corrected by the difference of the limited differences.
Without the limiter, that is, with centered differences, it is Fromm's scheme.

The exchange terms $s_{i+1/2}$ only involve the values of $w$ and $p$ at the same edge. They
are the centered approximations written above: the vertical velocity at the interface $\aph$
is $\frac12(w_\alpha+w_{\alpha+1})$, and the pressure in the layer $\alpha$ is
$p_\alpha=\frac12(p_\amh+p_\aph)$. The thicknesses are $h_{\alpha,i+1/2}$ and
$\hs_{\aph,i+1/2}$, and the exchange fluxes are computed at the edges by compact differences,
\[
\begin{aligned}
&\Gamma_{\aph,i+1/2}=\frac1{\Delta x}\Bigl(\sum_{\beta=\alpha+1}^N\delta_{i+1/2}q_\beta
-\frac{N-\alpha}N\sum_{\beta=1}^N\delta_{i+1/2}q_\beta\Bigr),\\
&\Gb_{\alpha,i+1/2}=\frac12\bigl(\Gamma_{\amh,i+1/2}+\Gamma_{\aph,i+1/2}\bigr).
\end{aligned}
\]
At the dry edges, $w=p=0$.

Some comments on the use of the non-conservative form are in order.
For smooth solutions, these equations are equivalent to the equations for $h_\alpha w_\alpha$ and $\hs_\aph p_\aph$ in \eqref{eq:H}.
Note that the latter are not in conservation form either, since their exchange terms are non-conservative products, as those of the momentum equations.
They could be discretized as the momenta, with the exchange terms computed from the same mass fluxes as the transport terms, but some care is needed, since $w$ and $p$ are located at the edges.
Firstly, the thicknesses at the edges are not unknowns of the scheme, but the means $h_{\alpha,i+1/2}$ of the values in the two neighboring cells.
A conservative scheme on the dual mesh would then need discrete mass equations for these means, with mass fluxes at the points $x_i$ and exchange fluxes at the edges compatible with the HLL fluxes of the cells.
Otherwise, a value of $w_\alpha$ or $p_\aph$ that is constant in space would be modified by the transport.
Secondly, $w_\alpha$ and $p_\aph$ would be recovered by dividing by these thicknesses, which may be arbitrarily small near the wet-dry fronts.
The scheme above avoids both difficulties: its transport terms vanish when $f$ is constant in space, and they do not involve the thicknesses.
Moreover, the limited differences $\delta^{\mathrm{mc}}_{i+1/2}f$ satisfy the same bound as in the cells, so that each stage of the Runge--Kutta method is, up to the exchange terms, a convex combination of $f_{i+1/2}$ and of the value at the upwind neighboring edge, provided that
\[
|a_{i+1/2}|\,\Delta t\le\Delta x .
\]
This condition holds with a large margin with the time step of Section~\ref{sec:split:time}, since $|a_{i+1/2}|\le\max_{i,\alpha}|u_{\alpha,i}|$ and $\mathrm{CFL}=0.45$.
Therefore, the transport does not create new extrema of $w$ and $p$.

\paragraph{Time integration} The semi-discrete equations for $W$, $\mathbf w$ and $\mathbf p$
are integrated together by the second-order SSP Runge--Kutta method: denoting by $Q$ the vector
of all the unknowns and by $\mathcal R(Q)$ the right-hand side, each half step of
\eqref{eq:strang} reads
\[
Q^{(1)}=Q+\frac{\Delta t}2\,\mathcal R(Q),\qquad
\mathcal S_{\mathrm H}\Bigl(\frac{\Delta t}2\Bigr)Q
=\frac12\,Q+\frac12\Bigl(Q^{(1)}+\frac{\Delta t}2\,\mathcal R\bigl(Q^{(1)}\bigr)\Bigr).
\]
The velocities that transport $w$ and $p$ are those of the current stage.

\paragraph{Two dimensions} On Cartesian meshes, the depth and the momenta
$(h_\alpha u_\alpha,h_\alpha v_\alpha)$ are cell averages, and the fluctuations are computed at
the four edges of each cell by the one-dimensional scheme in the direction normal to the edge
\cite{Castro2011CRM}. At each edge, the momenta are written in their normal and tangential
components. The tangential momenta $h_\alpha u_{\alpha,\mathrm t}$ are included in the same HLL
fluctuations as the normal ones: their flux is $h_\alpha u_{\alpha,\mathrm t}u_{\alpha,\mathrm n}$, they have
the same numerical viscosity $\rho_0$, and their exchange terms are
$\tilde u_{\mathrm t,\aph}\Gamma_\aph-\tilde u_{\mathrm t,\amh}\Gamma_\amh$, with the same exchange fluxes
$\Gamma_\aph$, computed from the normal momenta. The vertical velocities and the pressures are
located at the vertices (Figure~\ref{fig:stencil}(b)). Their transport terms are the sum of the
one-dimensional upwind terms above in the directions $x$ and $y$, computed along the lines of
vertices, with the transport velocities given by the mean of the values in the four cells of
the vertex. The exchange fluxes at the vertices are computed with the same compact differences
as the constraint, with the four cells of the vertex (Section~\ref{sec:gpu:disc}).

Let us finally describe the treatment of the wet-dry fronts, first in the hydrostatic step and then in the non-hydrostatic one.

\paragraph{Wet-dry fronts in the hydrostatic step} We follow the single-layer case. The velocities are computed with
the desingularization \cite{KurganovPetrova2007}
\[
u_\alpha=\frac{\sqrt2\,h\,(Nq_\alpha)}{\sqrt{h^4+\max\bigl(h^4,h_{\mathrm{dry}}^4\bigr)}},
\]
which reduces to $u_\alpha=q_\alpha/h_\alpha$ when $h\ge h_{\mathrm{dry}}$. The MUSCL reconstruction
acts on $h$, on the hydrostatic jumps of $\eta$ and on the velocities, with
$q_\alpha=h_\alpha u_\alpha$ at the edges (zero if the edge is dry), and the momenta and the edge
values are set to zero in dry cells and edges ($h<h_{\mathrm{dry}}$, with $h_{\mathrm{dry}}=10^{-6}$
in units of the reference depth). The lake at rest is kept also with emerged areas
if the initial depth is computed from the cell averages of $\eta$ and of the bottom.

\paragraph{Wet-dry fronts in the non-hydrostatic step} Near the shoreline, without \eqref{eq:limB}, the value of
$\lambda_{\max}$ in Theorem~\ref{thm:stab} can be very large, for two reasons. First, the
geometric terms of the constraint couple $u$ and $p$ with frequencies of order $c\,N|\dx z|/h$.
With the plain jumps of $z$ these frequencies have no bound in shallow water over a sloping
bottom, since $\delta z_b$ does not vanish when $h\to0$. The hydrostatic reconstruction
\eqref{eq:dzHR} does not see the step of the bottom from the dry side, so that
$|\delta_{i+1/2}z_\aph|\le\max(h_i,h_{i+1})$, but $N|\delta z|/h$ is still of order $N$ at the
front. Second, at an edge between a wet and an almost dry cell, the divergence terms, of size
$h_{i+1/2}/\Delta x$, act on the small inertia $h_i$ of the almost dry cell. The limited
coefficients \eqref{eq:limB} are designed to remove both effects, and only act in these regions; the
non-hydrostatic pressure is not switched off. At the dry edges ($h_{i+1/2}<h_{\mathrm{dry}}$) we set
$w=p=0$. In addition, the velocities of all the layers are replaced by their mean in the cells
with $h<2\Delta x$, after each part of the step. This collapse to a single layer conserves the
momentum and does not increase the energy. Note that there the depth is less than two cells, so
that the mesh cannot resolve the vertical structure of the flow anyway.

\paragraph{Wet-dry fronts in two dimensions} At the vertex $v$, let $h_v$ be the mean of the depths $h_j$ of its
four cells $j\in\mathcal N_v$, so that $h_{\alpha,v}=h_v/N$ (Section~\ref{sec:gpu:disc}). The thickness in the
coefficients of the constraint is
\[
h^{\mathrm{lim}}_v=\frac1N\min\Bigl(h_v,(1+\varepsilon_{\mathrm h})\min_{j\in\mathcal N_v}h_j\Bigr),
\]
and the jumps of each interface across the vertex,
$(\Delta x\,\partial_xz_\aph,\Delta y\,\partial_yz_\aph)_v$, are limited as a vector: they are
multiplied by $\min\bigl(1,h^{\mathrm{lim}}_v/|(\Delta x\,\partial_xz_\aph,\Delta y\,\partial_yz_\aph)_v|\bigr)$.
Limiting each component separately is not enough, since both components may then take their
largest value at the same time. As in one dimension, the velocities of all the layers are
replaced by their mean in the cells with $h<2\Delta x$, and the pressure is not switched off.
The
projection method of Section~\ref{sec:gpu} uses the same treatment.

\section{The projection method and the well-prepared initial data}\label{app:proj}\label{app:init}
The reference solutions of LDNH$_0$ are computed by a projection method on the same staggered
mesh as the relaxed scheme. Let us denote by $\mathcal R(Q)$ the right-hand side of the
discretized hydrostatic part of LDNH$_0$ (finite volumes for $h$ and $q_\alpha$, and transport of $w_\alpha$;
see \ref{app:hydro}), and by $\mathcal P_\theta[R]$ the projection of a state $R$, which computes
\[
Q=R+\theta\,\mathbb M^{-1}B^T\mathbf p\qquad\text{such that}\qquad\mathcal K(Q)=0,
\]
where $B$ and $\mathbb M$ are the matrices of Section~\ref{sec:scheme:stag}. Since $\mathcal K$
is linear in the velocities, this leads to the linear system
\[
B\,\mathbb M^{-1}B^T\boldsymbol\psi=-\mathcal K(R),\qquad \boldsymbol\psi=\theta\,\mathbf p,
\]
which is block tridiagonal in one dimension and has a nine-point stencil in two dimensions, with
$N\times N$ tridiagonal blocks in both cases.

In time, we use the two-stage, second-order, L-stable diagonally implicit Runge--Kutta scheme
LSDIRK2 \cite{BoscheriPareschi2021}, with
\[
\gamma=1-\frac{1}{\sqrt2},\qquad \zeta=\frac{1}{2\gamma},
\]
in which the pressure is implicit and the hydrostatic part is explicit. In two dimensions, the
scheme reads
\[
Q^{(1)}=\mathcal P_{\gamma\Delta t}\bigl[Q^n+\gamma\Delta t\,\mathcal R(Q^n)\bigr],\qquad
k_1=\frac{Q^{(1)}-Q^n}{\gamma\Delta t},
\]
\[
Q^{n+1}=\mathcal P_{\gamma\Delta t}\bigl[Q^n+(1-\gamma)\Delta t\,k_1+\gamma\Delta t\,\mathcal R(Q^n+\zeta\,\Delta t\,k_1)\bigr].
\]
The pressure force and the constraint use the depth of the projected state, so that the matrix is
symmetric and positive definite, and the system is solved by the conjugate gradient method
preconditioned with the block Jacobi method (see Section~\ref{sec:gpu}). The scheme is second-order accurate (Tables~\ref{tab:order:exact} and~\ref{tab:gpu:order}).

The solver of the projection method is also used to compute the well-prepared initial data of the relaxed scheme (Section~\ref{sec:scheme:init}).
We first project the initial state, so that $\mathcal K=0$, and then we set $p(0)$ equal to the
pressure of LDNH$_0$, for which $\dt\mathcal K=0$. Let us write the discrete constraint as
\[
\mathcal K=B(h)\,U,\qquad U=(\mathbf u,\mathbf w),
\]
where we stress the dependence of $B$ on the depth, and the pressure force as
$\mathbb M^{-1}B^T\mathbf p$, with the matrices of Section~\ref{sec:scheme:stag}. Let us also
denote by $\mathcal R_U$ and $\mathcal R_h$ the right-hand sides of the discretized hydrostatic
part for $U$ and $h$, so that, for LDNH$_0$,
\[
\dt h=\mathcal R_h(Q),\qquad \dt U=\mathcal R_U(Q)+\mathbb M^{-1}B^T\mathbf p .
\]
Then, the condition $\dt\mathcal K=0$ at $t=0$ reads
\[
B\Bigl(\mathcal R_U(Q^0)+\mathbb M^{-1}B^T\mathbf p(0)\Bigr)+\bigl(\partial_hB\;\mathcal R_h(Q^0)\bigr)U^0=0 .
\]
This pressure is computed with the solver of the projection method: one explicit step of the
hydrostatic part, $Q^0+\Delta t\,\mathcal R(Q^0)$, is projected with the weight $\Delta t$, that is,
$\mathbf p$ is such that $h^*=h^0+\Delta t\,\mathcal R_h(Q^0)$ and
$U^*=U^0+\Delta t\bigl(\mathcal R_U(Q^0)+\mathbb M^{-1}B^T\mathbf p\bigr)$ satisfy $B(h^*)U^*=0$. Since
$B(h^0)U^0=0$ after the first projection,
\[
B(h^*)U^*=\Delta t\Bigl[B\Bigl(\mathcal R_U(Q^0)+\mathbb M^{-1}B^T\mathbf p\Bigr)
+\bigl(\partial_hB\;\mathcal R_h(Q^0)\bigr)U^0\Bigr]+O(\Delta t^2),
\]
so that this pressure satisfies the condition above up to $O(\Delta t)$: it is the pressure of
LDNH$_0$ at $t=0$ up to $O(\Delta t)$. The auxiliary step is discarded: the first time step starts
from the projected state and from this pressure, and it is already a step of the relaxed scheme.

\section{Computation of the exact solitary waves}\label{app:tw}
The solitary waves of LDNH$_0$ and of the relaxed system are traveling waves $f(x-Ct)$. They are
computed by a spectral Newton method on the periodic domain $[-40,40]$ with $256$ Fourier modes. The
unknowns are the values of $h$, $u_\alpha$, $w_\alpha$ and $p_\aph$ on half of the domain (since $h$, $u$
and $p$ are even and $w$ is odd) and the speed $C$, with the amplitude condition $h(0)=H+A$ and the
conditions $u_\alpha=p_\aph=0$ far from the wave, and a continuation in $A$. For $N=1$ the method
recovers the explicit solitary wave of LDNH$_0$ with an error of $2\cdot10^{-13}$. The solitary waves of
the Euler equations are computed with the method of \cite{ClamondDutykh2013}, which we have checked
against the value of the Froude number given there, and the layer means of the velocity are obtained
from the conformal map.

\end{document}